\documentclass[a4paper,fleqn]{cas-sc}

\usepackage[authoryear,longnamesfirst]{natbib}
\usepackage{mathtools}
\usepackage{stmaryrd}
\usepackage{placeins} 
\usepackage{algorithm}
\usepackage{algpseudocode} 
\usepackage{graphicx}
\usepackage{tikz}
\usepackage{pgfplots}
\usepackage{float}
\usepackage{subcaption}
\def\tsc#1{\csdef{#1}{\textsc{\lowercase{#1}}\xspace}}
\tsc{WGM}
\tsc{QE}

\newcommand{\jump}{\par\bigskip}
\newtheorem{definition}{Definition}

\newdefinition{rmk}{Remark}

\begin{document}
\let\WriteBookmarks\relax
\def\floatpagepagefraction{1}
\def\textpagefraction{.001}

\shorttitle{Neural-manifold model order reduction for nonlinear variational inequalities}

\shortauthors{S. Ezzehi et~al.}

\title [mode = title]{Neural-Manifold Primal–Dual Model Order Reduction and Hyper-Reduction for Nonlinear Parametric Variational Inequalities}  



%

\author[1,2]{Sofiane Ezzehi}

\cormark[1]


\ead{sofiane.ezzehi@enpc.fr}



\affiliation[1]{organization={CERMICS, CNRS, ENPC, Institut Polytechnique de Paris},
            addressline={8 Av. Blaise Pascal}, 
            city={Champs-sur-Marne},
            postcode={77420}, 
            country={France}}

\author[1]{Virginie Ehrlacher}


\author[2]{Guillaume Enchéry}




\affiliation[2]{organization={IFP Energies nouvelles, Department of Mathematics},
            addressline={1 et 4 avenue de Bois-Préau}, 
            city={Rueil-Malmaison},
            postcode={92852}, 
            country={France}}

\author[2]{Thibault Faney}


\cortext[1]{Corresponding author}



\begin{abstract}
    We propose a model order reduction framework for nonlinear parametrized variational inequalities arising in computational mechanics. The high-dimensional model is written in mixed primal-dual form with projection-based complementarity conditions, leading to nonlinear nonsmooth algebraic systems solved by a semi-smooth Newton method. On this basis, reduced models are constructed by proper orthogonal decomposition (POD) of both primal and dual solution snapshots.    
    
    To address cases where low-dimensional linear spaces provide limited approximation efficiency, we introduce a neural-network-augmented reduced model. Two feedforward networks learn corrections in the truncated POD coordinates of the primal and dual variables from solutions of enlarged Galerkin ROMs. The learned reconstructions enter the coupled equilibrium and constraint equations directly, defining a nonlinear reduced system solved by approximate semi-smooth Newton iterations. The online cost associated with high-dimensional residual evaluations is reduced through hyper-reduction, using a sparse cubature approach based on greedy nonnegative least squares.
    
    The proposed methodology is assessed on two nonlinear variational inequalities with distinct sources of nonlinearity: a two-dimensional obstacle problem with a cubic nonlinearity in the state equation, and a three-dimensional frictional contact problem in which the Coulomb law induces a nonlinear projection in the constraint equation. Numerical results compare the high-dimensional model, the linear reduced model, the neural-network-augmented reduced model, and their hyper-reduced variants, demonstrating accurate approximations with substantial reductions in online computational cost. The performance study also evaluates feasibility over the complete constraint set and measures how errors in the reconstructed obstacle and contact forces affect the primal response.
\end{abstract}


    

\maketitle


\section{Introduction}\label{sec:introduction}
    The efficient numerical approximation of parametrized nonlinear variational inequalities (VIs) remains a challenging task in computational mechanics. Such formulations typically arise from variational problems involving additional constraints on the admissible set. In the standard setting, the solution is sought in a closed convex subset of a Hilbert space \cite{Grossmann2007}. VIs appear in a wide range of mechanical applications, including obstacle problems
    \cite{obstaclePenalty}, unilateral contact \cite{Wohlmuth2011}, and frictional contact \cite{Beaude2023}.
    
    Constraints in variational inequalities can be handled in several ways. One possibility is to replace the constrained problem by a relaxed or regularized formulation, for instance through penalty methods \cite{Grossmann2007,obstaclePenalty}. Another is to introduce dual variables associated with the constraints, leading to mixed, saddle-point, or primal--dual complementarity formulations \cite{MOR_Salomon,Charbel_VI}.  Recent latent-variable proximal formulations provide another approach to enforcing pointwise inequality constraints \cite{Dokken2025LVPP}. 
 In computational contact mechanics, these ideas are closely related to Lagrange-multiplier, augmented-Lagrangian \cite{AugmentedLagrangian}, and Nitsche-type weak-enforcement formulations \cite{Beaude2023,Idrissa_Nitsche,Nitsche}, which are widely used for contact and interface constraints \cite{Wohlmuth2011}.
    
    In the present work, we retain the primal-dual projected complementarity formulation. After discretization, the projection operators associated with the constraints yield nonsmooth algebraic residuals. Moreover, the problems considered here contain nonlinearities either in the state equation or in the constraint law, so that the resulting systems do not fall within the quadratic-programming framework used in several reduced-basis approaches for linear VIs \cite{MOR_Salomon,Charbel_VI}.
    
    At the algebraic level, related nonlinear complementarity systems are commonly solved using Uzawa-type projected-gradient iterations \cite{Giulia}, interior-point methods \cite{PotraYe1996}, smoothing and regularization techniques \cite{HaddouMaheux2014,Sun1999}, or semi-smooth Newton methods \cite{PDAS-SSN,SSN_for_VI}. In this work, we adopt the latter approach and solve the projected primal-dual systems by a semi-smooth Newton (SSN) method written in primal-dual form. This choice is natural for projected complementarity formulations, since the projection residual directly defines the active and inactive constraint sets. In our formulation, the linear systems solved at each SSN iteration are defined by generalized Jacobians obtained from slant derivatives of the nonsmooth projection operators \cite{PDAS_CNQ,PDAS_Q}.
    
    In parametrized settings, VI discretizations give rise to high-fidelity systems whose repeated solution over many parameter values is required in applications such as parametric studies, design, optimization, and uncertainty quantification. This many-query context motivates model order reduction (MOR), which replaces repeated solves of a high-dimensional model (HDM) by a low-dimensional reduced-order model (ROM) while preserving the relevant solution features or quantities of interest. For standard parametrized PDEs, MOR techniques have been extensively studied and applied in solid mechanics \cite{MOR_exple_solid_mechanics}, structural dynamics \cite{ECSW} and wave propagation \cite{MOR_exple_wave_propagation}. Parametrized VIs inherit many of these ideas, but introduce additional difficulties related to inequality constraints, primal-dual coupling, reduced stability, and nonsmooth complementarity or projection operators \cite{MOR_Salomon,Charbel_VI,Idrissa_Nitsche}.
    
    A classical MOR approach relies on low-dimensional \textit{linear} approximations of the solution manifold, for instance via proper orthogonal decomposition (POD) \cite{POD} or reduced basis (RB) methods \cite{RB}. However, linear reduced spaces may become inefficient when the solution manifold is poorly approximated by low-dimensional linear spaces, notably when the Kolmogorov $n$-width decays slowly \cite{Kolmogorov_decay}. This issue is particularly relevant in constrained mechanical problems, where moving contact regions, changing active sets or stick-slip transitions may induce strongly nonlinear solution manifolds. In such situations, nonlinear MOR strategies can provide more efficient approximations. Representative directions include local reduced bases \cite{local_RB}, quadratic manifold-based MOR techniques \cite{QM_ROM}, and learned nonlinear manifolds based on neural networks \cite{NN_ROM}.
    
    After projection onto reduced spaces, the online evaluation of reduced models may still involve operations that scale with the HDM dimension, notably in the presence of nonaffine parameter dependence, nonlinear residual terms, or nonsmooth projection operators. Hyper-reduction methods address this bottleneck by approximating full-dimensional contributions using only a small number of sampled entries, elements, or quadrature points. A first class of methods relies on interpolation of nonlinear or nonaffine terms, as in EIM and DEIM \cite{EIM,DEIM}, with QR-based variants such as Q-DEIM \cite{Q-DEIM}. A second class consists of cubature or sampling-weighting approaches, including greedy-NNLS cubature \cite{greedy_cubature}, the empirical cubature method (ECM) \cite{ECM}, and the energy-conserving sampling and weighting method (ECSW) \cite{ECSW}. 
    
    While projection-based MOR, nonlinear reduced representations, and hyper-reduction are now standard ingredients of efficient reduced models for nonlinear parametrized PDEs, their extension to variational inequalities raises additional difficulties: the reduced model must account not only for nonlinear and nonaffine operators, but also for primal-dual coupling, inequality constraints, and the nonsmooth structure induced by complementarity or projection operators.
    
    For parametrized VIs with linear state equations and linear constraints, several MOR strategies have been proposed. \cite{MOR_Salomon} construct separate reduced spaces for the primal variable and the Lagrange multiplier, enforcing dual feasibility through nonnegative combinations of multiplier snapshots. \cite{Charbel_VI} follow a related strategy but use nonnegative matrix factorization to obtain a more compact dual representation. In both cases, the linear structure of the problem leads to reduced formulations that can be handled within a quadratic-programming framework. More recently, \cite{Idrissa} developed a stable MOR strategy for linear VIs with parametric constraints, including guarantees on the inf-sup stability of the reduced formulation. Hyper-reduction has also been considered in related VI settings, including penalty-based formulations \cite{Bader2016,Balajewicz2017} and frictionless contact problems in mixed form \cite{Fauque2018}.
    
    Beyond the linear-VI setting, MOR strategies for variational inequalities with nonlinear state equations or nonlinear constraint laws remain comparatively scarce. In \cite{Amina}, the authors consider VIs with nonlinear constraints and solve the resulting problems by a Kačanov fixed-point method. The reduced spaces are constructed by POD for the primal variable and by a cone-projected greedy algorithm for the dual variable, in order to preserve non-negativity of the multipliers. For frictional contact problems, \cite{Idrissa_Nitsche} combine a reduced-basis approximation of a Nitsche formulation with empirical interpolation to treat the nonlinear contact and friction contributions.
    
    A complementary line of work concerns nonlinear reduced representations in the context of VIs. In \cite{Giulia}, the authors study a time-dependent parametrized VI arising from an agent-based model of crowd motion with obstacles. Although the underlying VI is linear, the work is particularly relevant here because it incorporates a nonlinear ROM inspired by \cite{NN_ROM}. This nonlinear reconstruction is, however, used as a post-processing step to recover the high-fidelity solution from a low-dimensional linear ROM, rather than being embedded directly into the reduced nonlinear solve.
    \jump
    In this work, we propose a model order reduction framework for parametrized variational inequalities with nonlinear state equations and/or nonlinear constraint laws. Both the neural-network augmentation and the hyper-reduction are embedded directly within the semi-smooth Newton iterations. In particular, the learned correction enters the coupled equilibrium and constraint equations, rather than being applied only after the solve as in post-processing reconstruction strategies such as that of \cite{Giulia}.

    The main contributions of this work can be summarized as follows:
    \begin{itemize}
    \item We formulate a semi-smooth-Newton-based primal-dual MOR framework for nonlinear parametrized VIs written in projection form. Starting from POD reduced spaces for both the primal and dual variables, two neural networks approximate the truncated POD components, in the spirit of \cite{NN_ROM}. The equilibrium and projected complementarity residuals are tested against the corresponding reconstruction tangents, and the resulting coupled system is solved by approximate semi-smooth Newton iterations.
    \item We construct the neural-network training data from Galerkin ROM solves in enlarged primal and dual POD spaces, rather than from independent projections of the HDM snapshots. The training targets therefore account for the coupling between the primal and dual variables in the reduced equations. We use an energy-norm loss for the primal network and a dual loss that measures the equilibrium action induced by multiplier-reconstruction errors.
    \item We hyper-reduce the online stage via sparse cubature rules constructed by greedy nonnegative least squares \cite{greedy_cubature}. Since the cubature procedure targets the projected residual contributions directly, its fitting targets have the lower dimension of the NN-augmented reduced equations, unlike interpolation strategies that approximate full-dimensional nonlinear fields before projection. 
    \end{itemize}
We assess the methodology on two nonlinear VI problems with distinct nonlinearities: a 2D obstacle problem with a cubic nonlinearity in the state equation and a 3D frictional contact problem governed by the Coulomb friction law. Numerical comparisons between the HDM, the linear ROM, the NN-augmented ROM, and their hyper-reduced variants demonstrate the accuracy and computational efficiency of the approach. The reduced multiplier coefficients remain unconstrained throughout the iterations, so pointwise feasibility is not imposed by construction. We therefore assess the reconstructed fields explicitly over the complete constraint set and quantify the effect of the dual reconstruction error on the primal mechanical response.

\section{Model Problems and Semi-smooth Newton Method}
\label{chap:SSN}
    We introduce in this section two parametric nonlinear variational inequality problems that will serve as running examples throughout the paper. They are chosen to illustrate two distinct sources of nonlinearity that commonly arise in variational inequalities. The first problem, presented in Section~\ref{sec:2D_obstacle_problem}, is a two-dimensional obstacle problem with a cubic nonlinearity in the state equation, inspired by the first application in \cite{obstaclePenalty}. This problem leads to a nonsmooth system of equations in which the nonlinearity appears in the primal equation. The second problem, presented in Section~\ref{sec:3D_contact_problem}, is a three-dimensional frictional contact problem inspired by the first numerical example of \cite{bubbleVEM}. We extend this setting to a parametric framework and incorporate frictional contact. In this case, the nonlinearity appears in the constraint equation, through the projection operator associated with the Coulomb friction law.
    
    We then introduce in Section~\ref{sec:unified_problem} a unified discrete formulation that covers both examples and will serve as a common reference point for the remainder of the paper. The semi-smooth Newton method used to solve the resulting nonlinear nonsmooth systems is presented in Section~\ref{sec:ssn}.
    
\subsection{Model problems and unified discrete formulation}
\subsubsection{2D Obstacle Problem with Cubic Nonlinearity}
\label{sec:2D_obstacle_problem}
  \paragraph{Problem setting}
    Let $\Omega\subset\mathbb{R}^2$ be an open bounded domain and let $\mathcal P\subset\mathbb{R}^p$ be the parameter space. For simplicity of presentation, we state the formulation with homogeneous Dirichlet boundary conditions; nonhomogeneous constant values are handled by the usual lifting and do not affect the developments below. For a given parameter $\mu\in\mathcal P$, let $\psi_\mu\in L^2(\Omega)$ denote the obstacle. The admissible set, assumed to be nonempty, is defined by
        \[
        \mathcal{G}_\mu := \left\{ v \in H^1_0(\Omega) : v \geq \psi_\mu \ \text{a.e. in } \Omega \right\}
        \]
        and encodes the non-penetration condition. The obstacle problem reads as the minimization problem:
    \begin{equation}
        \min_{v \in \mathcal{G}_\mu} \int_\Omega \left( \frac{1}{2} |\nabla v|^2 + \frac{\gamma}{4} v^4 \right) dx,
        \label{eq:obstacle_nonlinear_minimization}
    \end{equation}
    where $\gamma > 0$ controls the magnitude of the nonlinear term.
    
    Defining the closed convex cone
    \[
    \mathcal{K} := \{ w \in L^2(\Omega) : w \geq 0 \ \text{a.e. in } \Omega \},
    \]
    we introduce the Lagrange multiplier $\lambda \in \mathcal{K}$ associated with the non-penetration constraint. The optimality conditions of the minimization problem \eqref{eq:obstacle_nonlinear_minimization} can then be written in mixed form as follows: for a given parameter $\mu \in \mathcal{P}$, find $(u_\mu,\lambda_\mu) \in H^1_0(\Omega) \times \mathcal{K}$ such that, for all $(v,\eta)\in H^1_0(\Omega)\times\mathcal K$,
    \begin{equation}
        \begin{cases}
            \displaystyle
            \int_\Omega \nabla u_\mu \cdot \nabla v \, dx
            +
            \gamma \int_\Omega u_\mu^3 v \, dx
            -
            \int_\Omega \lambda_\mu v \, dx
            = 0,
            \\[2mm]
            \displaystyle
            \int_\Omega (u_\mu - \psi_\mu)(\eta-\lambda_\mu) \, dx
            \geq 0.
        \end{cases}
        \label{eq:obstacle_mix_form}
    \end{equation}
    The first equation corresponds to the nonlinear state equation, while the second one encodes the non-penetration constraint.
    
    Assuming that the admissible set $\mathcal{G}_{\mu}$ is nonempty, the
    coercivity and strict convexity of the energy functional imply existence and
    uniqueness of the primal minimizer of the minimization problem~\eqref{eq:obstacle_nonlinear_minimization}. In a general obstacle problem, the
    associated multiplier is naturally interpreted as an element of
    $H^{-1}(\Omega)$. In the present setting, we assume sufficient regularity for
    the multiplier to admit an $L^2(\Omega)$ representative, so that the
    primal-dual pairing can be written as an $L^2(\Omega)$ inner product, as in~\eqref{eq:obstacle_mix_form}. The model-reduction developments below start from
    the discrete complementarity system and do not require a separate continuous
    inf-sup analysis.
    
     We refer to \cite{Grossmann2007} for a general presentation of the theory of variational inequalities, and to \cite{obstaclePenalty} for a more specific discussion of obstacle problems with state-equation nonlinearities.
    
    \paragraph{Discretization.}
    
    We use a standard $\mathbb{P}_k$ finite element discretization, with the same finite element space for the primal and dual variables.  Let $\{\phi_i\}_{i=1}^N$ denote the corresponding basis functions. We introduce the coefficient vector $U_\mu \in \mathbb{R}^N$ of the discrete primal variable by $u^h_\mu = \sum_{i=1}^N (U_\mu)_i \phi_i$ and the dual vector of Lagrange multipliers $\Lambda_\mu \in \mathbb{R}^N$ by $\lambda^h_\mu = \sum_{i=1}^N (\Lambda_\mu)_i \phi_i$.
    
    The obstacle is represented by a vector $G_\mu \in \mathbb{R}^N$, associated with the finite element interpolant of $\psi_\mu$ in the finite element space. The complementarity conditions are then imposed coefficientwise on the nodal unknowns associated with the discrete primal and dual variables.
    
    We consider the following discrete mixed system: find
    $U_\mu \in \mathbb{R}^N$ and $\Lambda_\mu \in \mathbb{R}^N$ such that
    \begin{equation}
    \begin{cases}
        K U_\mu + \gamma M U_\mu^{\circ 3} - M \Lambda_\mu = 0,
        \\
        \Lambda_\mu \geq 0,
        \qquad
        U_\mu - G_\mu \geq 0,
        \qquad
        \Lambda_\mu^T (U_\mu - G_\mu) = 0,
    \end{cases}
    \label{eq:discrete_obstacle_mix_form_CP}
    \end{equation}
    
    where $K \in \mathbb{R}^{N \times N}$ and $M \in \mathbb{R}^{N \times N}$ are the stiffness and mass matrices, with entries $K_{ij} = \int_\Omega \nabla \phi_i \cdot \nabla \phi_j \, dx$ and $M_{ij} = \int_\Omega \phi_i \phi_j \, dx$. The notation $U_\mu^{\circ 3}$ denotes the componentwise cube of $U_\mu$, namely $(U_\mu^{\circ 3})_i = (U_\mu)_i^3$. In the discrete state equation, the nonlinear term is evaluated by nodal interpolation, so that the vector with entries $\left(\int_\Omega (u_\mu^h)^3 \phi_i \, dx\right)_{i=1}^N$ is approximated by $M U_\mu^{\circ 3}$.

    \begin{rmk}[Discrete interpretation of the inequality constraint]
        The discrete admissible cone used in this work is defined coefficientwise:
        \[\mathcal{K}_h = \left\{Z\in\mathbb{R}^{N}\;:\;Z_i\geq 0,\quad i=1,\ldots,N\right\}.\]
        Accordingly, the obstacle and multiplier inequalities in \eqref{eq:discrete_obstacle_mix_form_CP}
        are nodal, or more generally coefficientwise, discrete constraints. For
        finite elements of degree larger than one, nonnegativity of all degrees of
        freedom does not necessarily imply pointwise nonnegativity of the finite
        element function throughout every element. All feasibility and
        complementarity indicators reported below therefore refer to
        this discrete coefficientwise constraint.
    \end{rmk}
        
    The complementarity system given in the second line of \eqref{eq:discrete_obstacle_mix_form_CP} can equivalently be written in projection form. For any
    $\rho > 0$, one has
    \[ \Lambda_\mu - \max\left(0,\Lambda_\mu - \rho (U_\mu-G_\mu)\right) = 0,\]
    where the maximum is applied componentwise. Consequently, the discrete obstacle problem can be recast as the following system: find $U_\mu,\Lambda_\mu \in \mathbb{R}^N$ such that
    \begin{equation}
    \begin{cases}
        K U_\mu + \gamma M U_\mu^{\circ 3} - M \Lambda_\mu = 0,
        \\
        \Lambda_\mu
        -
        \max\left(0,\Lambda_\mu - \rho (U_\mu-G_\mu)\right)
        = 0.
    \end{cases}
    \label{eq:discrete_obstacle_nonsmooth_form}
    \end{equation}
    This projection form is the discrete nonlinear nonsmooth system solved by the semi-smooth Newton method and later projected in the reduced models.
    
\subsubsection{3D Contact Problem}
\label{sec:3D_contact_problem}
\paragraph{Problem setting.}

Let $\Omega = (-1,1)^3 \subset \mathbb{R}^3$ be the domain of interest, and let $\Gamma = \{0\}\times (-1,1)^2 \subset \Omega$ be a planar fracture interface at $x=0$ that splits $\Omega$ into two subdomains $\Omega^+$ and $\Omega^-$. We denote by $\gamma^\pm$ the trace operators on $\Gamma$ from each side, and define the jump operator
\[
    \llbracket v \rrbracket = \gamma^+ v - \gamma^- v,
    \qquad
    \forall v \in H^1(\Omega \setminus \overline{\Gamma})^3.
\]
Its normal and tangential components are denoted by $\llbracket v \rrbracket_n = \llbracket v \rrbracket \cdot n$ and $\llbracket v \rrbracket_\tau = \llbracket v \rrbracket - \llbracket v \rrbracket_n n$, where $n$ is a fixed unit normal on $\Gamma$. We further introduce the space
\[
    U_0 =
    \left\{
    v \in H^1(\Omega \setminus \overline{\Gamma})^3 :
    v = 0 \text{ on } \partial\Omega
    \right\}.
\]
Let $\mu_L>0$ and $\lambda_L>0$ be the Lamé parameters. For $u \in U_0$, we define the linear elastic stress tensor by
\[
    \sigma(u) = 2\mu_L \varepsilon(u) + \lambda_L \operatorname{div}(u) I,
    \qquad
    \varepsilon(u) = \frac12(\nabla u + \nabla u^T).
\]
We consider a parametric family indexed by $\mu\in\mathcal P$, in which the body force $g_\mu\in L^2(\Omega)^3$ and the Coulomb friction coefficient $F_\mu\geq0$ depend on $\mu$. 
For $u_\mu\in U_0$, the surface tractions on $\Gamma$ are defined by $ T_\mu^\pm = \left.\sigma(u_\mu^\pm)n^\pm\right|_\Gamma, $ with $n^+=n$ and $n^-=-n$. We decompose $T_\mu^+$ into its normal and tangential components $ T_{\mu,n}=T_\mu^+\cdot n$ and $T_{\mu,\tau}=T_\mu^+-T_{\mu,n}n. $

For each $\mu\in\mathcal P$, the static contact problem reads: find $u_\mu\in U_0$ such that
\[
\begin{cases}
  -\operatorname{div}\sigma(u_\mu)=g_\mu, & \text{in }\Omega\setminus\overline{\Gamma},\\
  T_\mu^++T_\mu^-=0, & \text{on }\Gamma,\\
  \llbracket u_\mu\rrbracket_n\leq0,\quad T_{\mu,n}\leq0,\quad T_{\mu,n}\llbracket u_\mu\rrbracket_n=0, & \text{on }\Gamma,\\
  |T_{\mu,\tau}|\leq-F_\mu T_{\mu,n}, & \text{on }\Gamma,\\
  T_{\mu,\tau}\cdot\llbracket u_\mu\rrbracket_\tau -F_\mu T_{\mu,n}|\llbracket u_\mu\rrbracket_\tau|=0, & \text{on }\Gamma.
\end{cases}
\]
The first two lines correspond to bulk equilibrium and the balance of surface tractions on the interface. The third line imposes nonpenetration and normal complementarity. The last two lines describe Coulomb friction: if $\llbracket u_\mu\rrbracket_\tau=0$, the interface sticks and $|T_{\mu,\tau}|\leq-F_\mu T_{\mu,n}$; otherwise, the interface slips, $|T_{\mu,\tau}|=-F_\mu T_{\mu,n}$, and $T_{\mu,\tau}$ is opposite to $\llbracket u_\mu\rrbracket_\tau$.

Following \cite{Wohlmuth2011,bubbleVEM}, for each $\mu\in\mathcal P$, we introduce a Lagrange multiplier $\lambda_\mu:\Gamma\to\mathbb R^3$ such that $\lambda_\mu=-T_\mu^+=T_\mu^-$ on $\Gamma$. Its normal and tangential components are $\lambda_{\mu,n}=\lambda_\mu\cdot n$ and $\lambda_{\mu,\tau}=\lambda_\mu-\lambda_{\mu,n}n$. For a given normal multiplier $\lambda_{\mu,n}$, we introduce the admissible set 
\[ C_{f,\mu}(\lambda_{\mu,n}) = \left\{ \eta\in H^{-1/2}(\Gamma)^3: \langle\eta,v\rangle_\Gamma \leq \langle F_\mu\lambda_{\mu,n},|v_\tau|\rangle_\Gamma, \quad \forall v\in H^{1/2}(\Gamma)^3,\ v_n\leq0 \right\}, \] 
which corresponds to the Coulomb friction law $|\lambda_{\mu,\tau}|\leq F_{\mu}\lambda_{\mu,n}$ together with $\lambda_{\mu,n} \geq 0$.

The frictional contact problem can then be written as follows: find  $(u_\mu,\lambda_\mu) \in U_0\times C_{f,\mu}(\lambda_{\mu,n})$ such that, for all $(v,\eta)\in U_0\times C_{f,\mu}(\lambda_{\mu,n})$, 

\[ 
  \begin{aligned} 
    a(u_\mu,v)+b(\lambda_\mu,v)&=l_\mu(v),\\
    b(\eta-\lambda_\mu,u_\mu)&\leq0, 
  \end{aligned} 
\] 
where 
\[ a(v,w)=\int_{\Omega\setminus\overline{\Gamma}} \sigma(v):\varepsilon(w)\,dx, \qquad b(\eta,v)=\langle\eta,\llbracket v\rrbracket\rangle_\Gamma, \qquad l_\mu(v)=\int_{\Omega}g_\mu\cdot v\,dx. \]

\paragraph{Discretization.}

We consider a tetrahedral mesh of $\Omega$ split by $\Gamma$ into two submeshes $\Omega^+$ and $\Omega^-$. The displacement is approximated with $(\mathbb{P}_2)^3$ finite elements on each side, yielding a discrete unknown $U=(U^+,U^-)\in\mathbb{R}^{2N}$. The interface multipliers are approximated with face-wise constant elements on $\Gamma$ and decomposed into normal and tangential components $\Lambda_n\in\mathbb{R}^{R}$ and $\Lambda_\tau\in\mathbb{R}^{2R}$, where $R$ is the number of fracture faces.

Following \cite{bubbleVEM,Beaude2023}, the resulting discrete system can be cast as the following nonlinear nonsmooth system: find 
$(U_\mu,\Lambda_{\mu,n},\Lambda_{\mu,\tau}) \in \mathbb{R}^{2N}\times\mathbb{R}^{R}\times\mathbb{R}^{2R}$ such that
\begin{equation}
\begin{aligned}
    K U_\mu + B_n^T \Lambda_{\mu,n} + B_\tau^T \Lambda_{\mu,\tau} - L_{\mu} &= 0, \\
    \Lambda_{\mu,n} - \max(0,\Lambda_{\mu,n}+\rho B_n U_\mu) &= 0, \\
    \Lambda_{\mu,\tau} - \Pi_{c_\mu(\Lambda_{\mu,n})} \left( \Lambda_{\mu,\tau}+\rho B_\tau U_\mu \right) &= 0,
\end{aligned}
\label{eq:non_smooth_eq_friction_modelproblem}
\end{equation}
where $U_\mu=(U_\mu^+,U_\mu^-)$ collects the displacement unknowns on the two subdomains. The matrices in \eqref{eq:non_smooth_eq_friction_modelproblem} are defined by
\[
\begin{aligned}
    \bigl(K^\pm\bigr)_{ij} &= a(\phi_j^\pm,\phi_i^\pm), & K &= \begin{pmatrix} K^+ & 0 \\ 0 & K^- \end{pmatrix} \in \mathbb{R}^{2N\times 2N}, \\
    \bigl(B_n^\pm \bigr)_{ji} &= b(\psi_j n,\phi_i^\pm), & B_n &= (B_n^+,B_n^-) \in \mathbb{R}^{R\times 2N}, \\
    \bigl(B_\tau^\pm \bigr)_{(j,l)i} &= b(\psi_j \tau_l,\phi_i^\pm),
    & B_\tau &= (B_\tau^+,B_\tau^-) \in \mathbb{R}^{2R\times 2N}, \\
    \bigl(L_{\mu}^\pm \bigr)_i &= l_\mu(\phi_i^\pm),
    & L_{\mu} &= (L_\mu^+,L_\mu^-) \in \mathbb{R}^{2N}.
\end{aligned}
\]

The nonlinear system \eqref{eq:non_smooth_eq_friction_modelproblem} is
solved using the semi-smooth Newton method presented in
Section~\ref{sec:ssn}. To make the tangential projection well defined at every iterate,
including trial iterates for which the normal multiplier may temporarily
become negative, we define the nonnegative friction radius
\[
c_{\mu,j}(\Lambda_n)
:=
F_{\mu}\,[\Lambda_{n,j}]_{+},
\qquad
[a]_{+}:=\max(a,0).
\]
The tangential operator is then understood as the face-wise projection onto the closed disk of radius $c_{\mu,j}(\Lambda_n)$. At any solution satisfying the normal contact conditions, one has $\Lambda_{n,j}\geq 0$, and this definition therefore coincides with the standard Coulomb disk of radius $F_\mu\Lambda_{n,j}$.

The mapping $\Pi_{c_\mu(\Lambda_n)}$ denotes the face-wise projection onto
the friction disks of radii $c_{\mu,j}(\Lambda_n)$. More precisely, for
$x=(x_1,\ldots,x_R)$ with $x_j\in\mathbb{R}^2$, we define
\[
\Pi_{c_\mu(\Lambda_n)}(x)
=
\left(
\Pi_{c_{\mu,1}(\Lambda_n)}(x_1)
\mid
\ldots
\mid
\Pi_{c_{\mu,R}(\Lambda_n)}(x_R)
\right)^T,
\]
with
\[
\Pi_{c_{\mu,j}(\Lambda_n)}(x_j)
=
\begin{cases}
x_j,
&
\text{if }
\lVert x_j\rVert_2
\leq c_{\mu,j}(\Lambda_n),
\\[1ex]
c_{\mu,j}(\Lambda_n)
\dfrac{x_j}{\lVert x_j\rVert_2},
&
\text{if }
\lVert x_j\rVert_2
>
c_{\mu,j}(\Lambda_n),
\end{cases}
\]
where $x_j \in \mathbb{R}^2$ is the tangential component of the multiplier on face $j$, and $\Lambda_{n,j}$ is the corresponding normal component.

In system \eqref{eq:non_smooth_eq_friction_modelproblem}, the nonsmoothness is localized in the two constraint equations: the normal complementarity condition through the $\max$ operator, and the tangential Coulomb law through the face-wise projection onto friction disks.

  \subsubsection{Unified formulation}
  \label{sec:unified_problem}

  The two preceding examples have different physical origins and different nonlinear mechanisms. Nevertheless, after discretization, both lead to a nonlinear nonsmooth primal-dual system involving a state equation coupled with a projected complementarity equation. We therefore introduce the following unified formulation, which will serve as the high-dimensional model throughout the rest of the paper.

  Let $\mathcal{P} \subset \mathbb{R}^p$ be the parameter space. We denote by $N$ and $M$ the HDM dimensions of the primal and dual variables, respectively. For a given parameter $\mu \in \mathcal{P}$, we consider the following problem: find $U_\mu \in \mathbb{R}^{N}$ and $\Lambda_\mu \in \mathbb{R}^{M}$ such that
  \begin{equation}
  \begin{cases}
      A_\mu(U_\mu) + B_\mu^T \Lambda_\mu - L_\mu = 0, \\
      \Lambda_\mu - \Pi_{\mathcal{K}_\mu(\Lambda_\mu)} \left( \Lambda_\mu + \rho g_\mu(U_\mu) \right) = 0.
  \end{cases}
  \label{eq:unified_problem}
  \end{equation}
  Here $A_\mu : \mathbb{R}^N \to \mathbb{R}^N$ is a possibly nonlinear operator, $B_\mu \in \mathbb{R}^{M \times N}$ is a primal-dual coupling matrix, $L_\mu \in \mathbb{R}^N$ is the load vector, and $g_\mu : \mathbb{R}^N \to \mathbb{R}^M$ maps the primal variable to the discrete gap or constraint
  residual. The set $\mathcal{K}_\mu(\Lambda_\mu) \subset \mathbb{R}^M$ denotes the admissible set for the dual variable; it may depend on $\Lambda_\mu$. The parameter $\rho>0$ is the projection parameter, and $\Pi_{\mathcal{K}_\mu(\Lambda_\mu)}$ denotes the projection onto $\mathcal{K}_\mu(\Lambda_\mu)$. In the frictional contact case, the dependence of the admissible set on the normal multiplier reflects the Coulomb threshold and gives the problem a quasi-variational structure. Accordingly, \eqref{eq:unified_problem} is used here only as a unifying discrete residual formulation, and no general existence or uniqueness result is claimed for this unified class of problems.
    
    \paragraph{Connection with the model problems.}
    
  For the obstacle problem of Section~\ref{sec:2D_obstacle_problem}, the abstract quantities are given by
  \[
      A_\mu(U) = K U + \gamma M U^{\circ 3}, \qquad B_\mu = -M, \qquad L_\mu = 0, \qquad g_\mu(U) = G_\mu - U.
  \]
  The admissible set is
  $\mathcal{K}_\mu = \left\{ \Lambda \in \mathbb{R}^N : \Lambda \geq 0 \right\},$
  which does not depend on $\Lambda_\mu$, and the projection
  $\Pi_{\mathcal{K}_\mu}$ is the componentwise operator $\max(0,\cdot)$.

  For the contact problem of Section~\ref{sec:3D_contact_problem}, the primal variable is
  $U\in\mathbb{R}^{2N}$ and the dual variable is
  $\Lambda=(\Lambda_n,\Lambda_\tau)\in\mathbb{R}^{R+2R}$. In this case,
  \[
      A_\mu(U) = K U,
      \qquad
      B_\mu^T\Lambda = B_n^T\Lambda_n + B_\tau^T\Lambda_\tau,
  \]
  and the constraint residual is
  \[
      g_\mu(U)
      =
      \begin{pmatrix}
          B_n U \\
          B_\tau U
      \end{pmatrix}
      \in \mathbb{R}^{R+2R}.
  \]
  The associated projection acts componentwise. Its normal component is the projection onto the positive
  orthant, corresponding to the operator $\max(0,\cdot)$, while its tangential component is the face-wise
  projection onto the friction disks of radii $c_\mu(\Lambda_{\mu,n})$, yielding the operator
  $\Pi_{c_\mu(\Lambda_{\mu,n})}$ in \eqref{eq:non_smooth_eq_friction_modelproblem}. Thus, the contact problem is
  included in the unified formulation through the multiplier-dependent admissible set
  $\mathcal{K}_\mu(\Lambda_\mu)$.
  
\subsection{Semi-smooth Newton Method}
\label{sec:ssn}

For the rest of the section, we drop the $\mu$ subscript from all variables and operators for readability. The unified problem \eqref{eq:unified_problem} can be written as the nonlinear nonsmooth system
\begin{equation}
  \mathcal{F}(U,\Lambda)
  \coloneqq
  \begin{pmatrix}
    A(U) + B^T\Lambda - L \\
    \Lambda - \Pi_{\mathcal{K}(\Lambda)}
    \bigl(\Lambda+\rho g(U)\bigr)
  \end{pmatrix}
  = 0,
  \label{eq:ssn_residual}
\end{equation}
where the nonsmoothness of $\mathcal{F}$ comes from the projection operator $\Pi_{\mathcal{K}(\Lambda)}$ associated with the constraint set. Because $\mathcal{F}$ is not Fréchet differentiable, we rely on a generalized notion of differentiability, slant-differentiability, which we briefly recall below; see \cite{PDAS_CNQ,PDAS_Q,PDAS-SSN} for a more detailed presentation.

\begin{definition}[Slant-differentiability]
  Let $\mathcal{V}$ and $\mathcal{W}$ be Banach spaces, $\mathcal{D} \subset \mathcal{V}$ an open set, and $\mathcal{F} : \mathcal{D} \to \mathcal{W}$ a nonlinear operator.
  We say that $\mathcal{F}$ is \emph{slant-differentiable} in an open subset $\mathcal{O} \subset \mathcal{D}$ if there exists a mapping $\mathcal{F}' : \mathcal{O} \to \mathcal{L}(\mathcal{V}, \mathcal{W})$ such that, for all $x \in \mathcal{O}$,
  \[
    \lim_{h \to 0}
    \frac{
      \|\mathcal{F}(x+h) - \mathcal{F}(x) - \mathcal{F}'(x+h)\,h
      \|_\mathcal{W}
    }{\|h\|_\mathcal{V}} = 0.
  \]
  The mapping $\mathcal{F}'$ is called a \emph{slanting function} for $\mathcal{F}$ in $\mathcal{O}$.
\end{definition}

A canonical example, relevant to both test cases, is the componentwise positive-part mapping $Q:\mathbb R^m\to\mathbb R^m$ defined by $Q(x)=[x]_+$. Although $Q$ is not Fréchet differentiable at points where at least one component $x_i$ is equal to zero, a slanting function is
  \begin{equation}
  Q'(x)=D_{>0}(x), \qquad \bigl[D_{>0}(x)\bigr]_{ii} =
  \begin{cases}
    1, & x_i>0,\\
    0, & x_i\leq0.
  \end{cases}
  \label{eq:slanting_function_max}
\end{equation}
  The semi-smooth Newton method for solving $\mathcal{F}(x) = 0$ is then defined exactly as the classical Newton method, with $\mathcal{F}'$ playing the role of the Jacobian. Under slant-differentiability of $\mathcal{F}$ near a solution, nonsingularity of the slanting functions, and local boundedness of their inverses, the standard theory \cite{PDAS_CNQ} guarantees local superlinear convergence.

At each step, the computational cost has two main contributions: the solution of the $(N+M)$-dimensional linear system involving $\mathcal{F}'(U^k,\Lambda^k)$, and the assembly of the residual and slanting function, which involve operations scaling with the HDM dimensions $N$ and $M$. As a byproduct of the slanting-function structure, an equivalent primal-dual active-set formulation emerges naturally at each iteration, which is the form adopted in the numerical experiments and later inherited by the reduced-order models.

\paragraph{Line search and stopping criteria.}
To improve robustness, each SSN step is damped using a backtracking line search based on a problem-dependent merit function. Convergence is assessed through a stopping residual that combines the relevant equilibrium and constraint equations. The precise line-search rules, merit functions, residual norms, and stopping criteria used are specified in Appendix~\ref{app:ssn_numerical_settings}.

\begin{algorithm}[h!]
\caption{Semi-smooth Newton method for the unified HDM problem
\eqref{eq:unified_problem}}
\begin{algorithmic}[1]
\State \textbf{Input:} Initial iterate $(U^0,\Lambda^0)$, prescribed
tolerance, and maximum number of iterations $K_{\max}$.
\State Set $k=0$ and evaluate the problem-specific stopping criterion
at $(U^0,\Lambda^0)$.
\While{$k<K_{\max}$ and not converged}
\State Compute the slanting function $\mathcal F'(U^k,\Lambda^k)$.
\State Compute the SSN direction $(d_U^k,d_\Lambda^k)$ by solving
\[ \mathcal F'(U^k,\Lambda^k) 
  \begin{pmatrix} 
    d_U^k\\ d_\Lambda^k 
  \end{pmatrix} 
= -\mathcal F(U^k,\Lambda^k). \]
\State Select $\alpha_k\in(0,1]$ using the problem-specific
backtracking rule described in
Appendix~\ref{app:ssn_numerical_settings}.
\State Set
\[ U^{k+1}=U^k+\alpha_kd_U^k, \qquad \Lambda^{k+1}=\Lambda^k+\alpha_kd_\Lambda^k. \]
\State Set $k=k+1$.
\State Evaluate the problem-specific stopping criterion at
$(U^k,\Lambda^k)$.
\EndWhile
\State \textbf{Output:} Approximate solution $(U^k,\Lambda^k)$.
\end{algorithmic}
\label{eq:alg:ssn}
\end{algorithm}

The explicit slanting functions for the two test cases of Sections~\ref{sec:2D_obstacle_problem} and~\ref{sec:3D_contact_problem}, together with their active-set reformulations, are derived in Appendix~\ref{sec:appendix_ssn_derivation}.

\section{Reduced-order modeling framework}
\label{chap:MOR}
  This section develops the reduced formulations used to reduce the online cost of the high-dimensional semi-smooth Newton solver introduced in Section~\ref{chap:SSN}. We first derive a projection-based reduced-order model in which separate POD spaces are constructed for the primal and dual variables. We then introduce a neural-network-augmented ROM (NN-ROM), inspired by \cite{NN_ROM}, which is embedded directly within the semi-smooth Newton iterations. Finally, we address the remaining online bottleneck associated with residual evaluations and assemblies that still scale with the HDM dimensions $N$ and $M$. For that purpose, we use a hyper-reduction strategy based on sparse cubature rules constructed by greedy nonnegative least squares \cite{greedy_cubature}, yielding hyper-reduced versions of both the ROM and the NN-ROM.

\subsection{POD-based reduced-order model}
\label{sec:linear_rom}

    We consider the high-dimensional primal--dual system introduced in Section~\ref{chap:SSN} and seek low-dimensional approximations of the corresponding primal and dual solution manifolds,
    \begin{equation}
          \label{eq:linear_manifolds}
          \mathcal{M}_U = \{U_\mu : \mu \in \mathcal{P}\},
          \qquad
          \mathcal{M}_\Lambda = \{\Lambda_\mu : \mu \in \mathcal{P}\},
      \end{equation}
    by low-dimensional linear spaces of dimensions $n \ll N$ and $m \ll M$, respectively.   
    
    Let $V\in\mathbb{R}^{N\times n}$ and $W\in\mathbb{R}^{M\times m}$ be matrices whose columns form bases of the reduced primal and dual spaces. The associated linear reconstruction maps are $\mathcal U_\mu(q_\mu)=Vq_\mu$ and $\mathcal D_\mu(\xi_\mu)=W\xi_\mu$, so that the high-dimensional variables are approximated as
    \begin{equation}
        \label{eq:linear_rom_approximation}
        U_\mu \approx \mathcal U_\mu(q_\mu)=Vq_\mu,
        \qquad
        \Lambda_\mu \approx \mathcal D_\mu(\xi_\mu)=W\xi_\mu,
    \end{equation}
    where $q_\mu \in\mathbb{R}^n$ and $\xi_\mu\in\mathbb{R}^m$ denote the primal and dual reduced coordinates. Possible affine offsets or lifting terms are omitted to simplify the notation.
            
    The reduced bases are constructed by POD from high-dimensional snapshots. Given a training set $\mathcal{P}_{\mathrm{POD}} = \{\mu^1,\ldots,\mu^{N_s}\}\subset \mathcal{P}$, we compute the corresponding HDM solutions with the SSN method and collect the primal and dual snapshots in    
    \[ S_U =
    \begin{pmatrix}
    U_{\mu^1} \mid \cdots \mid U_{\mu^{N_s}}
    \end{pmatrix}
    \in \mathbb{R}^{N\times N_s}, \qquad
    S_\Lambda =
    \begin{pmatrix}
    \Lambda_{\mu^1} \mid \cdots \mid \Lambda_{\mu^{N_s}}
    \end{pmatrix}
    \in \mathbb{R}^{M\times N_s}. \]
    The POD bases are then obtained by computing the singular value decompositions
    \[
    S_U = \widetilde{V}\Sigma_U \widetilde{Z}_U^T,
    \qquad
    S_\Lambda = \widetilde{W}\Sigma_\Lambda \widetilde{Z}_\Lambda^T,
    \]
    and retaining the first $n$ and $m$ left singular vectors, respectively:
    \[
    V=\widetilde{V}_{(:,1:n)},
    \qquad
    W=\widetilde{W}_{(:,1:m)}.
    \]
    The resulting bases verify $V^T V=I_n$ and $W^T W=I_m$.    

    For a new parameter value $\mu\in\mathcal{P}$, the reduced coordinates are determined by inserting the approximations~\eqref{eq:linear_rom_approximation} into the HDM residual equations and imposing Galerkin orthogonality with respect to the primal and dual reduced spaces. The reduced problem reads: find $(q_\mu,\xi_\mu)\in\mathbb{R}^n\times\mathbb{R}^m$ such that
    \begin{equation}
    \label{eq:reduced_problem}
        \begin{aligned}
        V^T A_\mu(Vq_\mu) + V^T B_\mu^T W\xi_\mu - V^T L_\mu &= 0,\\
        \xi_\mu - W^T \Pi_{\mathcal{K}_\mu(W\xi_\mu)}
        \bigl(W\xi_\mu + \rho g_\mu(Vq_\mu)\bigr) &= 0.
        \end{aligned}
        \end{equation}
    In the second equation, we have used the Euclidean orthonormality of the dual POD basis. The resulting nonlinear nonsmooth system has dimension $n+m$ and is solved by the same SSN strategy as the HDM.
    
    \begin{rmk}
        For variational inequalities, reduced representations of the multiplier are often designed to preserve dual feasibility, for instance by using nonnegative combinations of multiplier snapshots, nonnegative matrix factorization, or cone-greedy constructions \cite{MOR_Salomon,Charbel_VI,Amina}. Such approaches are particularly natural when the reduced formulation enforces nonnegative reduced coefficients, so that the lifted multiplier remains in the admissible cone. In the present SSN formulation, however, the reduced multiplier coefficients are unconstrained Newton unknowns and the complementarity conditions are enforced through the projected reduced residual. A cone-preserving basis would therefore not, by itself, guarantee pointwise feasibility of the reconstructed multiplier. We instead construct both the primal and dual spaces by POD, which minimizes the reconstruction error of the training snapshots for a prescribed reduced dimension. 
        Pointwise feasibility of the reconstructed multiplier is therefore not guaranteed. Its feasibility, together with the mechanical effect of the dual reconstruction error, is summarized in Sections~\ref{sec:results_obstacle} and~\ref{sec:results_contact} and numerically examined in detail in Appendix~\ref{app:constraint_diagnostics}. 
    \end{rmk}

     \begin{rmk}[Projection parameter in the reduced formulation]
        At the HDM level, the complementarity solution is independent of the projection parameter $\rho>0$. This property is generally lost after reduction, since only the projected complementarity residual is enforced. Therefore, $\rho$ becomes a parameter of the reduced formulation. In the numerical experiments, it is selected empirically from preliminary reduced-model tests based primarily on the primal validation error.
    \end{rmk}
\subsection{Neural Network-Augmented Reduced Bases}
\label{sec:nn_for_rom}

    The efficiency of the linear ROM depends on how well the primal and dual solution manifolds $\mathcal{M}_U$ and $\mathcal{M}_\Lambda$ in \eqref{eq:linear_manifolds} can be approximated by low-dimensional linear subspaces. This approximability is commonly characterized by the Kolmogorov $n$-width. For a set $\mathcal{M}$ in a normed space $X$, it is defined by
    \[d_n(\mathcal{M}, X) =
    \inf_{\substack{V_n \subset X \\ \dim(V_n)=n}} \sup_{v\in\mathcal{M}} \inf_{v_n\in V_n} \|v-v_n\|_X,\]
    which quantifies the best achievable worst-case error over all $n$-dimensional linear subspaces. Rapid decay of the Kolmogorov $n$-width indicates that accurate linear reduced models can be obtained with small reduced dimensions, whereas slow decay implies that many POD modes may be required to reach a given accuracy. In practice, the decay of the POD singular values provides a computable indicator of snapshot compressibility, although it does not by itself certify the Kolmogorov $n$-width of the full solution manifold.
    
    The limitations of linear representations motivate the adoption of nonlinear reduced-order manifolds. To this end, we extend the NN-augmented ROM framework of \cite{NN_ROM} to the primal-dual VI setting considered here. Reduced-order solutions are first computed in a sufficiently rich POD space to ensure satisfactory accuracy, and subsequently used as training data. Two feedforward neural networks, one per variable, are then trained to map the leading POD coordinates to the remaining ones, thereby defining a nonlinear manifold correction in which the trailing modes are expressed as a learned function of a small set of dominant modes. In the present work, parameter features are also provided as inputs to the networks.  
        
    The construction of the NN-augmented ROM is organized into an offline stage, in which the POD bases are fixed, the training data are generated, and the neural mappings are learned, followed by an online stage, in which the reduced nonlinear nonsmooth system is solved for new parameter values.
        
\subsubsection{Offline construction and training}
\label{sec:nn_offline}
    \begin{figure}[pos=t!]
        \centering
        \resizebox{\linewidth}{!}{%
            \begingroup

\definecolor{nnromPodBlue}{HTML}{1D4E91}
\definecolor{nnromPodBlack}{HTML}{1A1A1A}
\definecolor{nnromTeal}{HTML}{0F746C}
\definecolor{nnromDarkGray}{HTML}{404448}
\definecolor{nnromLightBlue}{HTML}{F4F7FB}
\definecolor{nnromLightTeal}{HTML}{F2F9F8}
\definecolor{nnromLightGray}{HTML}{F7F7F7}

\tikzset{
  nnrom data/.style={
    draw=none,
    fill=none,
    align=center,
    inner sep=1pt,
    font=\small
  },
  nnrom operation/.style={
    draw,
    rounded corners=2pt,
    line width=0.65pt,
    align=center,
    inner xsep=6pt,
    inner ysep=5pt,
    font=\small
  },
  nnrom pod operation/.style={
    nnrom operation,
    draw=nnromPodBlue,
    fill=nnromLightBlue
  },
  nnrom nn operation/.style={
    nnrom operation,
    draw=nnromTeal,
    fill=nnromLightTeal
  },
  nnrom online operation/.style={
    nnrom operation,
    draw=nnromDarkGray,
    fill=nnromLightGray
  },
  nnrom flow/.style={
    ->,
    draw=nnromDarkGray,
    line width=0.58pt,
    line cap=round,
    line join=round
  },
  nnrom pod flow/.style={nnrom flow,draw=nnromPodBlue,line width=0.68pt},
  nnrom guide/.style={draw=nnromDarkGray!55,line width=0.45pt},
  nnrom boundary/.style={draw=nnromDarkGray,dashed,line width=0.55pt}
}

\begin{tikzpicture}[x=1.08cm,y=1cm]

  \node[anchor=west,font=\bfseries\large,text=nnromPodBlue]
    at (-1.55,4.55) {POD space construction};

  \node[anchor=west,font=\bfseries\large,text=nnromTeal]
    at (-1.55,1.90) {NN training-data generation and learning};

  \draw[nnrom boundary] (17.20,-2.85) -- (17.20,4.72);
  \node[anchor=east,font=\bfseries\small] at (17.00,4.55) {OFFLINE};
  \node[anchor=west,font=\bfseries\small] at (17.40,4.55) {ONLINE};

  \node[nnrom data,text width=2cm] (pod-samples) at (0.00,3.32) {
    POD snapshot parameters\\[2pt]
    $\left\{\mu_{\mathrm{POD}}^i\right\}_{i=1}^{N_s}$
  };

  \node[nnrom pod operation,minimum width=1.85cm,minimum height=0.95cm]
    (hdm-solves) at (3.20,3.32) {HDM\\solver};

  \node[nnrom data,text width=1.5cm] (snapshots) at (6.00,3.32) {
    Snapshot matrices\\[2pt]
    $S_U, S_\Lambda$
  };

  \node[nnrom pod operation,minimum width=1.20cm,minimum height=0.88cm]
    (pod) at (8.55,3.32) {POD};

  \node[nnrom data,text width=3.0cm] (bases) at (11.95,3.32) {
    Enlarged POD bases\\[2pt]
    $V^{\mathrm{tot}}=\left[\,V^r\mid V^c\,\right]$\\[-1pt]
    $W^{\mathrm{tot}}=\left[\,W^r\mid W^c\,\right]$
  };

  \draw[nnrom flow] (pod-samples.east) -- (hdm-solves.west);
  \draw[nnrom flow] (hdm-solves.east) -- (snapshots.west);
  \draw[nnrom flow] (snapshots.east) -- (pod.west);
  \draw[nnrom flow] (pod.east) -- (bases.west);

  \node[nnrom data,text width=1.9cm] (nn-samples) at (0.00,-0.20) {
    NN-training parameters\\[2pt]
    $\left\{\mu_{\mathrm{NN}}^i\right\}_{i=1}^{N_{\mathrm{train}}}$
  };

  \node[nnrom nn operation,text width=2.35cm,minimum height=1.15cm]
    (rom-solves) at (3.55,0.52) {
      ROM \\ solver
    };

  \node[nnrom nn operation,text width=2.15cm,minimum height=0.92cm]
    (feature-map) at (3.55,-1.18) {
      Parameter-feature map\\[-1pt]
      $\mu\longmapsto\eta_\mu$
    };

  \node[nnrom data,text width=2.25cm] (coordinates) at (6.95,0.52) {
    Reduced coefficients\\[2pt]
    $q_\mu^{\mathrm{tot}}=(q_\mu^r,q_\mu^c)^T$\\[1pt]
    $\xi_\mu^{\mathrm{tot}}=(\xi_\mu^r,\xi_\mu^c)^T$
  };

  \node[nnrom data,text width=2.25cm] (features) at (6.95,-1.18) {
    Parameter features\\[2pt]
    $\eta_\mu$
  };

  \coordinate (sample-split) at (1.62,-0.20);
  \draw[draw=nnromDarkGray,line width=0.58pt]
    (nn-samples.east) -- (sample-split);
  \draw[nnrom flow] (sample-split) |- (rom-solves.west);
  \draw[nnrom flow] (sample-split) |- (feature-map.west);

  \draw[nnrom flow] (rom-solves.east) -- (coordinates.west);
  \draw[nnrom flow] (feature-map.east) -- (features.west);

  \draw[nnrom pod flow]
    (bases.south) |- (3.55,1.50) -- (rom-solves.north);

  \node[nnrom data,text width=2.3cm] (training-pairs) at (10.35,-0.20) {
    Training pairs\\[2pt]
    $(q_\mu^r,\eta_\mu)\longmapsto q_\mu^c$\\[1pt]
    $(\xi_\mu^r,\eta_\mu)\longmapsto\xi_\mu^c$
  };

  \node[nnrom nn operation,minimum width=1.55cm,minimum height=0.98cm]
    (nn-training) at (13.30,-0.20) {NN\\training};

  \node[nnrom data,text width=1.4cm] (trained-maps) at (15.80,-0.20) {
    NN maps\\[2pt]
    $\mathcal N_{\theta_u}^{u}$\\[-1pt]
    $\mathcal N_{\theta_\lambda}^{\lambda}$
  };

  \draw[nnrom flow]
    (coordinates.east) -- ++(0.15,0) |- ([yshift=0.34cm]training-pairs.west);
  \draw[nnrom flow]
    (features.east) -- ++(0.15,0) |- ([yshift=-0.34cm]training-pairs.west);
  \draw[nnrom flow] (training-pairs.east) -- (nn-training.west);
  \draw[nnrom flow] (nn-training.east) -- (trained-maps.west);

  \node[nnrom online operation,text width=2.50cm,minimum height=1.15cm]
    (online-solve) at (19.65,-0.20) {
      Online \\
      NN-ROM
      solver
    };

  \draw[nnrom flow] (trained-maps.east) -- (online-solve.west);
  \draw[nnrom pod flow]
    (bases.east) -- ++(1.85,0) -| (online-solve.north);

\end{tikzpicture}
\endgroup%
        }
        \caption{Offline construction of the NN-ROM. POD of the HDM snapshot
        matrices provides the enlarged primal and dual bases. Galerkin ROM
        solves in these spaces, together with the parameter features, provide
        the data used to train the correction maps.}
        \label{fig:nn_rom_offline_pipeline}
    \end{figure}
    
    \paragraph{NN-ROM approximation.}
    We first construct enlarged linear POD spaces for the primal and dual variables, with bases $V^{\mathrm{tot}}\in\mathbb{R}^{N\times n_{\mathrm{tot}}}$ and $W^{\mathrm{tot}}\in\mathbb{R}^{M\times m_{\mathrm{tot}}}$. For prescribed dimensions $n_r<n_{\mathrm{tot}}$ and $m_r<m_{\mathrm{tot}}$, we split these bases as
    \[V^{\mathrm{tot}} =
        \begin{pmatrix}
            V^r \mid V^c
        \end{pmatrix},
    \qquad
        W^{\mathrm{tot}} =
        \begin{pmatrix}
            W^r \mid W^c
        \end{pmatrix},\]
    where $V^r\in\mathbb{R}^{N\times n_r}$ and $W^r\in\mathbb{R}^{M\times m_r}$ contain the retained leading POD modes, while $V^c\in\mathbb{R}^{N\times n_c}$ and $W^c\in\mathbb{R}^{M\times m_c}$ contain the complementary modes, with $n_c=n_{\mathrm{tot}}-n_r$ and $m_c=m_{\mathrm{tot}}-m_r$.   

    The NN-ROM approximates the complementary coordinates as functions of the retained coordinates and of problem-dependent parameter features. More precisely, two neural networks
    \[
    \mathcal N_{\theta_u}^u:\mathbb{R}^{n_r+p_\eta}\to\mathbb{R}^{n_c},
    \qquad
    \mathcal N_{\theta_{\lambda}}^{\lambda}:\mathbb{R}^{m_r+p_\eta}\to\mathbb{R}^{m_c}
    \]
    are trained to approximate the maps $(q^r_\mu,\eta_\mu)\mapsto q^c_\mu$ and $(\xi^r_\mu,\eta_\mu)\mapsto\xi^c_\mu$, where $\eta_\mu\in\mathbb{R}^{p_\eta}$ denotes a vector of $p_\eta$ features constructed from the physical parameters. These features are chosen so as to provide the networks with a convenient representation of the parametric dependence, for instance by encoding periodic quantities or normalized parameters.
    
    The NN-ROM then replaces the linear reconstruction maps introduced in \eqref{eq:linear_rom_approximation} by the nonlinear reconstructions,
    \begin{equation}
        \label{eq:nn_rom_reconstruction}
        \mathcal U_\mu(q_\mu^r) = V^r q_\mu^r + V^c\,\mathcal N_{\theta_u}^u(q_\mu^r,\eta_\mu),
        \qquad
        \mathcal D_\mu(\xi_\mu^r) = W^r \xi_\mu^r + W^c\,\mathcal N_{\theta_{\lambda}}^{\lambda}(\xi_\mu^r,\eta_\mu).
    \end{equation}
    
\paragraph{Training data.}
    Training requires pairs of leading and complementary coordinates $(q_{\mu^i}^r,q_{\mu^i}^c)$ and $(\xi_{\mu^i}^r,\xi_{\mu^i}^c)$ for parameters $\mathcal P_{\mathrm{NN}} = \{\mu^1_{\mathrm{NN}},\ldots,\mu^{N_{\mathrm{train}}}_{\mathrm{NN}}\}\subset\mathcal P$.
    
    A simple possibility would be to project the HDM snapshots associated with $\mathcal P_{\mathrm{POD}}$ onto the enlarged POD bases, as done in \cite{NN_ROM} for nonlinear PDEs:
      \[ q_{\mu^i}^{\mathrm{tot}}=(V^{\mathrm{tot}})^T U_{\mu^i}, \qquad \xi_{\mu^i}^{\mathrm{tot}}=(W^{\mathrm{tot}})^T \Lambda_{\mu^i}. \]
    However, while these projection-based coordinates minimize the reconstruction errors in the chosen POD inner product, they are obtained by approximating the primal and dual snapshots separately, and are not required to satisfy the reduced primal--dual equations. In contrast, the coordinates computed by the reduced SSN solver are determined by the projected stationarity and complementarity conditions simultaneously: the primal coordinates are coupled to the dual variable through the term $B_\mu^T W\xi_\mu$, while the dual coordinates are coupled to the primal variable through $g_\mu(Vq_\mu)$. Consequently, the projected pair $((V^\mathrm{tot})^T U_\mu, (W^\mathrm{tot})^T \Lambda_\mu)$ does not in general coincide with the pair that solves the reduced system.
  
    For this reason, we generate the training low-fidelity snapshots from reduced solves rather than from direct projection. More precisely, for each training parameter $\mu \in \mathcal{P}_{\mathrm{NN}}$, we solve the Galerkin ROM in the enlarged spaces $\operatorname{span}(V^{\mathrm{tot}})$ and $\operatorname{span}(W^{\mathrm{tot}})$. The resulting coordinates
        \[ q_{\mu}^{\mathrm{tot}} = \begin{pmatrix} q^r_{\mu}\\ q^c_{\mu} \end{pmatrix}, \qquad \xi_{\mu}^{\mathrm{tot}} = \begin{pmatrix} \xi^r_{\mu}\\ \xi^c_{\mu} \end{pmatrix}, \qquad \mu \in \mathcal{P}_{\mathrm{NN}}, \]
    provide training targets that are consistent with the reduced equations solved in the online NN-ROM.

    The NN-ROM offline workflow is summarized in Figure~\ref{fig:nn_rom_offline_pipeline}.
    The NN-ROM introduces an additional offline cost, mainly associated with generating the low-fidelity reduced-solution database rather than with the neural-network optimization itself. Its magnitude relative to the HDM snapshot campaign is problem-dependent and is quantified separately for each test case in the numerical results.  
    
    \paragraph{Loss function.}
    In the NN-ROM approach of \cite{NN_ROM}, the neural networks are trained by minimizing a mean squared error on the complementary reduced coordinates. Here, we instead use an energy-norm loss for the primal network, so that errors in the predicted complementary coordinates are weighted according to their contribution to the reconstructed field. Let $K\in\mathbb{R}^{N\times N}$ denote the stiffness matrix defining the primal energy norm, $\|z\|_K^2=z^T K z$. For the obstacle problem, $K$ is the scalar elliptic stiffness matrix appearing in \eqref{eq:discrete_obstacle_mix_form_CP}, while for the contact problem it is the block elasticity matrix appearing in \eqref{eq:non_smooth_eq_friction_modelproblem}. The primal network is trained by minimizing
    
    \begin{equation}
        \label{eq:primal_nn_loss}
        \theta_u
        =
        \arg\min_{\theta}
        \frac{1}{N_{\mathrm{train}}}
        \sum_{i=1}^{N_{\mathrm{train}}}
        \left\|
        V^c\left(
        \mathcal N_{\theta}^{u}(q_{\mu^i}^r,\eta_{\mu^i})
        -
        q^c_{\mu^i}
        \right)
        \right\|_K^2.
    \end{equation}
    
For the dual network, we use the same residual-aware principle. Since the dual variable enters the primal residual through $B_\mu^T\Lambda$, the complementary-coordinate error produced by the network induces the primal residual error
\[
B_{\mu^i}^T W^c
\left(
\mathcal N_{\theta}^{\lambda}(\xi_{\mu^i}^r,\eta_{\mu^i})
-
\xi^c_{\mu^i}
\right).
\]
Here, $B_\mu^T$ is the primal-dual coupling operator: for the obstacle problem, $B_\mu^T=-M$, whereas for the contact problem, writing $\Lambda=(\Lambda_n,\Lambda_\tau)$, one has $B_\mu^T\Lambda = B_n^T\Lambda_n + B_\tau^T\Lambda_\tau.$

We measure this residual error in the dual norm associated with the primal energy norm, namely $\|z\|_{K^{-1}}^2=z^T K^{-1}z$, and train the dual network by minimizing
    \begin{equation}
    \label{eq:dual_nn_loss}
    \theta_\lambda
    =
    \arg\min_{\theta}
    \frac{1}{N_{\mathrm{train}}}
    \sum_{i=1}^{N_{\mathrm{train}}}
    \left\|
    B_{\mu^i}^T W^c
    \left(
    \mathcal N_{\theta}^{\lambda}(\xi_{\mu^i}^r,\eta_{\mu^i})
    -
    \xi^c_{\mu^i}
    \right)
    \right\|_{K^{-1}}^2
    \end{equation}
    Thus, the dual loss penalizes the neural correction according to the size of the primal residual error it induces, rather than according to a componentwise error on the reduced dual coefficients.

\subsubsection{Online stage: NN-ROM solve}
\label{sec:nn_online}

    For a new parameter value $\mu\in\mathcal{P}$, the online unknowns are the retained coordinates $(q_\mu^r,\xi_\mu^r)\in\mathbb{R}^{n_r}\times\mathbb{R}^{m_r}$. 
    Using the nonlinear reconstructions defined in \eqref{eq:nn_rom_reconstruction}, the tangent matrices with respect to the online reduced coordinates are
    \begin{equation}
    \label{eq:tangent_matrices_NN_ROM}
        T_{\mathcal U_\mu}(q^r) := \frac{\partial \mathcal U_\mu}{\partial q^r}(q^r) = V^r + V^c\,J_{\mathcal N_{\theta_u}^u}(q^r,\eta_\mu) \in\mathbb{R}^{N\times n_r},
        \quad
        T_{\mathcal{D}_\mu}(\xi^r) := \frac{\partial \mathcal{D}_\mu}{\partial \xi^r}(\xi^r) = W^r + W^c\,J_{\mathcal N_{\theta_\lambda}^{\lambda}}(\xi^r,\eta_\mu) \in\mathbb{R}^{M\times m_r},
    \end{equation}
    where $J_{\mathcal N^u_{\theta_u}}(q^r,\eta_\mu)$ and
    $J_{\mathcal N^\lambda_{\theta_\lambda}}(\xi^r,\eta_\mu)$ denote the Jacobians of the neural networks with respect to the reduced-coordinate arguments, with the parameter features $\eta_\mu$ fixed.  
    
      Let 
    \[ r_\mu(U_\mu,\Lambda_\mu) := A_\mu(U_\mu)+B_\mu^T\Lambda_\mu-L_\mu \in\mathbb{R}^N, \qquad c_\mu(U_\mu,\Lambda_\mu) := \Lambda_\mu-\Pi_{\mathcal K_\mu(\Lambda_\mu)} \bigl(\Lambda_\mu+\rho g_\mu(U_\mu)\bigr) \in\mathbb{R}^M \]
  denote the state and constraint residuals. The NN-ROM is obtained by enforcing the orthogonality of these residuals to the corresponding tangent spaces, yielding a Petrov-Galerkin closure. Thus, for a parameter $\mu \in \mathcal P$, we seek $(q_\mu^r,\xi_\mu^r)\in\mathbb{R}^{n_r}\times\mathbb{R}^{m_r}$ such that,
    \begin{equation} \label{eq:nn_rom_problem} 
      \begin{aligned} 
        T_{\mathcal U_\mu}(q^r_\mu)^T r_\mu\bigl(\mathcal U(q^r_\mu),\mathcal{D}(\xi^r_\mu)\bigr) &=0,\\ 
        T_{\mathcal{D}_\mu}(\xi^r_\mu)^T c_\mu\bigl(\mathcal U(q^r_\mu),\mathcal{D}(\xi^r_\mu)\bigr) &=0. 
      \end{aligned} 
    \end{equation} 

The tangent-space system \eqref{eq:nn_rom_problem} is solved by an
approximate semi-smooth Newton method. At a current iterate
$(q^{r},\xi^{r})$, we use the block linearization
\[
\widetilde{\mathcal{J}}_{\mu}
=
\begin{bmatrix}
T_{\mathcal{U},\mu}^{T}
\,\partial_{U}r_{\mu}\,
T_{\mathcal{U},\mu}
&
T_{\mathcal{U},\mu}^{T}
\,\partial_{\Lambda}r_{\mu}\,
T_{\mathcal{D},\mu}
\\[1ex]
T_{\mathcal{D},\mu}^{T}
\,\partial_{U}c_{\mu}\,
T_{\mathcal{U},\mu}
&
T_{\mathcal{D},\mu}^{T}
\,\partial_{\Lambda}c_{\mu}\,
T_{\mathcal{D},\mu}
\end{bmatrix},
\] where all residual derivatives and tangent matrices are evaluated at the current reconstructed state. The derivatives of the nonsmooth projection operators are understood in the slanting sense.

This matrix omits the product-rule terms arising from the dependence of $T_{\mathcal{U},\mu}$ and $T_{\mathcal{D},\mu}$ on the reduced coordinates, namely terms involving derivatives of the tangent matrices multiplied by the current residuals. These terms are omitted in all reported computations because the resulting iterations converge reliably to the prescribed residual tolerances, while the reduced solutions retain satisfactory accuracy in both test cases.
Therefore, $\widetilde{\mathcal{J}}_{\mu}$ is not, in general, an exact slanting function of the tangent-projected residual \eqref{eq:nn_rom_problem}. The resulting iteration should consequently be interpreted as an approximate semi-smooth Newton, or quasi-Newton, method. The local superlinear convergence result recalled for the exact HDM semi-smooth Newton method does not directly apply to this approximation; its robustness and convergence are assessed numerically.
        
    \begin{rmk}
    When the nonlinear correction is removed, i.e. $V^c=0$ and $W^c=0$, one has
    $\mathcal U_\mu(q^r)=V^r q^r$,
    $\mathcal D_\mu(\xi^r)=W^r\xi^r$,
    $T_{\mathcal U,\mu}(q^r)=V^r$, and
    $T_{\mathcal D,\mu}(\xi^r)=W^r$; hence \eqref{eq:nn_rom_problem} reduces to the linear Galerkin ROM \eqref{eq:reduced_problem}.
    \end{rmk}

\subsection{Hyper-reduction by Greedy-NNLS Cubature}
\label{sec:hyper_reduction}

    The ROM and NN-ROM formulations introduced in Sections~\ref{sec:linear_rom} and~\ref{sec:nn_for_rom} reduce the dimension of the nonlinear systems solved by SSN. However, their online evaluation may still involve operations scaling with the HDM dimensions $N$ and $M$.
    
    This can for instance be seen on the two-dimensional obstacle problem. The Galerkin ROM residual reads 
    \begin{equation} 
      \label{eq:reduced_problem_obstacle} 
      \mathcal{F}_{r,\mu}^{(1)}(q,\xi) = 
        \begin{pmatrix} K_n q + \gamma V^T M (V q)^{\circ 3} - M_n \xi \\ 
        \xi - W^T \max\left(0, W\xi - \rho(Vq-G_\mu)\right) \end{pmatrix}
      =0,
    \end{equation} 
    where $K_n=V^T K V$ and $M_n=V^T M W$. Although this system has dimension $n+m$, its online evaluation is not independent of the HDM dimension. Indeed, the residual requires forming the lifted quantities $Vq$ and $W\xi$, evaluating the nonlinear vector $(Vq)^{\circ 3}$, applying the componentwise projection in $\mathbb{R}^N$, and projecting the resulting vectors back with $V^T$ and $W^T$. The same issue arises when forming the linear system solved at each SSN iteration. Differentiating the cubic term requires the diagonal matrix $\operatorname{diag}((Vq)^{\circ 2}) \in \mathbb{R}^{N\times N}$, while differentiating the projection term requires applying the pointwise derivative of the $\max$ operator to the vector $W\xi-\rho(Vq-G_\mu) \in \mathbb{R}^N$. Both operations are defined at the HDM level before being projected by $V^T$ and $W^T$. As shown in the numerical results, these HDM-scale residual and Jacobian evaluations constitute the main online bottleneck of the Galerkin ROM. 
    
    Hyper-reduction methods address this issue by replacing full-dimensional nonlinear or nonaffine contributions with approximations whose online cost depends only on a small number of sampled entries, elements, or quadrature points. Broadly speaking, interpolation-based methods such as EIM, DEIM, and Q-DEIM approximate nonlinear fields in the HDM space before projection \cite{EIM,EIM_2,DEIM,Q-DEIM}. In contrast, empirical quadrature or cubature methods approximate directly the reduced projected quantities entering the Newton or SSN iterations, typically through sparse weighted sums with nonnegative weights \cite{greedy_cubature,ECSW,ECM}.
    
    In this work, we focus on the greedy-NNLS cubature approach of \cite{greedy_cubature}. This choice is well suited to the present formulation because it targets the projected residual and tangent contributions directly, rather than reconstructing full-dimensional nonlinear fields. This distinction will also be important when discussing the interaction between hyper-reduction and the NN-augmented ROM.

\subsubsection{Sparse cubature formulation}
\label{sec:cubature_setup}
    We describe the cubature construction for projected quantities whose direct evaluation scales with the HDM dimension. Specifically, we consider a vector-valued nonlinear contribution of the form $f(\mathcal U_\mu(q))\in\mathbb{R}^N$, where the nonlinear mapping $f:\mathbb{R}^N\to\mathbb{R}^N$ is evaluated at the reconstructed primal state and $q$ denotes a generic reduced coordinate of the ROM or the NN-ROM. This setting covers, for example, the cubic contribution of the obstacle problem and provides a convenient setting to present the sparse weighted-sum construction used below. Other nonlinear contributions, for instance depending on the dual reconstruction or on nonsmooth projection operators, are treated in the same way.
        
    The associated reduced equations involve projected nonlinear terms of the form
    \begin{equation}
        \label{eq:projected_nonlinearity}
        \Phi_\mu(q)^T f\bigl(\mathcal U_\mu(q)\bigr)\in\mathbb{R}^n,
    \end{equation}
    where $\Phi_\mu(q)\in\mathbb{R}^{N\times n}$ denotes the projection matrix associated with this contribution, including any HDM operator applied before projection. For the cubic term in the obstacle problem, for instance, $f(\mathcal U_\mu(q))=\mathcal U_\mu(q)^{\circ 3}$. In the linear ROM, this gives $\Phi_\mu^T=V^T M$, while in the NN-ROM, $\Phi_\mu(q)^T=T_{\mathcal U,\mu}(q)^T M$. 
    \jump
    We assume that $f$ is sample-evaluable: for any index set $\mathcal I\subset\{1,\ldots,N\}$, the entries $f(U)_{\mathcal I}$ can be computed at cost $\mathcal O(|\mathcal I|)$ without assembling the full vector $f(U)$. This is the case for the pointwise nonlinearities considered here, for which $f(U)_i=f_i(U_i)$.

    The projected quantity \eqref{eq:projected_nonlinearity} is a matrix-vector product that can therefore be expanded columnwise as
    \[\Phi_\mu(q)^T f\bigl(\mathcal U_\mu(q)\bigr)
    =
    \sum_{i=1}^{N}
    \bigl(\Phi_\mu(q)^T\bigr)_{(:,i)}
    \bigl(f(\mathcal U_\mu(q))\bigr)_i.\]
    The greedy-NNLS cubature approximation consists in approximating this full sum by a sparse weighted sum over a subset $\mathcal S=\{s_1,\ldots,s_{c_{\mathrm{hr}}}\}\subset\{1,\ldots,N\}$, with $c_{\mathrm{hr}}\ll N$:
    \begin{equation}
    \label{eq:cubature_approx}
    \Phi_\mu(q)^T f\bigl(\mathcal U_\mu(q)\bigr)
    \approx
    \sum_{j=1}^{c_{\mathrm{hr}}}
    \omega_j
    \bigl(\Phi_\mu(q)^T\bigr)_{(:,s_j)}
    \bigl(f(\mathcal U_\mu(q))\bigr)_{s_j},
    \qquad 
    \omega_j\ge 0.
    \end{equation}
    Thus, the method does not approximate the full nonlinear field $f(\mathcal U_\mu(q))\in\mathbb{R}^N$ as in interpolation-based methods, but rather approximates directly its reduced projection $\Phi_\mu(q)^T f(\mathcal U_\mu(q))$.
    
\subsubsection{Offline construction of the cubature rule}
\label{sec:cubature_offline}
The offline phase determines the sampling set $\mathcal S$ and the nonnegative weights $\omega$ used in the cubature approximation \eqref{eq:cubature_approx}. For a set of cubature-training parameters $\mathcal{P}_{\mathrm{cub}} = \{\mu^1_\mathrm{cub},\ldots,\mu^{N_{\mathrm{cub}}}_\mathrm{cub}\}\subset\mathcal P$, we first solve the non-hyper-reduced ROM or NN-ROM, depending on the model to be hyper-reduced. This provides reduced coordinates $(q_{\mu})_{\mu \in \mathcal{P}_{\mathrm{cub}}}$. For each $k \in \{1, \ldots, N_\mathrm{cub} \} $, we define
\[
\Phi_k:=\Phi_{\mu^k_\mathrm{cub}}(q_{\mu^k_\mathrm{cub}}),
\qquad
F_k:=f\bigl(\mathcal U_{\mu^k_\mathrm{cub}}(q_{\mu^k_\mathrm{cub}})\bigr).
\]
The corresponding projected term and its matrix-vector expansion are given by
\[
y_k \coloneqq \Phi_k^T F_k
=
\sum_{i=1}^N
\bigl(\Phi_k^T\bigr)_{(:,i)}\bigl(F_k\bigr)_i
\in\mathbb{R}^n.
\]
Collecting these vectors in
\[
Y=\bigl(y_1\ \mid \cdots\ \mid y_{N_{\mathrm{cub}}}\bigr)
\in\mathbb{R}^{n\times N_{\mathrm{cub}}},
\]
the full contribution can be written as
\[
Y=\sum_{i=1}^N G_i,
\qquad
G_i :=
\begin{pmatrix}
\bigl(\Phi_1^T\bigr)_{(:,i)}\bigl(F_1\bigr)_i
\mid
\cdots
\mid
\bigl(\Phi_{N_{\mathrm{cub}}}^T\bigr)_{(:,i)}\bigl(F_{N_{\mathrm{cub}}}\bigr)_i
\end{pmatrix}
\in\mathbb{R}^{n\times N_{\mathrm{cub}}}.\]
The cubature problem consists in approximating this full sum by a sparse nonnegative weighted sum:
\[
Y \approx \sum_{i=1}^N \widetilde\omega_i G_i,
\qquad
\widetilde\omega\in\mathbb{R}^N_+,
\qquad
\|\widetilde\omega\|_0\le c_{\mathrm{max}},
\]
where $c_{\max}$ is a prescribed maximum number of cubature points.
The selected set is then $\mathcal S=\{i:\widetilde\omega_i\neq 0\}$, and the cubature weights are the corresponding nonzero entries of $\widetilde\omega$.
\jump
Equivalently, after vectorization, let
\[
\mathbf y:=\operatorname{vec}(Y)\in\mathbb{R}^{nN_{\mathrm{cub}}},
\qquad
\mathbf g_i:=\operatorname{vec}(G_i)\in\mathbb{R}^{nN_{\mathrm{cub}}},
\qquad
\mathbf G =
\begin{pmatrix}
\mathbf g_1 \mid \cdots \mid \mathbf g_N
\end{pmatrix}
\in\mathbb{R}^{nN_{\mathrm{cub}}\times N}.
\]
The cubature weights are sought as a sparse nonnegative solution of
\begin{equation}
\label{eq:cubature_sparse_nnls}
\min_{\widetilde\omega\ge 0}
\|\mathbf y-\mathbf G\widetilde\omega\|_2^2
\qquad
\text{subject to}
\qquad
\|\widetilde\omega\|_0\le c_{\mathrm{max}}.
\end{equation}

We solve \eqref{eq:cubature_sparse_nnls} using a greedy non-negative least-squares procedure, following the nonnegative orthogonal matching pursuit (NNOMP) strategy of \cite{NNOMP}. The procedure is controlled by a prescribed residual tolerance $\varepsilon_{\mathrm{tol}}$ and by a maximum number of selected cubature points $c_{\max}$. Starting from $\mathcal S_0=\emptyset$ and $\mathbf r_0=\mathbf y$, at greedy iteration $\ell$ we select
\[
s_\ell
=
\arg\max_{i\notin\mathcal S_{\ell-1}}
\langle \mathbf r_{\ell-1},\mathbf g_i\rangle,
\]
update $\mathcal S_\ell=\mathcal S_{\ell-1}\cup\{s_\ell\}$, and recompute the nonnegative weights on the selected set by solving
\[
\omega_\ell
=
\arg\min_{\omega\in\mathbb{R}^{|\mathcal S_\ell|}_+}
\|\mathbf y-\mathbf G_{\mathcal S_\ell}\omega\|_2^2,
\qquad
\mathbf G_{\mathcal S_\ell}
=
\bigl(\mathbf g_i\bigr)_{i\in\mathcal S_\ell}.
\]
The residual is then updated as
\[
\mathbf r_\ell
=
\mathbf y-\mathbf G_{\mathcal S_\ell}\omega_\ell.
\]
The iteration is stopped when either $|\mathcal S_\ell|=c_{\max}$ or the relative residual $\|\mathbf r_\ell\|_2/\|\mathbf y\|_2$ falls below $\varepsilon_{\mathrm{tol}}$.

The final selected set $\mathcal S$ and weights $\omega$ define the cubature rule used online, and the resulting hyper-reduced dimension is
\[
c_{\mathrm{hr}} = |\mathcal S|.
\]

The sampling set and nonnegative cubature weights are held fixed during the online solve and used for both the hyper-reduced residual and its linearization. In particular, no independent cubature rule is fitted for the generalized Jacobian.

The cubature rule is fitted using converged ROM or NN-ROM solutions. The intermediate semi-smooth Newton iterates are not included in the present cubature-training database. Therefore, the cubature approximation is optimized near the reduced solution manifold, but no uniform accuracy guarantee is claimed away from that manifold.

\subsubsection{Online evaluation of the hyper-reduced terms}
  \label{sec:cubature_online}

In the online phase, the sampling set $\mathcal S=\{s_1,\ldots,s_{c_{\mathrm{hr}}}\}\subset\{1,\ldots,N\}$ and the associated nonnegative weights $\omega\in\mathbb{R}^{c_{\mathrm{hr}}}_+$ are those obtained in the offline cubature construction. We denote by $P_{\mathcal S}\in\mathbb{R}^{N\times c_{\mathrm{hr}}}$ the corresponding selection matrix,
\[
P_{\mathcal S}=\bigl(e_{s_1}\ \mid \cdots\ \mid e_{s_{c_{\mathrm{hr}}}}\bigr).
\]
For a fixed online parameter $\mu$, only the sampled entries of the reconstructed state are required. In the linear ROM, for $q\in\mathbb{R}^n$,
\[
U_{\mathcal S}
:=
P_{\mathcal S}^T\mathcal U_\mu(q)
=
P_{\mathcal S}^T Vq
=
V_{\mathcal S}q,
\qquad
V_{\mathcal S}:=P_{\mathcal S}^T V\in\mathbb{R}^{c_{\mathrm{hr}}\times n}.
\]
In the NN-ROM, for $q^r\in\mathbb{R}^{n_r}$,
\[
U_{\mathcal S}
:=
P_{\mathcal S}^T\mathcal U_\mu(q^r)
=
V_{\mathcal S}^r q^r
+
V_{\mathcal S}^c\mathcal N_{\theta_u}^u(q^r,\eta_\mu),
\]
with
\[
V_{\mathcal S}^r:=P_{\mathcal S}^T V^r\in\mathbb{R}^{c_{\mathrm{hr}}\times n_r},
\qquad
V_{\mathcal S}^c:=P_{\mathcal S}^T V^c\in\mathbb{R}^{c_{\mathrm{hr}}\times n_c}.
\]
Since $f$ is sample-evaluable, the nonlinear term is evaluated only at the selected entries:
\[
P_{\mathcal S}^T f\bigl(\mathcal U_\mu(q)\bigr)
=
f(U_{\mathcal S})
\in\mathbb{R}^{c_{\mathrm{hr}}},
\]
where the equality follows from the pointwise structure $f(U)_i=f_i(U_i)$.

The cubature approximation of the projected nonlinear term is then
\begin{equation}
\label{eq:ecm_online_main}
\Phi_\mu(q)^T f\bigl(\mathcal U_\mu(q)\bigr)
\approx
\Phi_\mu(q)^T P_{\mathcal S}\operatorname{diag}(\omega)P_{\mathcal S}^T
f\bigl(\mathcal U_\mu(q)\bigr)
=
\bigl(\Phi_\mu(q)^T\bigr)_{(:,\mathcal S)}
\bigl(\omega\odot f(U_{\mathcal S})\bigr)
\in\mathbb{R}^n,
\end{equation}
where $\odot$ denotes the componentwise product.

In the parameter-independent case, the sampled projection matrix can be precomputed offline. In particular, for the linear ROM, if $\Phi_\mu(q)=\Phi$ is independent of both $q$ and $\mu$, then
\[
\bigl(\Phi^T\bigr)_{(:,\mathcal S)}\in\mathbb{R}^{n\times c_{\mathrm{hr}}}
\]
is assembled once and reused throughout the online stage. The online cost of \eqref{eq:ecm_online_main} therefore scales with $c_{\mathrm{hr}}$ and $n$, but not with $N$.

In the NN-ROM, $\Phi_\mu(q^r)$ may depend on $q^r$ through the tangent matrix. For example, for the cubic nonlinearity in the obstacle problem,
\[
\Phi_\mu(q^r)^T=T_{\mathcal U,\mu}(q^r)^T M.
\]
Therefore,
\begin{equation}
\label{eq:ecm_online_TtM}
\begin{aligned}
\bigl(\Phi_\mu(q^r)^T\bigr)_{(:,\mathcal S)}
&=
T_{\mathcal U,\mu}(q^r)^T M_{(:,\mathcal S)} \\
&=
\bigl(V^r+V^c J_{\mathcal N_{\theta_u}^u}(q^r,\eta_\mu)\bigr)^T M_{(:,\mathcal S)} \\
&=
(V^r)^T M_{(:,\mathcal S)}
+
J_{\mathcal N_{\theta_u}^u}(q^r,\eta_\mu)^T
(V^c)^T M_{(:,\mathcal S)}.
\end{aligned}
\end{equation}
The matrices $(V^r)^T M_{(:,\mathcal S)}\in\mathbb{R}^{n_r\times c_{\mathrm{hr}}}$ and $(V^c)^T M_{(:,\mathcal S)}\in\mathbb{R}^{n_c\times c_{\mathrm{hr}}}$ are precomputed offline. Online, one only evaluates the neural network and its Jacobian with respect to the reduced coordinate $q^r$, with $\eta_\mu$ fixed.

Thus, the online evaluation of the hyper-reduced projected nonlinearity depends on the reduced dimensions, the number of cubature points $c_{\mathrm{hr}}$, and the cost of evaluating the neural network and its Jacobian, but not on the HDM dimension $N$.

\subsubsection{Interaction with the NN-augmented ROM}
\label{sec:hr_nn_compatibility}

The distinction between interpolation-based hyper-reduction and the present cubature approach is particularly important for the NN-ROM. Linear interpolation-based methods, such as DEIM or Q-DEIM, approximate the nonlinear term $f(\mathcal U_\mu(q))$ in a linear collateral space in $\mathbb R^N$ before projection. Their accuracy therefore depends on the ability of this collateral space to represent the nonlinear terms encountered during the reduced solves. In the linear ROM, these terms are evaluated at reconstructions of the form $V^{\mathrm{tot}}q^{\mathrm{tot}}$, whereas in the NN-ROM they are evaluated at reconstructions of the form $V^r q^r+V^c\mathcal N_{\theta_u}^u(q^r,\eta_\mu)$. Although the NN-ROM reduces the number of online coordinates from $n_{\mathrm{tot}}$ to $n_r$, DEIM or Q-DEIM would still approximate nonlinear terms evaluated at full-dimensional reconstructions in $\mathbb R^N$. Thus, this coordinate reduction does not by itself imply that the nonlinear terms are more easily approximated by a low-dimensional interpolation basis.

By contrast, the greedy-NNLS cubature rule is fitted directly to projected quantities such as
\[
\Phi_\mu(q^r)^T f(\mathcal U_\mu(q^r))\in\mathbb R^{n_r}.
\]
The fitting problem is posed at the level of the reduced equations, rather than for the full nonlinear term in the HDM space. Consequently, when the NN-ROM replaces a linear ROM of dimension $n_{\mathrm{tot}}$ by a nonlinear-manifold ROM with $n_r$ retained coordinates, the cubature construction is also performed in a lower-dimensional projected space. While this does not guarantee a smaller cubature rule, it is consistent with the reduction in the number of cubature points observed in the numerical experiments.

\begin{rmk}
One could also consider applying hyper-reduction interpolation methods to the projected quantity itself. However, this would generally not lead to an online-efficient procedure: each component of $\Phi_\mu(q)^T f(\mathcal U_\mu(q))$ is a sum over HDM entries, so evaluating even a few components of this projected vector still requires HDM-scale operations unless additional structure is available.
\end{rmk}
\section{Numerical results}
\label{sec:numerical_results}

\paragraph{Computational environment and parallelization.}

All online timing experiments were performed on a Linux workstation equipped with two Intel Xeon E5-2667 v2 processors at $3.30\,\mathrm{GHz}$, providing 16 physical cores and 32 hardware threads in total, and with $188.8\,\mathrm{GiB}$ of RAM. The implementations use Python~3.13.9. The HDMs for both test cases were implemented in FEniCS, while the reduced-order models were implemented primarily using NumPy and SciPy.

The offline HDM snapshot generation was parallelized over parameter values using Python multiprocessing with 32 independent worker processes on the same shared-memory node for the two-dimensional obstacle problem and 24 independent worker processes for the three-dimensional frictional-contact problem. This difference is due to the larger memory footprint of the three-dimensional HDM, which limited the number of concurrent processes. Each worker handled one parameter instance at a time, and no MPI parallelism was used. The offline generation of the low-fidelity snapshots required for neural-network training and of the reduced states required for cubature training was likewise parallelized over parameter values using Python multiprocessing.

The neural networks were trained offline using JAX and JAXopt on an NVIDIA Tesla P100 GPU with $16\,\mathrm{GB}$ of memory. This GPU acceleration is restricted to the offline training stage and is entirely decoupled from the online evaluation of the reduced-order models. During online evaluation, all operations, including neural-network inference, are performed on the CPU.

All reported online timings were obtained sequentially, using a single process and a single numerical-library thread.

  \subsection{Two-dimensional nonlinear obstacle problem}
  \label{sec:results_obstacle}
  \subsubsection{Problem setup}
  \label{sec:obstacle_hdm_snapshots}
    We consider the two-dimensional nonlinear obstacle problem introduced in Section~\ref{sec:2D_obstacle_problem}. In the numerical experiments, the computational domain is
    \[\Omega=(-1,1)\times(-1,1).\]
    On its boundary, we impose the constant Dirichlet condition
    \[u_{\mu}=6 \qquad \text{on } \partial\Omega. \]
    The obstacle is chosen as a curved, anisotropic Gaussian profile, yielding banana-shaped configurations.
    More precisely, for $x=(x_1,x_2)\in\Omega$, we define
    \[\psi_{\mu}(x)=A\,\frac{\widehat{\psi}_{\mu}(x)}
    {\displaystyle
     \max_{y\in\mathcal{N}_h}\widehat{\psi}_{\mu}(y)},\]
    where $\mathcal{N}_h$ denotes the set of degrees of freedom and
    \[\widehat{\psi}_{\mu}(x)=\exp\left[
    -\frac{1}{2}
    \left(\left(\frac{s_{\mu}(x)}{\alpha}\right)^{2}+\left(t_{\mu}(x)-\kappa s_{\mu}(x)^{2}\right)^{2}\right)\right].\]
    The coordinates are given by
    \[s_{\mu}(x)=\cos(\theta)\frac{x_1-c_x}{w}+\sin(\theta)\frac{x_2-c_y}{w},\]
    and
    \[t_{\mu}(x)=-\sin(\theta)\frac{x_1-c_x}{w}+\cos(\theta)\frac{x_2-c_y}{w}.\]
  
  Here, $(c_x,c_y)$ controls the obstacle center, $\theta$ its orientation, $\alpha$ its longitudinal aspect ratio, and $\kappa$ its curvature. The amplitude and width are fixed in all experiments to
  \[
  A=5.5,\qquad w=0.3.
  \]
  The sampled parameters are the five geometric parameters $(c_x,c_y,\theta,\alpha,\kappa)$, together with the nonlinear coefficient $\gamma$, whose ranges are reported in Table~\ref{tab:obstacle_parameters}.  In practice, the nonlinear coefficient is sampled through a normalized parameter $\widehat{\gamma}\in[0,1]$, which is mapped logarithmically to $\gamma$.
  
\begin{table}[pos=h!]
  \centering
  \begin{tabular}{lll}
    \toprule
    \textbf{Parameter} & \textbf{Range} & \textbf{Description} \\
    \midrule
    $c_x$, $c_y$ & $[-0.5,0.5]$ & coordinates of the obstacle center \\
    $\theta$ & $[-\pi,\pi]$ & orientation angle \\
    $\alpha$ & $[1.0,1.5]$ & longitudinal aspect ratio \\
    $\kappa$ & $[0,0.5]$ & curvature strength \\
    $\gamma$ & $[0.1,1]$ & nonlinear coefficient \\
    \bottomrule
  \end{tabular}
  \caption{Parameter description and ranges for the two-dimensional nonlinear obstacle problem with banana-shaped obstacle.}
  \label{tab:obstacle_parameters}
\end{table}

\begin{figure}[pos=t!]
    \centering
    \includegraphics[width=\textwidth]{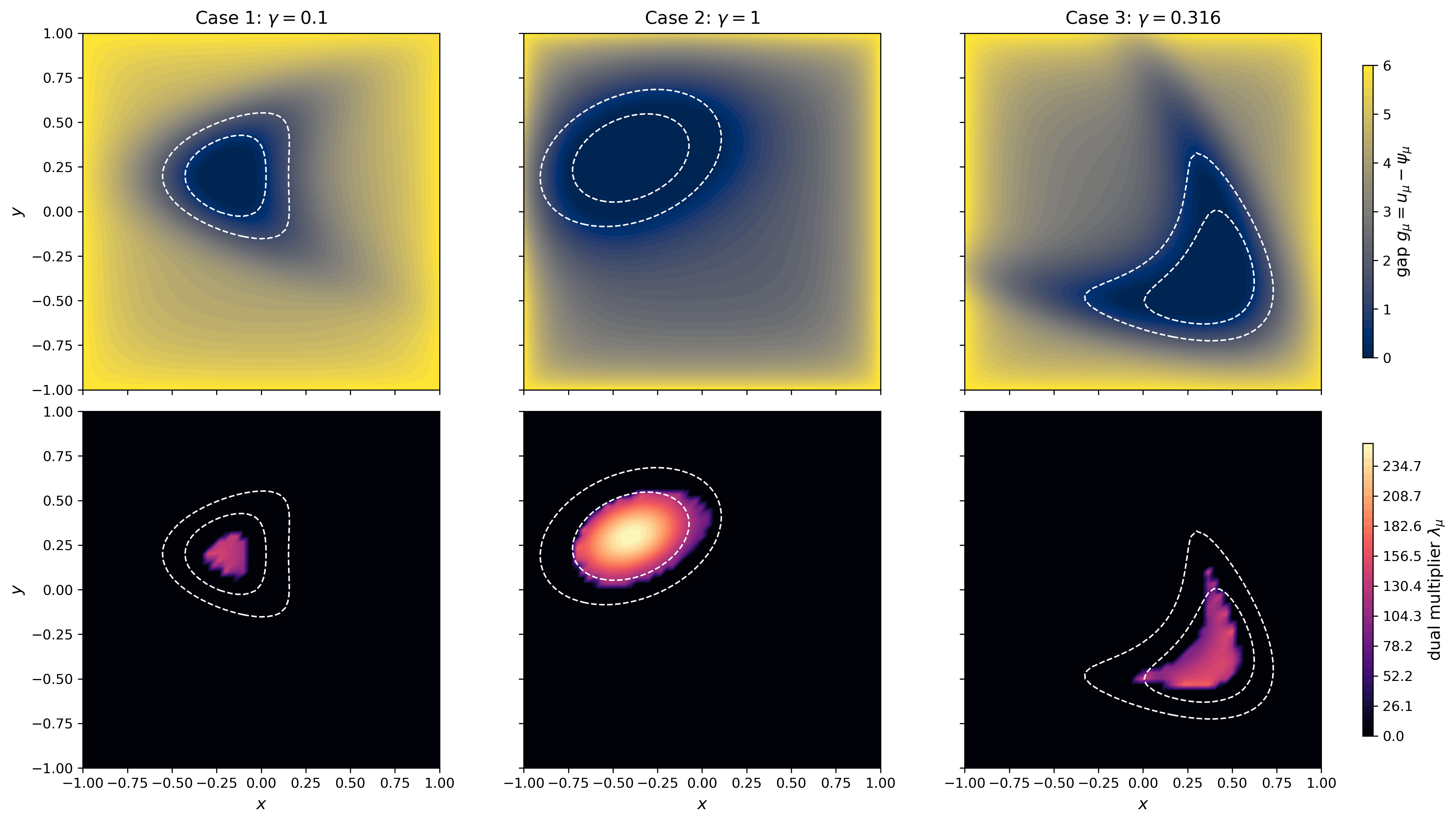}
    \caption{Representative obstacle-contact solutions for three parameter instances. The top row shows the separation gap $g_\mu=u_\mu-\psi_\mu$, and the bottom row shows the corresponding dual multiplier $\lambda_\mu$. White dashed contours indicate level sets of the prescribed obstacle.}
    \label{fig:obstacle_examples}
\end{figure}

\subsubsection{HDM snapshot generation}
  The HDM is discretized using $\mathbb P_2$ Lagrange finite elements on a structured $140\times140$ mesh. The same finite element space is used for the primal and dual variables, leading to
  \[ N=78961,\qquad M=78961, \]
  and therefore to a coupled primal-dual system with $N+M=157922$ degrees of freedom. For each parameter value, the nonlinear nonsmooth algebraic system is solved by the semi-smooth Newton method described in Section~\ref{sec:ssn} with $\rho = 1$. The line-search strategy and HDM stopping criteria are detailed in Appendix~\ref{app:ssn_numerical_settings}.

  Figure~\ref{fig:obstacle_examples} shows representative parameter instances over the prescribed ranges. For each representative parameter instance, we report the discrete separation gap $g_{\mu}=u_{\mu}-\psi_{\mu}$, together with level-set contours of the prescribed obstacle $\psi_{\mu}$.
    
All parameter samples are generated using Sobol sequences in the unit hypercube and are then mapped to the physical ranges reported in Table~\ref{tab:obstacle_parameters}. Independent validation and test sets of $200$ samples each are generated separately from the training set. The validation set is used for all model-selection decisions, including the snapshot budget, reduced dimensions, neural-network architecture, and cubature dimensions, whereas the test set is reserved exclusively for the final numerical comparisons.

  For the training set, we first generate an ordered Sobol sample of size $2^{12}=4096$. Smaller training sets are then obtained by taking the first samples of this list, so that increasing the training budget only requires adding new HDM snapshots without discarding the previous ones. Starting from $250$ snapshots, the training budget is increased progressively. For each budget, the primal POD basis is recomputed and assessed on the validation set through its projection energy error. The snapshot budget is increased until this error reaches the saturation criterion detailed in Appendix~\ref{app:snapshot_budget}. The criterion detects saturation at $N_{\mathrm{snap}}=1500$, which we therefore retain for the construction of the POD bases.

\subsubsection{POD compressibility and linear ROM performance}
\label{sec:obstacle_linear_rom}

\begin{figure}[pos = H]
  \centering
  \includegraphics[width=0.5\textwidth]{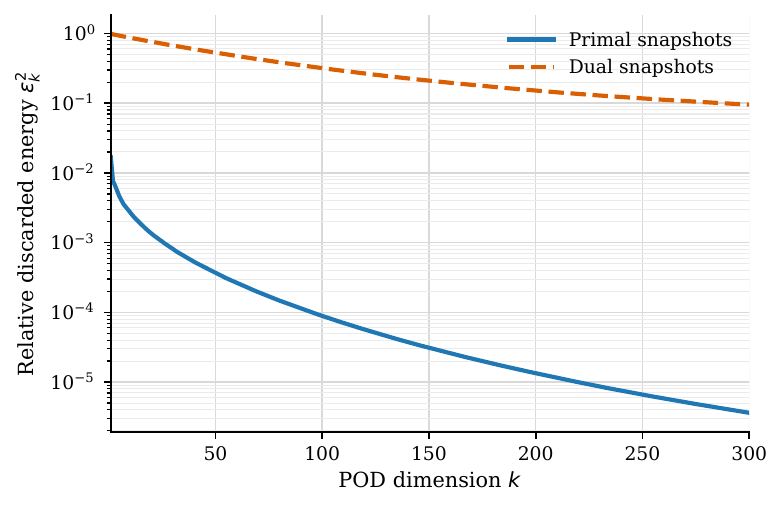}
  \caption{Relative discarded POD energy for the primal and dual snapshot matrices.}
  \label{fig:obstacle_svd}
\end{figure}

We first examine the POD spectra associated with the primal and dual snapshot matrices, denoted by $S_U$ and $S_\Lambda$. Figure~\ref{fig:obstacle_svd} shows the relative discarded energy after retaining the first $k$ POD modes,
\[
\varepsilon_k^2
=
1-\frac{\sum_{i=1}^k \sigma_i^2}{\sum_{i=1}^{N_s}\sigma_i^2},
\]
where $\sigma_i$ are the singular values of the corresponding snapshot matrix. The dual spectrum decays much more slowly than the primal one, showing that the multiplier snapshots are less compressible by a linear basis. This is consistent with the localized nature of the multiplier and with its sensitivity to changes in the active set.

We first evaluate the Galerkin ROM for several primal and dual reduced dimensions $(n,m)$ on the validation set. For each validation parameter, the reduced nonlinear nonsmooth system is solved by the SSN method and the reconstructed primal solution is compared with the corresponding HDM reference.

For the validation comparisons reported below, the relative primal energy error is defined, for each parameter $\mu$, by
\[
e_u(\mu)=
\frac{\left\|U^{\mathrm{red}}_{\mu}-U_{\mu}\right\|_{K}}
{\left\|U_{\mu}\right\|_{K}},
\qquad
\lVert z\rVert_K:=\sqrt{z^{T}Kz}.
\]
The values shown in Figure~\ref{fig:obstacle_rom_heatmap} are the arithmetic means of $e_u(\mu)$ over the validation set. The online speed-up is defined by
$S=\frac{\overline{t}_{\mathrm{HDM}}}{\overline{t}_{\mathrm{red}}},$
where both times are averaged over the same parameter set.
\begin{figure}[pos=H]
  \centering
  \resizebox{0.6\linewidth}{!}{
    \begin{tikzpicture}
\begin{axis}[
    width=9cm,
    height=9cm,
    scale only axis,
    xlabel={Number of primal modes $n$},
    ylabel={Number of dual modes $m$},
    xlabel style={font=\small},
    ylabel style={font=\small},
    tick label style={font=\small},
    xtick={1,2,3,4,5},
    xticklabels={25,50,100,150,200},
    ytick={1,2,3,4,5},
    yticklabels={25,50,100,150,200},
    xmin=0.5, xmax=5.5,
    ymin=0.5, ymax=5.5,
    enlargelimits=false,
    axis on top,
    point meta min=1.95,
    point meta max=8.85,
    colorbar,
    colorbar style={
        ylabel={Mean relative energy error (\%)},
        ylabel style={
            font=\small,
            rotate=0,
            at={(axis description cs:6.0,0.5)},
            anchor=south,
        },
        tick label style={font=\small},
    },
    colormap={paperfriendly}{
        rgb255=(250,248,225)
        rgb255=(235,235,180)
        rgb255=(205,215,145)
        rgb255=(160,190,145)
        rgb255=(115,160,150)
        rgb255=(80,125,155)
    },
]

\addplot[
    matrix plot*,
    mesh/rows=5,
    draw=none,
    point meta=explicit,
] coordinates {
    (1,1) [8.13]
    (2,1) [7.61]
    (3,1) [7.57]
    (4,1) [7.57]
    (5,1) [7.57]

    (1,2) [8.38]
    (2,2) [4.90]
    (3,2) [4.50]
    (4,2) [4.49]
    (5,2) [4.49]

    (1,3) [8.65]
    (2,3) [5.01]
    (3,3) [2.98]
    (4,3) [2.80]
    (5,3) [2.80]

    (1,4) [8.77]
    (2,4) [5.11]
    (3,4) [2.91]
    (4,4) [2.32]
    (5,4) [2.21]

    (1,5) [8.85]
    (2,5) [5.16]
    (3,5) [2.93]
    (4,5) [2.24]
    (5,5) [1.95]
};


\node[font=\small, align=center] at (axis cs:1,1)
    {8.13\%\\($\times$ 56.28)};
\node[font=\small, align=center] at (axis cs:2,1)
    {7.61\%\\($\times$ 31.03)};
\node[font=\small, align=center] at (axis cs:3,1)
    {7.57\%\\($\times$ 12.14)};
\node[font=\small, align=center] at (axis cs:4,1)
    {7.57\%\\($\times$ 8.08)};
\node[font=\small, align=center] at (axis cs:5,1)
    {7.57\%\\($\times$ 5.21)};

\node[font=\small, align=center] at (axis cs:1,2)
    {8.38\%\\($\times$ 47.10)};
\node[font=\small, align=center] at (axis cs:2,2)
    {4.90\%\\($\times$ 22.41)};
\node[font=\small, align=center] at (axis cs:3,2)
    {4.50\%\\($\times$ 12.96)};
\node[font=\small, align=center] at (axis cs:4,2)
    {4.49\%\\($\times$ 7.83)};
\node[font=\small, align=center] at (axis cs:5,2)
    {4.49\%\\($\times$ 5.13)};

\node[font=\small, align=center] at (axis cs:1,3)
    {8.65\%\\($\times$ 30.95)};
\node[font=\small, align=center] at (axis cs:2,3)
    {5.01\%\\($\times$ 15.93)};
\node[font=\small, align=center] at (axis cs:3,3)
    {2.98\%\\($\times$ 9.00)};
\node[font=\small, align=center] at (axis cs:4,3)
    {2.80\%\\($\times$ 6.83)};
\node[font=\small, align=center] at (axis cs:5,3)
    {2.80\%\\($\times$ 4.58)};

\node[font=\small, align=center] at (axis cs:1,4)
    {8.77\%\\($\times$ 23.13)};
\node[font=\small, align=center] at (axis cs:2,4)
    {5.11\%\\($\times$ 11.05)};
\node[font=\small, align=center] at (axis cs:3,4)
    {2.91\%\\($\times$ 7.57)};
\node[font=\small, align=center] at (axis cs:4,4)
    {2.32\%\\($\times$ 5.02)};
\node[font=\small, align=center] at (axis cs:5,4)
    {2.21\%\\($\times$ 4.12)};

\node[font=\small, align=center] at (axis cs:1,5)
    {8.85\%\\($\times$ 17.54)};
\node[font=\small, align=center] at (axis cs:2,5)
    {5.16\%\\($\times$ 8.16)};
\node[font=\small, align=center] at (axis cs:3,5)
    {2.93\%\\($\times$ 6.55)};
\node[font=\small, align=center] at (axis cs:4,5)
    {2.24\%\\($\times$ 4.55)};
\node[font=\small, align=center] at (axis cs:5,5)
    {1.95\%\\($\times$ 3.35)};

\end{axis}
\end{tikzpicture}}
  \caption{Mean relative energy error of the Galerkin ROM for different primal and dual reduced dimensions on the validation set. The values in parentheses indicate the online speed-up with respect to the HDM.}
  \label{fig:obstacle_rom_heatmap}
\end{figure}

\begin{rmk}
The POD spectra in Figure~\ref{fig:obstacle_svd} suggest that the primal space is significantly more compressible than the dual space. However, the mean energy error of the coupled ROM in Figure~\ref{fig:obstacle_rom_heatmap} shows a more balanced dependence on the primal and dual dimensions. This is because the two figures measure different quantities: the POD spectra characterize the separate projection compressibility of the primal and dual snapshots, whereas the heatmap reports the primal energy error of the \emph{coupled} nonlinear ROM. Since the reconstructed multiplier enters directly into the primal equilibrium through the primal--dual coupling term, an under-resolved dual space limits the achievable displacement accuracy even when the primal space is enriched. The error therefore decreases mainly when the two reduced spaces are enlarged together.
\end{rmk}
\subsubsection{NN-ROM training and model selection}
\label{sec:nn_rom_2D_obstacle}
Consistently with the offline construction described in Section~\ref{sec:nn_offline}, the NN-ROM training data are not obtained from additional HDM solves. Once the POD bases have been fixed, we generate an additional independent set of parameter values using a Sobol sequence. For each of these parameters, we solve the non-hyper-reduced Galerkin ROM with $n_{\mathrm{tot}}=m_{\mathrm{tot}}=150$ and use the resulting primal and dual reduced coefficients as training data for the neural networks. The NN-ROM training set consists of $4500$ low-fidelity samples.

The network architecture is selected by a validation-based hyperparameter search. For this search, we fix the enlarged POD dimensions to $n_{\mathrm{tot}}=m_{\mathrm{tot}}=150$ and the retained dimensions to $n_r=m_r=25$. The parameter-feature vector is chosen as
\[\eta_\mu = \bigl( \widehat{\gamma},\ c_x,\ c_y,\ \cos\theta,\ \sin\theta,\ \alpha,\ \kappa \bigr),\]
where the orientation angle is encoded through its sine and cosine to avoid the artificial discontinuity at $\theta=\pm\pi$.
The input dimension of each network is therefore $25+7$, corresponding to the retained reduced coordinates and the parameter features, while the output dimension is $125$. The neural-network architecture is selected using the validation set according to the procedure detailed in Appendix~\ref{app:nn_architecture_selection}. For the two-dimensional problem, the networks are trained using $4500$ low-fidelity ROM samples, and the retained architecture uses SiLU activations and three hidden layers of widths $256$, $512$, and $256$, providing a good compromise between validation accuracy and model size.

Since the numbers of retained coordinates $n_r$ and $m_r$ directly affect both the information available to the neural correction and the dimension of the online NN-ROM system, we also examine their influence in Appendix~\ref{app:nn_input_dimension_sensitivity}.

\subsubsection{Hyper-reduction configuration}
\label{sec:obstacle_hyper_reduction}

For both the ROM and the NN-ROM, hyper-reduction is applied to the remaining full-order contributions in the online reduced solve. These are the cubic nonlinear term and the complementarity projection. Separate sparse cubature rules are therefore constructed for these two contributions, following the procedure described in Section~\ref{sec:hyper_reduction}.

In the offline phase, we set the cubature tolerance to $\varepsilon_\mathrm{cub}^\mathrm{cubic}=10^{-2}$ for the cubic nonlinear term and to $\varepsilon_{\mathrm{cub}}^\mathrm{comp}=2\times 10^{-2}$ for the complementarity projection. We impose a maximum number of cubature points $c_{\max}=1250$. The resulting NN-HR-ROM cubature rules contain $221$ points for the cubic nonlinear contribution and $763$ points for the complementarity projection. For the HR-ROM, the corresponding rules contain $255$ and $1236$ points, respectively. The cubature error curves are reported in Appendix~\ref{app:cubature_convergence}.

\subsubsection{Results}
Table~\ref{tab:ROM_results_2D} summarizes the performance of the different models for the two-dimensional obstacle problem. The final assessment is performed on the independent test set of $N_{\mathrm{test}}=200$ parameter values, using the same error and speed-up definitions.

The standard ROM, with $150$ primal and dual coordinates, achieves the lowest mean relative energy error, $2.34\%$, with a speed-up factor of $5.35$ relative to the HDM. The NN-ROM reduces the number of coordinates to $25$ for each variable and gives a mean error of $2.70\%$, with a speed-up factor of $6.88$. The moderate computational gain of the NN-ROM is expected because its residuals are still assembled at the full-order level.

Hyper-reduction provides substantially larger gains while introducing only a small additional error. The HR-ROM achieves a mean error of $2.44\%$ with a speed-up factor of $236$. The proposed NN-HR-ROM is the fastest configuration, with a speed-up factor of $619$, while retaining a mean relative energy error of $2.77\%$. It uses $25$ primal and dual coordinates together with $221$ and $763$ cubature points for the cubic state and projected complementarity contributions, respectively.

\paragraph{Constraint feasibility and mechanical response.} 
The primal energy error does not by itself assess whether the reconstructed primal and dual fields satisfy the obstacle conditions. We therefore perform the full-space audit detailed in Appendix~\ref{app:constraint_diagnostics}. The corresponding results are reported in Table~\ref{tab:obstacle_constraint_analysis}. Across the four models, the mean residual penetration ranges from $0.164\%$ to $0.225\%$, while the mean obstacle-constraint residual remains between $0.167\%$ and $0.203\%$. The main feasibility limitation concerns multiplier nonnegativity, for which the mean indicator ranges from $6.544\%$ to $7.413\%$. The Jacobian-based mechanical-response error remains between $2.703\%$ and $3.269\%$. Most of the increase in the negative-multiplier indicator is associated with the NN augmentation, whereas hyper-reduction produces the systematic increase in the mechanical-response error. Overall, the primal response remains only moderately affected.
{ \scriptsize
  \begin{table}[pos=H]
    \centering
    \setlength{\tabcolsep}{3pt}
    \renewcommand{\arraystretch}{1.8}
    \begin{tabular}{|l||c|c||c|c||c||c||c|}
      \hline
      \multirow{2}{*}{\textbf{Method}} &
        \multicolumn{2}{c||}{\textbf{Dimension}} &
        \multicolumn{2}{c||}{\makecell{\textbf{Cubature points}}} &
        \multirow{2}{*}{\makecell{\textbf{Mean relative}\\\textbf{energy error}}} &
        \multirow{2}{*}{\makecell{\textbf{Mean}\\\textbf{solve}\\\textbf{time (s)}}} &
        \multirow{2}{*}{\makecell{\textbf{Speed--up}\\\textbf{factor}\\\textbf{vs HDM}}} \\
      \cline{2-5}
        & \textbf{Primal}
        & \textbf{Dual}
        & \textbf{Cubic}
        & \textbf{Projection}
        & & & \\
      \hline

      \textbf{HDM}
        & $78{,}961$
        & $78{,}961$
        & --
        & --
        & --
        & $130$
        & -- \\ \hline

      \textbf{ROM}
        & $150$
        & $150$
        & --
        & --
        & $2.34\%$
        & $24.3$
        & \textbf{5.35} \\ \hline

      \textbf{NN-ROM}
        & $25$
        & $25$
        & --
        & --
        & $2.70\%$
        & $18.9$
        & \textbf{6.88} \\ \hline

      \textbf{HR-ROM}
        & $150$
        & $150$
        & $255$
        & $1236$
        & $2.44\%$
        & $0.55$
        & \textbf{236} \\ \hline

      \textbf{NN-HR-ROM}
        & $25$
        & $25$
        & $221$
        & $763$
        & $2.77\%$
        & $0.21$
        & \textbf{619} \\ \hline

    \end{tabular}
    \caption{Performance of the HDM and reduced models on the independent test set. The reported error is the mean relative primal energy error and the computation time is the mean online solution time. The ``Cubic'' and ``Projection'' columns give the numbers of points selected for the cubic nonlinear term and the complementarity projection, respectively. Offline costs are not included in the reported online times.}
    \label{tab:ROM_results_2D}
  \end{table}
}
\paragraph{Offline data-generation and learning costs.} For the two-dimensional problem, the main offline cost is the generation of the $1500$ HDM snapshots, which required $5.01\,\mathrm{h}$. The additional low-fidelity database used for the NN-ROM required $3.09\,\mathrm{h}$ for $5000$ ROM solves, so its cost is non-negligible relative to the high-fidelity snapshot campaign. On the other hand, training the selected neural networks took $8.8\,\mathrm{min}$. Thus, the NN-specific offline overhead is driven almost entirely by low-fidelity data generation rather than by the neural-network training itself.

\subsection{Three-dimensional frictional contact problem}
\label{sec:results_contact}

\subsubsection{Problem setup}
\label{sec:contact_hdm_snapshots}

\begin{figure}[pos=h!]
    \centering
    \includegraphics[width=0.6\textwidth]{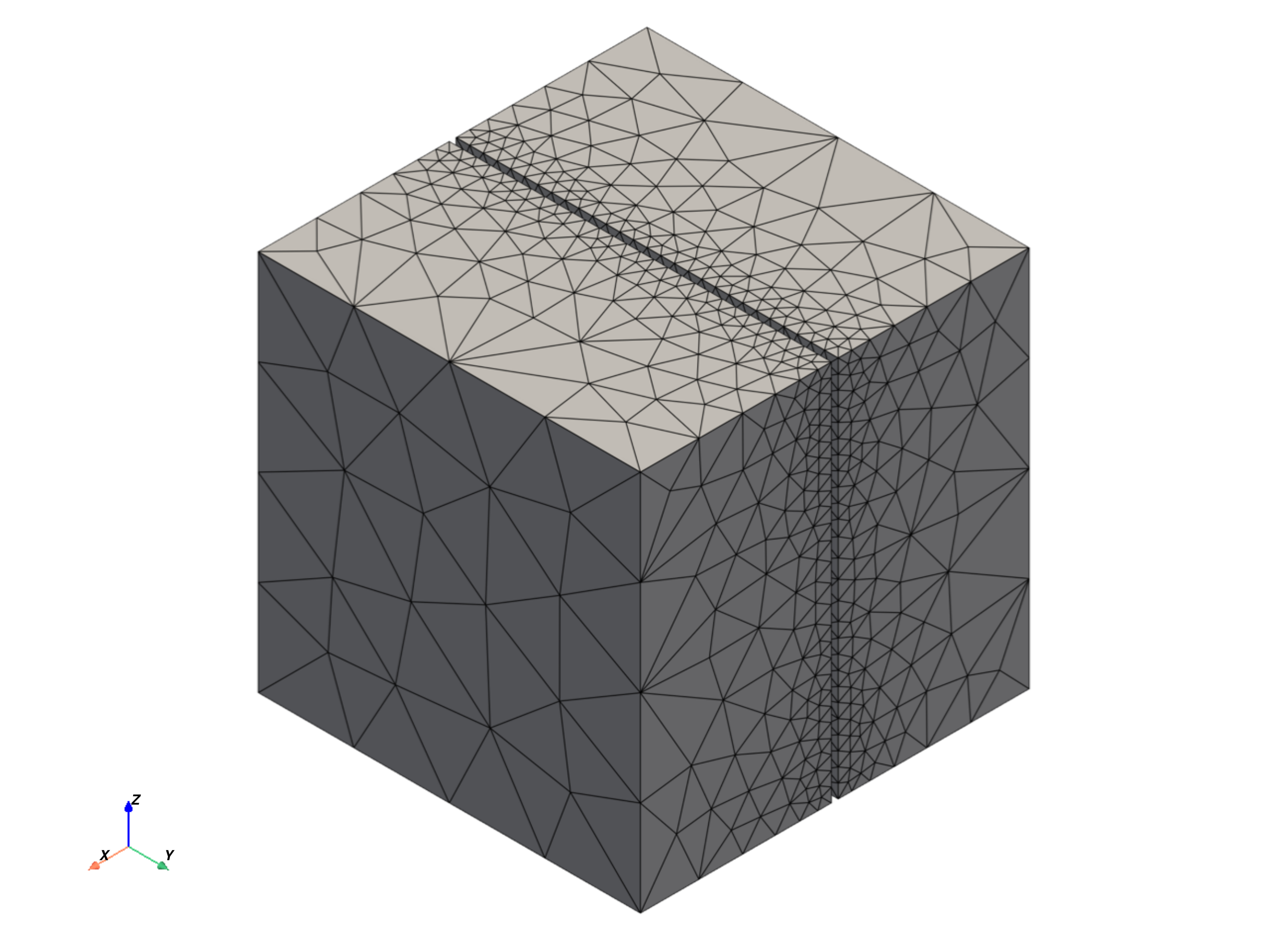}
    \caption{Conforming tetrahedral mesh used for the three-dimensional frictional contact problem. The mesh is locally refined near the internal interface $\Gamma=\{0\}\times(-1,1)^2$ and gradually coarsened away from the interface. The two subdomains are slightly separated in the visualization only.}
    \label{fig:contact_mesh}
\end{figure}

We now consider the three-dimensional frictional contact problem introduced in Section~\ref{sec:3D_contact_problem}. The computational domain and fracture interface are
\[
\Omega=(-1,1)^3,
\qquad
\Gamma=\{0\}\times(-1,1)^2,
\]
which separates the two elastic bodies $\Omega^-$ and $\Omega^+$. Homogeneous displacement boundary conditions are imposed on the exterior boundary $\partial\Omega$. The domain is discretized by a conforming tetrahedral mesh matching the internal interface $\Gamma$. The mesh is locally refined near $\Gamma$ and gradually coarsened away from the interface, with target mesh sizes $h_\Gamma=1/14$ close to the fracture and $h_{\mathrm{far}}=1/2$ near the outer boundary. The resulting mesh is shown in Figure~\ref{fig:contact_mesh}. The displacement field is approximated with vector-valued $\mathbb P_2$ finite elements on the two subdomains, while the normal and tangential contact multipliers are discretized by face-wise constant functions on $\Gamma$.

The two bodies are assigned different Lamé parameters:
\[
\mu_L^-=\lambda_L^-=0.08,
\qquad
\mu_L^+=\lambda_L^+=0.80.
\]
These material parameters remain fixed throughout the sampling campaign. 

For a center $c=(c_1,c_2,c_3)\in\mathbb R^3$, we define
\[ \zeta_c(x) = 
  \begin{cases}
    \displaystyle \exp\left(-\frac{1}{1-d_c(x)^2}+1\right), & d_c(x)<1,\\[2mm] 0, & d_c(x)\geq1,
  \end{cases}
\qquad d_c(x)=\frac{\|x-c\|_2}{\ell_f}, \]

where $\ell_f>0$ controls the width of the profile.

The force-patch centers have fixed positions in the direction normal to the interface and parameter-dependent transverse coordinates:
\[ c_\mu^-=(x_c^-,Y_{CL},Z_{CL}), \qquad c_\mu^+=(x_c^+,Y_{CR},Z_{CR}). \]
We parameterize the problem by
\[
\mu=(Y_{CL},Z_{CL},Y_{CR},Z_{CR},F)\in\mathcal P,
\]
where $F$ denotes the Coulomb friction coefficient entering the tangential projection operator. The parameter domain $\mathcal P$, defined by the sampling ranges of these five parameters, is reported in Table~\ref{tab:contact_parameters}.

The corresponding body forces are
\[ g_\mu^-(x) = \zeta_{c_\mu^-}(x)
  \begin{pmatrix}
  A_N\\ A_T\\ 0
  \end{pmatrix},
\qquad g_\mu^+(x) = -\zeta_{c_\mu^+}(x)
  \begin{pmatrix}
  A_N\\ A_T\\ 0
  \end{pmatrix}. \]
The fixed quantities are
\[ x_c^-=-0.35, \qquad x_c^+=0.35, \qquad A_N=100, \qquad A_T=20, \qquad \ell_f=0.30. \]
Thus, the two bodies are subjected to opposite force patches whose transverse locations vary independently.

\begin{table}[pos=h!]
  \centering
  \begin{tabular}{lll}
    \toprule
    \textbf{Parameter} & \textbf{Range} & \textbf{Description} \\
    \midrule
    $Y_{CL}$, $Z_{CL}$ & $[-0.65,0.65]$ & transverse coordinates of the force patch in $\Omega^-$ \\
    $Y_{CR}$, $Z_{CR}$ & $[-0.65,0.65]$ & transverse coordinates of the force patch in $\Omega^+$ \\
    $F$ & $[0.15,0.80]$ & Coulomb friction coefficient \\
    \bottomrule
  \end{tabular}
  \caption{Parameter description and ranges for the three-dimensional frictional contact problem.}
  \label{tab:contact_parameters}
\end{table}
\subsubsection{HDM snapshot generation}
For each parameter value, the HDM system is solved by the primal--dual semi-smooth Newton method described in Section~\ref{sec:ssn}. The projection parameter is set to $\rho=10^{-2}$.

\begin{figure}[pos=t!]
    \centering
    \includegraphics[width=0.98\textwidth]{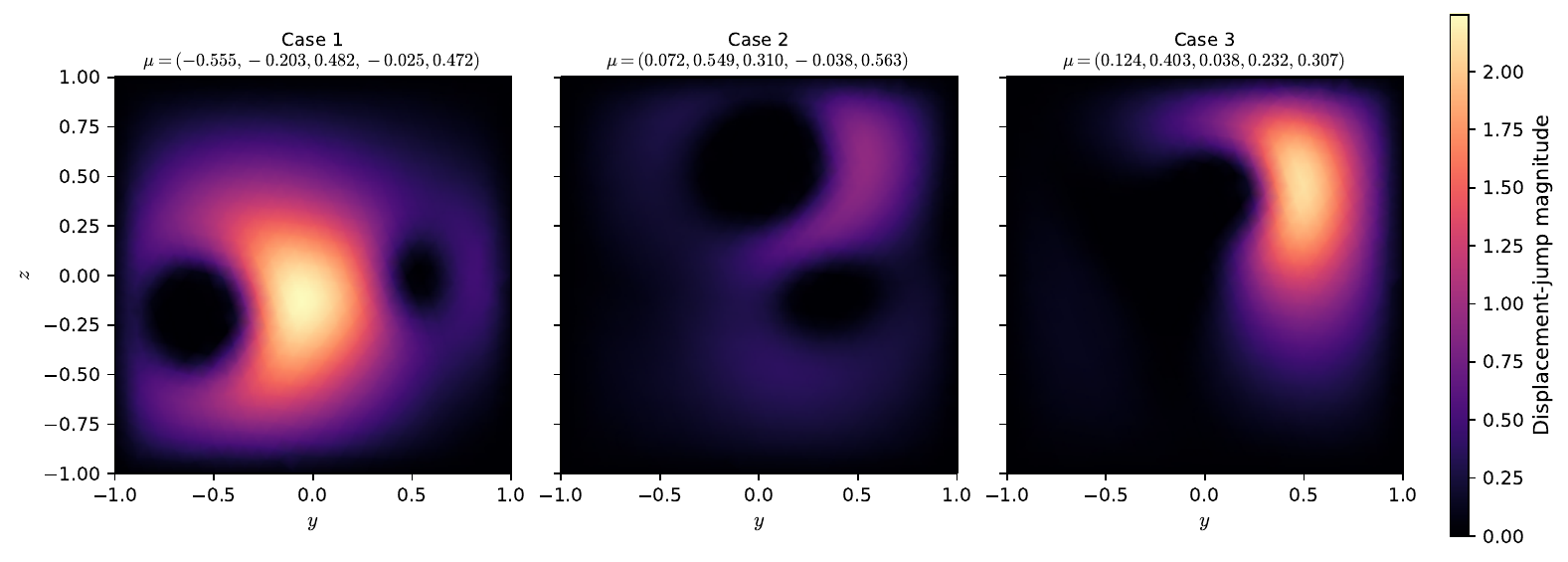}
    \caption{
    Interface displacement-jump magnitude for three HDM snapshots of the three-dimensional frictional contact problem. The complete parameter point is reported below each panel title in the order $\mu=(Y_{CL},Z_{CL},Y_{CR},Z_{CR},F)$. For visualization only, the face-averaged jump is displayed through its continuous piecewise-linear $L^2$-projection on the interface mesh.
    }
    \label{fig:contact_jump_examples}
\end{figure}

The snapshot campaign follows the same Sobol-based sampling strategy as in the two-dimensional test case. We generate $2000$ training samples, together with independent validation and test sets of $200$ samples each. The snapshot-budget criterion described in Appendix~\ref{app:snapshot_budget} again selects $N_{\mathrm{snap}}=1500$ HDM snapshots for the construction of the POD bases. For every parameter value, the HDM solve provides the displacement field and the normal and tangential contact multipliers. 
Figure~\ref{fig:contact_jump_examples} shows the displacement-jump magnitude
$\left\|\llbracket u_\mu\rrbracket\right\|_2$ on the interface $\Gamma$
for three HDM snapshots associated with different parameter values.

These snapshots are used below to construct the POD spaces and to assess the reduced models.

\subsubsection{Correction-based POD spaces}
\label{sec:contact_correction_pod}

Before constructing the reduced spaces, we first adapt the primal variable to the structure of the three-dimensional contact problem. The displacement field contains a dominant unconstrained elastic response induced by the body forces. If the POD basis were built directly from the full displacement snapshots, a significant part of the reduced space would be used to represent this forced linear response, whereas the main nonlinear effect to be captured is the contact-induced correction.

For each parameter value, we therefore compute the unconstrained elastic predictor $Y_\mu$ defined by
\[
K Y_\mu=L_\mu,
\]
and decompose the HDM displacement as
\[
U_\mu = Y_\mu+\delta U_\mu .
\]
The primal POD basis $V$ is constructed from the correction snapshots
\[
\delta U_\mu=U_\mu-Y_\mu,
\]
rather than from the full displacement snapshots. The reduced primal approximation is then written as
\[
U_\mu \approx Y_\mu+Vq_\mu .
\]

The contact conditions are imposed on the total displacement $Y_\mu+Vq_\mu$. The normal and tangential multipliers are represented with separate POD bases,
\[
\Lambda_{\mu,n}\approx W_n \xi_{\mu,n}, \qquad \Lambda_{\mu,\tau}\approx W_\tau \xi_{\mu,\tau}.
\]
This leads to a Galerkin ROM in the reduced variables $(q_\mu,\xi_{\mu,n},\xi_{\mu,\tau}),$ whose shifted form is given in Appendix~\ref{app:contact_shifted_rom} and solved with the semi-smooth Newton linearization. The predictor $Y_\mu$ is not an additional unknown of the reduced problem: after the shift, it only appears through $B_nY_\mu$ and $B_\tau Y_\mu$, which are computed before the semi-smooth Newton iterations.

\begin{figure}[pos = h!]
    \centering
    \includegraphics[width=0.9\textwidth]{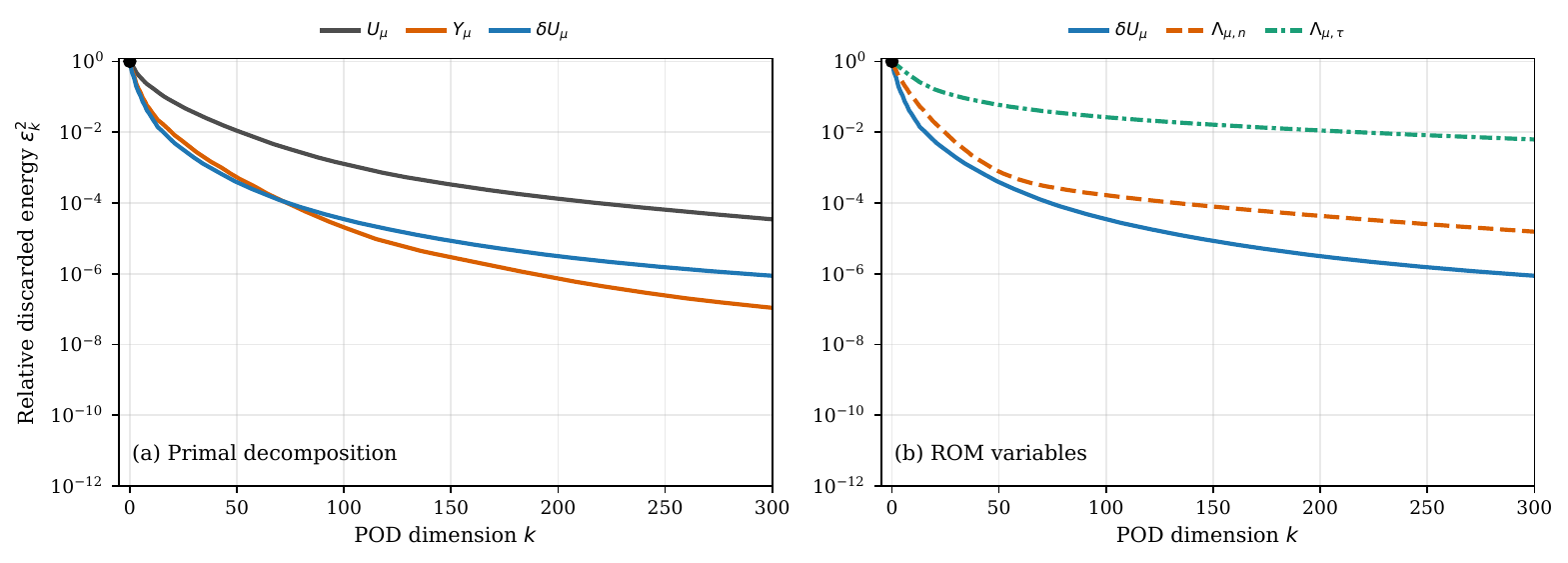}
    \caption{Relative discarded POD energy for the three-dimensional contact problem. The left panel compares the full displacement $U_\mu$, the unconstrained elastic predictor $Y_\mu$, and the correction $\delta U_\mu=U_\mu-Y_\mu$. The faster decay of the correction snapshots motivates the use of $\delta U_\mu$ to build the primal reduced basis. The right panel reports the discarded-energy curves of the variables used in the reduced contact system, namely the primal correction, the normal multiplier, and the tangential multiplier.}
    \label{fig:contact_svd}
\end{figure}

Figure~\ref{fig:contact_svd} reports the POD discarded-energy curves. For the primal variable, we compare the spectra obtained from the full displacement snapshots $U_\mu$, the unconstrained elastic predictors $Y_\mu$, and the correction snapshots $\delta U_\mu$. The correction snapshots exhibit a faster decay, confirming that subtracting the unconstrained elastic predictor removes a dominant linear component and yields a more compact reduced representation of the contact-induced response. The same figure also reports the POD spectra of the normal and tangential multiplier snapshots.

\subsubsection{Online approximation of the elastic predictor}
\label{sec:contact_elastic_predictor_surrogate}

The correction-based ROM requires the unconstrained elastic predictor $Y_\mu$, or at least its interface contributions in the contact equations. In principle, $Y_\mu$ is obtained by solving
\[
K Y_\mu = L_\mu .
\]
Although $K$ is parameter-independent, computing $Y_\mu$ exactly for a new parameter value would still require assembling the parameter-dependent load vector $L_\mu$ on the HDM space and applying a full-order linear solve. Once the contact problem has been reduced, this auxiliary predictor computation becomes the dominant remaining online cost, as reflected by the timing results in Table~\ref{tab:ROM_results_3D}. We therefore replace it by a reduced surrogate.

The surrogate is built at the level of the POD coefficients of the elastic predictor. By linearity, the responses to the two loads can be treated separately:
\[
Y_\mu = Y^-_\mu+Y^+_\mu .
\]
For each contribution, we construct a POD basis from the corresponding predictor snapshots and write
\[
Y^-_\mu \approx V_Y^- a^-_\mu,
\qquad
Y^+_\mu \approx V_Y^+ a^+_\mu .
\]
These predictor snapshots require no separate full-order snapshot-generation campaign. They are generated at negligible additional cost during the HDM computations and saved for later use, so constructing $V_Y^-$ and $V_Y^+$ only requires a POD of already available snapshots.

The POD coefficients $a_\mu^-$ and $a_\mu^+$ are then approximated from the corresponding load parameters by a local radial basis function (RBF) approximation. Denoting by $\chi$ the load-parameter vector of either subdomain, normalized to $[0,1]^2$ over the training set, the approximation at $\chi$ is built from the $k=16$ training points closest to $\chi$, with index set $\mathcal{N}_k(\chi)$:
\[
\widehat{a}(\chi)=\sum_{i\in\mathcal{N}_k(\chi)} \alpha_i(\chi)\,\varphi\bigl(\lVert\chi-\chi_i\rVert\bigr)+p_\chi(\chi),
\]
where $\varphi(r)=r^2\log r$ is the thin-plate spline kernel and $p_\chi$ is a polynomial of degree at most $3$.
This interpolation is used only to approximate the auxiliary elastic predictor.

The resulting online approximation is denoted by
\[
\widehat Y_\mu = V_Y^- \widehat a^-_\mu + V_Y^+ \widehat a^+_\mu .
\]
The surrogate $\widehat Y_\mu$ is used only to approximate the shifted interface traces $B_nY_\mu$ and $B_\tau Y_\mu$, not as a direct substitute in the equilibrium equation. Appendix~\ref{app:contact_shifted_rom} details this choice and reports validation metrics for the predictor, interface traces, and projected equilibrium defect.

\begin{rmk}
A conventional reduced approximation of the elastic predictor would write $Y_\mu\approx V_Ya_\mu$ and determine the coefficients online from
\[
V_Y^TKV_Ya_\mu=V_Y^TL_\mu.
\]
Although the reduced operator can be precomputed, this approach still requires assembling the full-order load vector $L_\mu$ for each new parameter. To address this, we investigated an EIM approximation of $L_\mu$, but the force profiles required too many interpolation modes to achieve sufficient accuracy at a reasonable online cost. This is why we adopted a different strategy and approximated the predictor coefficients directly by RBF interpolation. This has the advantage that the resulting POD--RBF surrogate depends only on the parameter values and avoids full-order load assembly. Its accuracy is assessed in Appendix~\ref{app:contact_shifted_rom}, while its effect on the online cost is reported in Table~\ref{tab:ROM_results_3D}.
\end{rmk}

\subsubsection{POD compressibility and linear ROM performance}
\label{sec:contact_linear_rom}

\begin{figure}[pos=!htbp]
  \centering

  \begin{subfigure}[c]{0.40\textwidth}
    \centering
    \includegraphics[width=\linewidth]{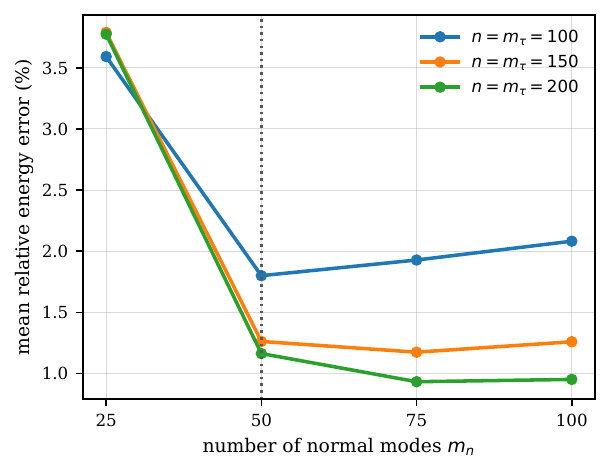}
    \caption{Selection of the normal-dual dimension $m_n$.}
    \label{fig:normal_dimension_selection}
  \end{subfigure}
  \hfill
  \begin{subfigure}[c]{0.56\textwidth}
    \centering
    \resizebox{\linewidth}{!}{
      \begin{tikzpicture}
\begin{axis}[
    width=9cm,
    height=9cm,
    scale only axis,
    xlabel={Number of primal modes $n$},
    ylabel={Number of tangential dual modes $m_\tau$},
    xlabel style={font=\small},
    ylabel style={font=\small},
    tick label style={font=\small},
    xtick={1,2,3,4,5},
    xticklabels={25,50,100,150,200},
    ytick={1,2,3,4,5},
    yticklabels={25,50,100,150,200},
    xmin=0.5, xmax=5.5,
    ymin=0.5, ymax=5.5,
    enlargelimits=false,
    axis on top,
    point meta min=1.16,
    point meta max=9.42,
    colorbar,
    colorbar style={
        ylabel={Mean relative energy error (\%)},
        ylabel style={
            font=\small,
            rotate=0,
            at={(axis description cs:6.0,0.5)},
            anchor=south,
        },
        tick label style={font=\small},
    },
    colormap={paperfriendly}{
        rgb255=(250,248,225)
        rgb255=(235,235,180)
        rgb255=(205,215,145)
        rgb255=(160,190,145)
        rgb255=(115,160,150)
        rgb255=(80,125,155)
    },
]

\addplot[
    matrix plot*,
    mesh/rows=5,
    draw=none,
    point meta=explicit,
] coordinates {
    (1,1) [9.40] (2,1) [4.33] (3,1) [3.12] (4,1) [3.50] (5,1) [3.52]

    (1,2) [9.42] (2,2) [4.23] (3,2) [2.30] (4,2) [2.36] (5,2) [2.38]

    (1,3) [9.41] (2,3) [4.11] (3,3) [1.80] (4,3) [1.53] (5,3) [1.54]

    (1,4) [9.40] (2,4) [4.14] (3,4) [1.70] (4,4) [1.26] (5,4) [1.27]

    (1,5) [9.42] (2,5) [4.16] (3,5) [1.69] (4,5) [1.20] (5,5) [1.16]
};


\node[font=\small, align=center] at (axis cs:1,1) {9.40\%\\($727.29\times$)};
\node[font=\small, align=center] at (axis cs:2,1) {4.33\%\\($703.10\times$)};
\node[font=\small, align=center] at (axis cs:3,1) {3.12\%\\($927.54\times$)};
\node[font=\small, align=center] at (axis cs:4,1) {3.50\%\\($873.41\times$)};
\node[font=\small, align=center] at (axis cs:5,1) {3.52\%\\($828.28\times$)};
\node[font=\small, align=center] at (axis cs:1,2) {9.42\%\\($610.85\times$)};
\node[font=\small, align=center] at (axis cs:2,2) {4.23\%\\($594.79\times$)};
\node[font=\small, align=center] at (axis cs:3,2) {2.30\%\\($693.13\times$)};
\node[font=\small, align=center] at (axis cs:4,2) {2.36\%\\($658.84\times$)};
\node[font=\small, align=center] at (axis cs:5,2) {2.38\%\\($606.76\times$)};
\node[font=\small, align=center] at (axis cs:1,3) {9.41\%\\($432.65\times$)};
\node[font=\small, align=center] at (axis cs:2,3) {4.11\%\\($398.24\times$)};
\node[font=\small, align=center] at (axis cs:3,3) {1.80\%\\($392.38\times$)};
\node[font=\small, align=center] at (axis cs:4,3) {1.53\%\\($387.38\times$)};
\node[font=\small, align=center] at (axis cs:5,3) {1.54\%\\($366.83\times$)};
\node[font=\small, align=center] at (axis cs:1,4) {9.40\%\\($305.18\times$)};
\node[font=\small, align=center] at (axis cs:2,4) {4.14\%\\($275.86\times$)};
\node[font=\small, align=center] at (axis cs:3,4) {1.70\%\\($263.85\times$)};
\node[font=\small, align=center] at (axis cs:4,4) {1.26\%\\($259.92\times$)};
\node[font=\small, align=center] at (axis cs:5,4) {1.27\%\\($241.38\times$)};
\node[font=\small, align=center] at (axis cs:1,5) {9.42\%\\($212.95\times$)};
\node[font=\small, align=center] at (axis cs:2,5) {4.16\%\\($188.20\times$)};
\node[font=\small, align=center] at (axis cs:3,5) {1.69\%\\($185.59\times$)};
\node[font=\small, align=center] at (axis cs:4,5) {1.20\%\\($178.98\times$)};
\node[font=\small, align=center] at (axis cs:5,5) {1.16\%\\($163.88\times$)};

\end{axis}
\end{tikzpicture}
    }
    \caption{Validation error for the remaining primal and tangential dimensions.}
    \label{fig:contact_rom_heatmap}
  \end{subfigure}

  \caption{Selection of the reduced dimensions for the three-dimensional frictional-contact problem. Left: mean relative displacement energy error as a function of the normal-dual dimension $m_n$, for $n=m_\tau\in\{100,150,200\}$; the dotted line indicates the selected value $m_n=50$. Right: mean relative displacement energy error of the Galerkin ROM on the validation set for the considered primal and tangential dimensions, with $m_n=50$. The values in parentheses indicate the online speed-up with respect to the HDM.}
  \label{fig:contact_dimension_selection}
\end{figure}

We evaluate the Galerkin ROM equipped with the POD--RBF elastic-predictor surrogate on the validation set for several choices of the primal, normal-dual, and tangential-dual reduced dimensions. The reported error is the mean relative energy error on the reconstructed displacement,
\[e_u(\mu) = \frac{ \left\| \widehat Y_\mu+Vq_\mu-U_\mu \right\|_K}{\left\| U_\mu \right\|_K}.\]
Based on preliminary validation tests, we fix the normal-dual dimension to $m_n=50$. As shown in Figure~\ref{fig:normal_dimension_selection}, increasing $m_n$ beyond $50$ does not consistently improve the displacement error and can even lead to a slight degradation as for $n=m_\tau=100$.
Figure~\ref{fig:contact_rom_heatmap} then reports the mean relative energy error obtained by varying the number of primal correction modes and tangential-dual modes.

\subsubsection{NN-ROM training and model selection}
\label{sec:contact_nnrom_training}

Following the NN-ROM training procedure described above, the low-fidelity training samples are generated at independent Sobol parameter values by solving the Galerkin ROM with the POD--RBF elastic-predictor surrogate with dimensions
\[
n_{\mathrm{tot}}=150, \qquad m_n=50, \qquad m_{\tau,\mathrm{tot}}=150.
\]

In the contact problem, this nonlinear correction is applied to the primal correction and tangential-dual coordinates, while the normal multiplier remains represented by the linear reduced coordinates $\xi_{\mu,n}\in\mathbb R^{50}$. We retain $25$ online coordinates for both the primal correction and the tangential multiplier, and learn the remaining $125$ coordinates through
\[
q_\mu^c \approx \mathcal N^u_{\theta_u}(q_\mu^r,\eta_\mu),
\qquad
\xi_{\mu,\tau}^c \approx \mathcal N^\tau_{\theta_\tau}(\xi_{\mu,\tau}^r,\eta_\mu),
\]
where
\[
\eta_\mu=(Y_{CL},Z_{CL},Y_{CR},Z_{CR},F_\mu).
\]
The online NN-ROM unknowns therefore are $
(q_\mu^r,\xi_{\mu,n},\xi_{\mu,\tau}^r)
\in
\mathbb R^{25}\times\mathbb R^{50}\times\mathbb R^{25}.$
The two networks are trained from $7000$ low-fidelity Galerkin-ROM samples.

For the three-dimensional problem, the selected primal network uses SiLU activations and two hidden layers of width $512$, while the tangential-dual network uses SiLU activations and two hidden layers of width $64$. Both networks are trained according to the procedure described in Appendix~\ref{app:nn_architecture_selection}.

\subsubsection{Hyper-reduction configuration}
\label{sec:contact_hyper_reduction}

In this problem, the main remaining online cost comes from the tangential Coulomb-projection term. We therefore apply hyper-reduction only to this contribution, while keeping the primal equilibrium and normal-contact residuals in their non-hyper-reduced form for both the HR-ROM and the NN-HR-ROM.

In both cases, the cubature rule is fitted over contact faces from $250$ converged reduced solutions. The greedy NNLS procedure uses a relative tolerance of $\varepsilon_{\mathrm{cub}}^{\mathrm{contact}}=2\times10^{-2}$ and a maximum budget of $750$ sampled contact faces.

The resulting cubature rules contain $662$ contact faces for the HR-ROM and $275$ contact faces for the NN-HR-ROM.

\subsubsection{Results}
Table~\ref{tab:ROM_results_3D} summarizes the performance of the different models for the three-dimensional frictional contact problem on the independent test set of $N_{\mathrm{test}}=200$ parameter values. The standard ROM achieves the lowest mean relative energy error, $1.17\%$, with a speed-up factor of $40.3$. Introducing the surrogate for the elastic predictor increases the speed-up to $294$, while maintaining a comparable error of $1.28\%$. The HR-ROM further reduces the mean solve time to $0.59\,\mathrm{s}$, corresponding to a speed-up factor of $861$, with an error of $1.37\%$ and $662$ sampled contact faces. The NN-ROM reduces the primal and tangential reduced coordinates from $150$ to $25$, reaching a speed-up factor of $1494$ with a mean error of $1.54\%$. Finally, the proposed NN-HR-ROM achieves the largest speed-up factor, $2209$, with a mean energy error of $1.55\%$ and $275$ sampled contact faces. In this test case, the main acceleration thus comes from the NN augmentation, with hyper-reduction providing an additional factor of about $1.5$.

{ \scriptsize
  \begin{table}[pos=H]
    \centering
    \setlength{\tabcolsep}{3pt}
    \renewcommand{\arraystretch}{1.8}
    \begin{tabular}{|l||c|c|c||c||c||c||c|}
      \hline
      \multirow{2}{*}{\textbf{Method}} &
        \multicolumn{3}{c||}{\textbf{Dimension}} &
        \multicolumn{1}{c||}{\textbf{Cubature}} &
        \multirow{2}{*}{\makecell{\textbf{Mean relative}\\\textbf{energy error}}} &
        \multirow{2}{*}{\makecell{\textbf{Mean}\\\textbf{solve}\\\textbf{time (s)}}} &
        \multirow{2}{*}{\makecell{\textbf{Speed--up}\\\textbf{factor}\\\textbf{vs HDM}}} \\
      \cline{2-5}
        & \textbf{Primal}
        & \makecell{\textbf{Normal}\\\textbf{dual}}
        & \makecell{\textbf{Tangential}\\\textbf{dual}}
        & \textbf{Faces}
        & & & \\
      \hline

      \textbf{HDM}
        & $54{,}210$
        & $1{,}702$
        & $3{,}404$
        & --
        & --
        & $508$
        & -- \\ \hline

      \textbf{ROM}
        & $150$
        & $50$
        & $150$
        & --
        & {$1.17 \%$}        
        & {$12.6$}
        & {\textbf{40.32}} \\ \hline

      \textbf{ROM + surr}
        & $150$
        & $50$
        & $150$
        & --
        & {$1.28 \%$}
        & {$1.73$}
        & {\textbf{293.6}} \\ \hline

      \textbf{HR-ROM}
        & $150$
        & $50$
        & $150$
        & $662$
        & {$1.37 \%$}
        & {$0.59$}
        & {\textbf{861.0}} \\ \hline

      \textbf{NN-ROM}
        & $25$
        & $50$
        & $25$
        & --
        & {$1.54 \%$}
        & {$0.34$}
        & {\textbf{1494}} \\ \hline

      \textbf{NN-HR-ROM}
        & $25$
        & $50$
        & $25$
        & $275$
        & {$1.55 \%$}
        & {$0.23$}
        & {\textbf{2209}} \\ \hline

    \end{tabular}
    \caption{Performance of the HDM and reduced models on the independent test set. The reported error is the mean relative primal energy error and the computation time is the mean online solution time. Because of its runtime, the HDM time is averaged over five test points, solved sequentially on a single thread. For the NN-augmented models, the primal and tangential-dual entries denote the retained online coordinates; the normal-dual dimension remains $50$. The ``Faces'' column gives the number of contact faces selected by the face-wise cubature rule; each selected face carries two tangential components. Offline costs are not included in the reported online times.}
    \label{tab:ROM_results_3D}
  \end{table}
}

\paragraph{Constraint feasibility and mechanical response.}
We also check how well each model satisfies the contact conditions over the entire interface; the results are reported in Appendix~\ref{app:constraint_diagnostics} and Table~\ref{tab:contact_constraint_analysis}. The mean residual penetration remains between $0.288\%$ and $0.321\%$ for all reduced variants. The exact-predictor ROM and the ROM using the POD--RBF predictor surrogate give essentially identical feasibility indicators, showing that the predictor approximation has negligible influence on the reconstructed contact constraints. The NN-augmented models show larger violations of the contact-force conditions, with mean negative-normal-traction indicators of approximately $1.05\%$ and friction-bound excesses of approximately $4.1$--$4.3\%$. Their similar values with and without hyper-reduction indicate that these deviations are mainly associated with the NN augmentation. Nevertheless, the mean mechanical-response error remains below $1.82\%$ for all five reduced variants. Thus, the combined normal and tangential contact-force error has only a limited effect on the displacement represented by each model.

\paragraph{Offline data-generation and learning costs.} The offline cost of the three-dimensional problem is dominated by HDM snapshot generation, which required $21.1\,\mathrm{h}$ for the $1500$ retained snapshots. In comparison, low-fidelity data generation, POD--RBF surrogate construction, and neural-network training required $14.5$, $9.53$, and $10.2\,\mathrm{min}$ respectively. These additional stages therefore amount to about $34\,\mathrm{min}$ in total, i.e. only $2.7\%$ of the HDM snapshot-generation wall time on the stated hardware.

\section{Conclusion}
\label{sec:conclusion}

We have introduced a model-order-reduction framework for nonlinear parametrized variational inequalities written as projected primal--dual systems. The high-dimensional problems are solved by a semi-smooth Newton method in primal--dual form. Separate POD spaces are constructed for the primal and dual variables, and a neural-network augmentation is used to express the complementary POD coordinates as nonlinear functions of a small number of retained coordinates and parameter features. The learned reconstruction is embedded directly in the reduced nonlinear solve rather than being applied as a post-processing correction.

To reduce the remaining dependence of the online residual and generalized Jacobian evaluations on the high-dimensional discretization, we combined the reduced formulations with a sparse nonnegative cubature rule. The cubature construction acts directly on the projected nonlinear and complementarity contributions. This feature is particularly relevant for the NN-ROM because the cubature fitting problem is then posed in the smaller tangent-projected equation space associated with the nonlinear manifold.

The numerical study on the two-dimensional nonlinear obstacle problem shows that the linear ROM provides accurate approximations but retains a significant full-order assembly cost. The NN-ROM reduces the number of online primal and dual coordinates from $150$ to $25$ while preserving a comparable primal energy error. Hyper-reduction produces the largest computational gains: the NN-HR-ROM reaches a speed-up factor of approximately $619$ with respect to the HDM in the reported test configuration.

The three-dimensional frictional contact problem further demonstrates the applicability of the framework to Coulomb contact. The correction-based primal representation and the surrogate approximation of the elastic predictor substantially reduce the online cost before hyper-reduction is introduced. The NN-ROM then reduces the retained primal and tangential coordinates from $150$ to $25$, while the normal multiplier remains represented in a linear reduced space. Unlike in the obstacle problem, this reduction provides the main computational gain, with a speed-up factor of $1494$ without hyper-reduction. In the reported configuration, the NN-HR-ROM reaches a speed-up factor of approximately $2209$ with respect to the HDM, compared with $861$ for the HR-ROM, while the mean relative energy error increases moderately from $1.17\%$ for the Galerkin ROM to $1.55\%$.

Residual penetration and mechanical-response errors remain small on average, but the reconstructed dual fields do not preserve pointwise feasibility. Future work will therefore focus on feasibility-preserving dual approximations, particularly when local multipliers or contact forces are quantities of interest.


\subsection*{Acknowledgements}
Virginie Ehrlacher and Sofiane Ezzehi acknowledge financial support from the European Research Council (ERC) under the European Union's Horizon 2020 Research and Innovation Programme, Grant Agreement No.~101077204 (HighLEAP).

\subsection*{Declaration of generative AI and AI-assisted technologies in the manuscript preparation process}
During the preparation of this work, the authors used ChatGPT to assist with improving the language and clarity of primarily author-written manuscript drafts, as well as with code debugging, restructuring, and refinement. All AI-assisted content and code changes were reviewed and tested by the authors, who take full responsibility for the manuscript and reported results.

\subsection*{Data availability}
    Data will be made available on request.

\FloatBarrier
\appendix
  \counterwithin{equation}{section}
  \counterwithin{table}{section}
  \counterwithin{figure}{section}
\section{Active-set SSN formulation of the two test cases}\label{sec:appendix_ssn_derivation}
    \subsection{Active-set form for the $2D$ obstacle problem with cubic nonlinearity}
    In the case of the $2D$ obstacle problem, the residual function $\mathcal{F}^{(1)}$ is given by:
    $$\mathcal{F}^{(1)}(U, \Lambda) = \begin{pmatrix} K U + \gamma M U^{\circ 3} - M \Lambda \\ \Lambda - \max(0, \Lambda - \rho (U - G)) \end{pmatrix}.$$
    We see that the state equation residual is smooth, while the constraint residual is nonsmooth due to the presence of the pointwise $\max$ operator. Using the slanting function $Q'(X)=D_{>0}(X)$ defined in \eqref{eq:slanting_function_max}, a slanting function of $\mathcal F^{(1)}$ is:
    \begin{equation*}
    \mathcal{F}^{(1)'}(U, \Lambda) = \begin{pmatrix} K + 3\gamma M \text{diag}(U^{\circ 2}) & -M \\ \rho D_{>0}(\Lambda - \rho (U - G)) & I - D_{>0}(\Lambda - \rho (U - G)) \end{pmatrix},
    \end{equation*}
    where $\text{diag}(U^{\circ 2})$ is the diagonal matrix with entries $[\text{diag}(U^{\circ 2})]_{ii} = (U_i)^2$.
    
    Therefore, the primal--dual active-set formulation is obtained directly from the SSN system. At iteration $k$, the active and inactive sets $\mathcal A^k = \{i:\Lambda_i^k-\rho(U_i^k-G_i)>0\}$ and $\mathcal I^k = \{i:\Lambda_i^k-\rho(U_i^k-G_i)\leq0\}$ are first updated. The Newton direction then satisfies
    \[ (d_U^k)_{\mathcal A^k} = (G)_{\mathcal A^k}-(U^k)_{\mathcal A^k}, \qquad (d_\Lambda^k)_{\mathcal I^k} = -(\Lambda^k)_{\mathcal I^k}. \]
    The accepted iterate is then set to
    \[ X^{k+1}=X^k+\alpha_k d^k. \]
    where $\alpha_k$ is selected by the line-search procedure described in Appendix~\ref{app:ssn_numerical_settings}.
        
        \subsection{Active-set form for the three-dimensional contact problem}

    During an SSN iteration, the parameter $\mu$ is fixed, and we omit it from the notation. In the following, the residual of the three-dimensional contact problem is
    \[
    \mathcal F^{(2)}(U,\Lambda_n,\Lambda_\tau)\coloneqq\begin{pmatrix}R_U(U,\Lambda_n,\Lambda_\tau)\\R_n(U,\Lambda_n)\\R_\tau(U,\Lambda_n,\Lambda_\tau)\end{pmatrix}=
    \begin{pmatrix}
    KU+B_n^T\Lambda_n+B_\tau^T\Lambda_\tau-L\\
    \Lambda_n-\max\bigl(0,\Lambda_n+\rho B_nU\bigr)\\
    \Lambda_\tau-\Pi_{c(\Lambda_n)}\bigl(\Lambda_\tau+\rho B_\tau U\bigr)
    \end{pmatrix}.
    \]
    The state-equation residual is smooth, whereas the normal and tangential constraint residuals are nonsmooth because of the componentwise $\max$ operator and the facewise projection onto the friction disks, respectively.

    The SSN direction $d^k=(d_U^k,d_{\Lambda_n}^k,d_{\Lambda_\tau}^k)$ satisfies
    \[
    \mathcal F^{(2)\prime}(U^k,\Lambda_n^k,\Lambda_\tau^k)d^k=-\mathcal F^{(2)}(U^k,\Lambda_n^k,\Lambda_\tau^k),
    \]
    or, equivalently,
    \[
    \begin{aligned}
    Kd_U^k+B_n^Td_{\Lambda_n}^k+B_\tau^Td_{\Lambda_\tau}^k&=-R_U(U^k,\Lambda_n^k,\Lambda_\tau^k),\\
    R_n'(U^k,\Lambda_n^k)[d_U^k,d_{\Lambda_n}^k]&=-R_n(U^k,\Lambda_n^k),\\
    R_\tau'(U^k,\Lambda_n^k,\Lambda_\tau^k)[d_U^k,d_{\Lambda_n}^k,d_{\Lambda_\tau}^k]&=-R_\tau(U^k,\Lambda_n^k,\Lambda_\tau^k).
    \end{aligned}
    \]
    Here, $R_n'$ and $R_\tau'$ denote the selected slanting functions of the normal and tangential residuals, respectively.
    
    \subsubsection*{Normal-contact equations}
    At iteration $k$, let
    \[
    a_n^k=\Lambda_n^k+\rho B_nU^k,\qquad \mathcal A_n^k=\{j:a_{n,j}^k>0\},\qquad \mathcal I_n^k=\{j:a_{n,j}^k\leq0\}.
    \]
    Equality is thus assigned to the inactive set, consistently with the choice $D_{>0}$ in \eqref{eq:slanting_function_max}. Define the diagonal matrix
    \[
    D_n^k\coloneqq D_{>0}(a_n^k),\qquad (D_n^k)_{jj}=\begin{cases}1,&j\in\mathcal A_n^k,\\0,&j\in\mathcal I_n^k.\end{cases}
    \]
    Since $\max(0,a_n^k)=D_n^ka_n^k$, the normal residual and its selected slanting function are
    \[
    R_n(U^k,\Lambda_n^k)=(I-D_n^k)\Lambda_n^k-\rho D_n^kB_nU^k
    \]
    and
    \[
    R_n'(U^k,\Lambda_n^k)[d_U^k,d_{\Lambda_n}^k]=(I-D_n^k)d_{\Lambda_n}^k-\rho D_n^kB_nd_U^k.
    \]
    Consequently, the normal-contact row of the SSN system can be written as
    \[
    (I-D_n^k)(\Lambda_n^k+d_{\Lambda_n}^k)-\rho D_n^kB_n(U^k+d_U^k)=0.
    \]
    Restricting this equation to the active and inactive sets gives
    \[
    \boxed{(B_nd_U^k)_{\mathcal A_n^k}=-(B_nU^k)_{\mathcal A_n^k}}\quad \text{and}\quad \boxed{(d_{\Lambda_n}^k)_{\mathcal I_n^k}=-(\Lambda_n^k)_{\mathcal I_n^k}.}
    \]

    \subsubsection*{Tangential-contact equations}
    Denote the projected term in the tangential residual by
    \[
    \phi_\tau(U,\Lambda_n,\Lambda_\tau)\coloneqq\Pi_{c(\Lambda_n)}\bigl(\Lambda_\tau+\rho B_\tau U\bigr),
    \]
    so that $R_\tau(U,\Lambda_n,\Lambda_\tau)=\Lambda_\tau-\phi_\tau(U,\Lambda_n,\Lambda_\tau)$. For each interface face $j$, let
    \[
    s_j^k=\Lambda_{\tau,j}^k+\rho(B_\tau U^k)_j,\qquad r_j^k=c_j(\Lambda_n^k)=F[\Lambda_{n,j}^k]_+,
    \]
    where $s_j^k\in\mathbb R^2$ and $r_j^k\geq0$. The sticking and slipping sets are defined on all interface faces by
    \[
    \mathcal S_\tau^k=\{j:\|s_j^k\|_2\leq r_j^k\},\qquad \mathcal T_\tau^k=\{j:\|s_j^k\|_2>r_j^k\}.
    \]
    Thus, faces satisfying $\|s_j^k\|_2=r_j^k$ are assigned to the sticking set. We also introduce the selected slanting derivative of the positive-part mapping,
    \[
    \chi_j^k=\begin{cases}1,&\Lambda_{n,j}^k>0,\\0,&\Lambda_{n,j}^k\leq0.\end{cases}
    \]
    The facewise projection and its selected slanting function are described as follows:
    \begin{itemize}
    \item \textbf{Sticking faces.} For $j\in\mathcal S_\tau^k$, the projection is the identity $\phi_{\tau,j}(U^k,\Lambda_n^k,\Lambda_\tau^k)=s_j^k$.
    Consequently,
    \[
    \phi_{\tau,j}'(U^k,\Lambda_n^k,\Lambda_\tau^k)[d_U^k,d_{\Lambda_n}^k,d_{\Lambda_\tau}^k]=\rho(B_\tau d_U^k)_j+d_{\Lambda_\tau,j}^k.
    \]
    It follows that
    \[
    R_{\tau,j}'(U^k,\Lambda_n^k,\Lambda_\tau^k)[d_U^k,d_{\Lambda_n}^k,d_{\Lambda_\tau}^k]=-\rho(B_\tau d_U^k)_j.
    \]
    And since $R_{\tau,j}(U^k,\Lambda_n^k,\Lambda_\tau^k)=-\rho(B_\tau U^k)_j$, the tangential-contact row of the SSN system gives
    \[
    \boxed{(B_\tau d_U^k)_{\mathcal S_\tau^k}=-(B_\tau U^k)_{\mathcal S_\tau^k}.}
    \]

    \item \textbf{Slipping faces.} For $j\in\mathcal T_\tau^k$, we define
    \[
    \widehat s_j^k=\frac{s_j^k}{\|s_j^k\|_2} \in \mathbb R^2,
    \qquad A_j^k=\frac{r_j^k}{\|s_j^k\|_2}\left(I_2-\widehat s_j^k(\widehat s_j^k)^T\right) \in\mathbb R^{2\times2},
    \qquad v_j^k=F\chi_j^k\widehat s_j^k \in\mathbb R^2.
    \]
    The projection is
    \[
    \phi_{\tau,j}(U^k,\Lambda_n^k,\Lambda_\tau^k)=r_j^k\widehat s_j^k,
    \]
    and the action of its selected slanting function is
    \[
    \phi_{\tau,j}'(U^k,\Lambda_n^k,\Lambda_\tau^k)[d_U^k,d_{\Lambda_n}^k,d_{\Lambda_\tau}^k]=A_j^k\bigl(\rho(B_\tau d_U^k)_j+d_{\Lambda_\tau,j}^k\bigr)+v_j^kd_{\Lambda_n,j}^k.
    \]
    \end{itemize}

    To assemble the slipping equations, we need to collect the slipping faces. We write $\mathcal T_\tau^k=\{j_1,\ldots,j_{N_T}\}$ and define
    \[
    J_1^k=\begin{pmatrix}
    A_{j_1}^k&\mathbf0_{2\times2}&\cdots&\mathbf0_{2\times2}\\
    \mathbf0_{2\times2}&A_{j_2}^k&\ddots&\vdots\\
    \vdots&\ddots&\ddots&\mathbf0_{2\times2}\\
    \mathbf0_{2\times2}&\cdots&\mathbf0_{2\times2}&A_{j_{N_T}}^k
    \end{pmatrix}\in\mathbb R^{2N_T\times2N_T},
    \]
    and
    \[
    J_2^k=\begin{pmatrix}
    v_{j_1}^k&\mathbf0_{2\times1}&\cdots&\mathbf0_{2\times1}\\
    \mathbf0_{2\times1}&v_{j_2}^k&\ddots&\vdots\\
    \vdots&\ddots&\ddots&\mathbf0_{2\times1}\\
    \mathbf0_{2\times1}&\cdots&\mathbf0_{2\times1}&v_{j_{N_T}}^k
    \end{pmatrix}\in\mathbb R^{2N_T\times N_T}.
    \]
    The slipping part of the selected slanting function of the tangential residual is therefore
    \[
    \begin{aligned}
    R_{\tau,\mathcal T_\tau^k}'(U^k,\Lambda_n^k,\Lambda_\tau^k)[d_U^k,d_{\Lambda_n}^k,d_{\Lambda_\tau}^k]
    ={}&(I_{2N_T}-J_1^k)(d_{\Lambda_\tau}^k)_{\mathcal T_\tau^k}\\
    &-\rho J_1^kB_{\tau,\mathcal T_\tau^k}d_U^k-J_2^k(d_{\Lambda_n}^k)_{\mathcal T_\tau^k}.
    \end{aligned}
    \]
    Accordingly, the slipping rows of the SSN system are
    \[
    \boxed{(I_{2N_T}-J_1^k)(d_{\Lambda_\tau}^k)_{\mathcal T_\tau^k}-\rho J_1^kB_{\tau,\mathcal T_\tau^k}d_U^k-J_2^k(d_{\Lambda_n}^k)_{\mathcal T_\tau^k}=-R_{\tau,\mathcal T_\tau^k}(U^k,\Lambda_n^k,\Lambda_\tau^k).}
    \]

    \begin{rmk}
    In practice, we make the additional choice to apply the sticking and slipping equations only on faces currently classified as normal active. We therefore define
    \[
    \mathcal S_{\tau,A}^k=\mathcal S_\tau^k\cap\mathcal A_n^k,\qquad \mathcal T_{\tau,A}^k=\mathcal T_\tau^k\cap\mathcal A_n^k,
    \]
    and assemble $J_1^k$ and $J_2^k$ over $\mathcal T_{\tau,A}^k$. On the normal inactive set, we directly set the tangential multiplier of the full SSN candidate to zero:
    \[
    \boxed{(d_{\Lambda_\tau}^k)_{\mathcal I_n^k}=-(\Lambda_\tau^k)_{\mathcal I_n^k}.}
    \]
    This choice enforces the physical condition that an open face cannot transmit tangential traction and constitutes an additional active-set elimination rule.
    \end{rmk}

    Together with the state-equation row, the boxed contact relations form the coupled linear system used to compute the SSN direction. In the implementation, this system is reduced by eliminating the displacement variables and solving the resulting Schur-complement system for the active multipliers.
       
  \section{Numerical settings for the SSN solver}
    \label{app:ssn_numerical_settings}

    All HDM and reduced solvers use a damped SSN strategy with backtracking, with model-dependent merit functions and stopping tolerances. Let $X^k$ denote the current iterate, $d^k$ the SSN direction, and $\mathcal M$ the model-dependent merit function. Starting from $\alpha=1$, the step lengths
    \[ \alpha\in\{1,\beta_{\mathrm{ls}},\ldots, \beta_{\mathrm{ls}}^{J_{\max}}\} \]
    are tested until
    \[ \mathcal M(X^k+\alpha d^k) \leq (1-c_\mathrm{ls}\alpha)\mathcal M(X^k), \]
    or until the trial merit function is below a prescribed tolerance. If no trial satisfies either condition, the tested step with the lowest merit score is retained:
    \[ \alpha_k\in \operatorname*{arg\,min}_{\alpha\in \{1,\beta_{\mathrm{ls}},\ldots,\beta_{\mathrm{ls}}^{J_{\max}}\}} \mathcal M(X^k+\alpha d^k). \]
    If several trials attain the same minimum, the largest step length is retained. Furthermore, this fallback is allowed to be non-monotonic. For each solver, the iteration is stopped when the corresponding model-dependent residual $R(X^k)$ is below the prescribed tolerance.
    
    \subsection{Two-dimensional obstacle}
    \subsubsection{High-dimensional model}

    For the HDM obstacle solver, the line-search merit function and stopping residual coincide:
    \[ \mathcal M_{\mathrm{obs}}^{h}(U,\Lambda) = R_{\mathrm{obs}}^{h}(U,\Lambda) = \max\left\{ \left\|KU+\gamma M U^{\circ 3}-M\Lambda\right\|_2,\, \left\| \Lambda-\max\left(0,\Lambda-\rho(U-G)\right) \right\|_\infty \right\}. \]

    \subsubsection{Reduced-order models}

    For the Galerkin ROM, the NN-ROM, and their hyper-reduced counterparts, the line-search merit function and stopping residual are
    \[ \mathcal M_{\mathrm{obs}}^{r}(X_r) = R_{\mathrm{obs}}^{r}(X_r) = \max\left\{ \|r_u^r(X_r)\|_2,\, \|r_c^r(X_r)\|_\infty \right\}, \]
    where $r_u^r$ and $r_c^r$ are the equilibrium and projected complementarity residuals of the considered reduced model. For the NN-ROM, these are the corresponding tangent-projected residuals, while the hyper-reduced solvers use their hyper-reduced counterparts.

      \paragraph{Numerical values.} In the two-dimensional test case, all line searches start from $\alpha_0=1$ and use $\beta_{\mathrm{ls}}=0.75$, $c_{\mathrm{ls}}=10^{-2}$, and $J_{\max}=10$. For the HDM, the projection parameter is $\rho=1$, the maximum number of SSN iterations is $100$, and convergence is declared when $R_{\mathrm{obs}}^{h}(U,\Lambda)\leq 10^{-10}$. For all reduced models, we use $\rho=500$, selected empirically in preliminary reduced-model experiments based on the primal approximation accuracy on the validation set. The maximum number of SSN iterations is $100$, and convergence is declared when $R_{\mathrm{obs}}^{r}(X_r)\leq 10^{-8}$.

    \subsection{Three-dimensional contact}
    \subsubsection{High-dimensional model}

    For the HDM contact solver, the line-search merit function is the tangential projection residual
    \[ \mathcal M_{\tau}(U,\Lambda_n,\Lambda_\tau) = \left\| \Lambda_\tau- \Pi_{c_\mu(\Lambda_n)} \left(\Lambda_\tau+\rho B_\tau U\right) \right\|_{\infty,2}, \]
    where $\|\cdot\|_{\infty,2} \coloneqq \max_{j=1,\ldots,R} \|Z_j\|_2$ denotes the maximum Euclidean norm over contact faces of a vector $Z\in\mathbb{R}^{2R}$ with facewise components $Z_j\in\mathbb{R}^2$.
    
    This quantity was chosen as the line-search merit function because the tangential projection residual was empirically the most difficult residual component to reduce. Using it substantially improved the robustness of convergence in the numerical experiments.

    The stopping residual is
    \[ R_{\mathrm{ct}}^{h}(U,\Lambda_n,\Lambda_\tau) = \max\left\{ r_{\mathrm{eq}},\,r_n,\,r_\tau \right\}, \]
    with
    \[ r_{\mathrm{eq}} = \frac{ \|KU-L_\mu+B_n^T\Lambda_n+B_\tau^T\Lambda_\tau\|_2 }{ \max\{1,\|L_\mu\|_2\} }, \] \[ r_n = \left\| \Lambda_n-\max\left(0,\Lambda_n+\rho B_nU\right) \right\|_\infty, \]
    and
    \[ r_\tau = \left\| \Lambda_\tau- \Pi_{c_\mu(\Lambda_n)} \left(\Lambda_\tau+\rho B_\tau U\right) \right\|_{\infty,2}. \]

    \subsubsection{Reduced-order models}

    For the Galerkin ROM, the NN-ROM, and their hyper-reduced counterparts, the line-search merit function and stopping residual coincide:
    \[ \mathcal M_{\mathrm{ct}}^{r}(X_r) = R_{\mathrm{ct}}^{r}(X_r) = \max\left\{ \|r_u^r(X_r)\|_\infty,\, \|r_n^r(X_r)\|_\infty,\, \|r_\tau^r(X_r)\|_\infty \right\}. \]
    Here $r_u^r$, $r_n^r$, and $r_\tau^r$ are the reduced equilibrium, normal-contact, and tangential-friction residuals of the considered model. For the hyper-reduced variants, they are replaced by their hyper-reduced counterparts. Since the HDM and reduced residuals use different scalings and norms, their stopping tolerances are not directly comparable and were calibrated separately; the numerical values used in each formulation are reported below.

    \paragraph{Numerical values.} In the three-dimensional contact test case, all line searches start from $\alpha_0=1$ and use $\beta_{\mathrm{ls}}=0.75$, $c_{\mathrm{ls}}=10^{-2}$, and $J_{\max}=10$. For the HDM, the projection parameter is $\rho=10^{-2}$, the maximum number of SSN iterations is $75$, and convergence requires $r_{\mathrm{eq}}\leq 10^{-6}$, $r_n\leq 10^{-6}$, and $r_\tau\leq 10^{-4}$. The tangential projection residual is also used to control the HDM line search, with tolerance $10^{-4}$. For all reduced models, the projection parameter is $\rho=10^{2}$, the maximum number of SSN iterations is $100$, and convergence is declared when $R_{\mathrm{ct}}^{r}(X_r)\leq 10^{-4}$.

\section{Constraint feasibility and mechanical response}
\label{app:constraint_diagnostics}
The relative primal energy error does not by itself quantify whether the reconstructed reduced fields satisfy the inequality and admissibility conditions of the underlying variational inequality. We therefore report two complementary types of diagnostics. The problem-specific feasibility indicators measure violations of the relevant constraint conditions. The relative mechanical-response error measures instead the displacement-energy response caused by the error in the reconstructed obstacle or contact force, relative to the response caused by the corresponding HDM force. It therefore quantifies the influence of the dual reconstruction error on the primal equilibrium.

For every reduced model, the primal and dual fields are first reconstructed on the high-dimensional discretization. The diagnostics are then evaluated using the complete set on which the constraints are imposed: all constrained degrees of freedom for the two-dimensional obstacle problem and all contact faces for the three-dimensional frictional contact problem. Consequently, the analyses of the HR-ROM and NN-HR-ROM are not restricted to the entries or faces selected by their cubature rules.

These quantities are a posteriori offline evaluation diagnostics, since their computation requires the corresponding HDM reference fields. Furthermore, for every parameter sample and reduced model considered in this evaluation, all normalization and reference-response denominators appearing below are strictly positive. Hence, no denominator regularization is required and all reported relative indicators are well defined.

\subsection{Two-dimensional obstacle problem}
\label{app:obstacle_constraint_analysis}

For a nodal scalar field $Z\in\mathbb{R}^{N}$, we use the mass-weighted norm
\[ \lVert Z\rVert_M=\left(Z^TMZ\right)^{1/2}, \]
where $M$ is the finite element mass matrix introduced in Section~\ref{sec:2D_obstacle_problem}. We denote by $U_\mu^{\mathrm{red}}$ and $\Lambda_\mu^{\mathrm{red}}$ the reconstructed reduced primal and multiplier vectors, respectively, and by $U_\mu$ and $\Lambda_\mu$ the corresponding HDM solutions.

\paragraph{\textbf{Relative residual penetration.}} The first metric measures the remaining violation of the discrete obstacle inequality $U_\mu^{\mathrm{red}}-G_\mu\geq0$. For each parameter $\mu$, it is defined as
\[ 100\,\frac{\left\lVert\left[G_\mu-U_\mu^{\mathrm{red}}\right]_+\right\rVert_M}{\left\lVert G_\mu\right\rVert_M}. \]
The numerator retains only the degrees of freedom at which the reconstructed solution penetrates the obstacle, while the denominator provides a reference scale based on the prescribed obstacle.

\paragraph{\textbf{Relative negative multiplier.}} The second metric checks whether the reconstructed multiplier remains nonnegative. For each parameter $\mu$, it is defined as
\[ 100\,\frac{\left\lVert\left[-\Lambda_\mu^{\mathrm{red}}\right]_+\right\rVert_M}{\left\lVert\Lambda_\mu\right\rVert_M}. \]
The numerator measures the magnitude of the inadmissible negative part of the reconstructed multiplier, while the HDM multiplier norm provides the corresponding reference scale.

\paragraph{\textbf{Relative obstacle-constraint residual.}} The third metric measures the overall obstacle-constraint violation using the projection residual in \eqref{eq:discrete_obstacle_nonsmooth_form}. For each parameter $\mu$, it is defined as
\[ 100\,\frac{\left\lVert\Lambda_\mu^{\mathrm{red}}-\left[\Lambda_\mu^{\mathrm{red}}-\rho\left(U_\mu^{\mathrm{red}}-G_\mu\right)\right]_+\right\rVert_M}{\left\lVert\Lambda_\mu\right\rVert_M+\rho\left\lVert\left[U_\mu-G_\mu\right]_+\right\rVert_M}. \]
The numerator vanishes if and only if the reconstructed fields satisfy coefficientwise nonpenetration, multiplier nonnegativity, and complementarity. The denominator combines the HDM multiplier scale and the admissible HDM separation-gap scale.

\paragraph{\textbf{Relative mechanical-response error.}}
The primal equation of the obstacle problem is written as
\[ KU_\mu+\gamma MU_\mu^{\circ 3}=M\Lambda_\mu. \]
The term $M\Lambda_\mu$ is the force exerted by the obstacle on the primal solution. Therefore, an error in the reconstructed multiplier produces the force error $M\Delta\Lambda_\mu$, where
$\Delta\Lambda_\mu=\Lambda_\mu^{\mathrm{red}}-\Lambda_\mu$.

Let $V$ be the primal ROM basis. For clarity of presentation, the metric is derived below for the linear ROM. For the NN-ROM, the same construction is obtained by replacing $V$ with the primal reconstruction tangent $T_{\mathcal U_\mu}(q_\mu^r)$ evaluated at the computed reduced solution.

Since the ROM represents displacement only through the directions contained in $V$, it responds to the projected force error\[ V^TM\Delta\Lambda_\mu. \]
Because the primal equation is nonlinear, the displacement produced by this force error is determined by the Jacobian of the primal operator. Evaluated at the HDM solution, this Jacobian is $J_\mu = K+3\gamma M\operatorname{diag}\!\left(U_\mu^{\circ 2}\right)$ and its projection onto the ROM space is
\[ J_{\mu,n}=V^TJ_\mu V. \]

We denote by $d_{\mu,n}^{\mathrm{err}}$ the displacement response produced by the multiplier-force error and by $d_{\mu,n}^{\mathrm{ref}}$ the response produced by the HDM obstacle force. They are given by
\[ d_{\mu,n}^{\mathrm{err}} = J_{\mu,n}^{-1}V^TM\Delta\Lambda_\mu, \qquad d_{\mu,n}^{\mathrm{ref}} = J_{\mu,n}^{-1}V^TM\Lambda_\mu. \]
Thus, the Jacobian converts each force into the corresponding response of the nonlinear primal equation.

Writing $K_n=V^TKV$, the full-space displacement associated with the
reduced response $d_n$ is $Vd_n$. Its elastic energy norm satisfies
\[ \lVert Vd_n\rVert_K = \left[(Vd_n)^TK(Vd_n)\right]^{1/2} = \left(d_n^TK_nd_n\right)^{1/2}. \]
The Jacobian $J_{\mu,n}$ determines the displacement response,
whereas $K_n$ is used only to measure its magnitude.

The relative mechanical-response error is then defined as
\[ 100\, \frac{ \left[ \left(d_{\mu,n}^{\mathrm{err}}\right)^T K_n d_{\mu,n}^{\mathrm{err}} \right]^{1/2} }{ \left[ \left(d_{\mu,n}^{\mathrm{ref}}\right)^T K_n d_{\mu,n}^{\mathrm{ref}} \right]^{1/2} }. \]
Consequently, a small value means that the multiplier reconstruction error produces only a small displacement response relative to that produced by the reference obstacle force.

This diagnostic therefore helps explain why noticeable local multiplier-feasibility errors can coexist with an accurate primal solution. 

\begin{table}[pos=H]
\centering
\begin{tabular}{lcccc}
\toprule
& \multicolumn{3}{c}{\textbf{Constraint feasibility}}
& \multicolumn{1}{c}{\textbf{Mechanical effect}} \\
\cmidrule(lr){2-4}\cmidrule(lr){5-5}
\textbf{Method}
& \shortstack{\textbf{Residual}\\\textbf{penetration}}
& \shortstack{\textbf{Negative}\\\textbf{multiplier}}
& \shortstack{\textbf{Obstacle-constraint}\\\textbf{residual}}
& \shortstack{\textbf{Mechanical-response}\\\textbf{error}} \\
\midrule
ROM       & $0.164$ & $6.544$ & $0.167$ & $2.703$ \\
HR-ROM    & $0.183$ & $6.685$ & $0.177$ & $3.269$ \\
NN-ROM    & $0.210$ & $7.259$ & $0.196$ & $2.757$ \\
NN-HR-ROM & $0.225$ & $7.413$ & $0.203$ & $3.245$ \\
\bottomrule
\end{tabular}
\caption{Mean full-space constraint-feasibility and mechanical-response diagnostics over the test set. All entries are percentages.}
\label{tab:obstacle_constraint_analysis}
\end{table}

Across all four models, the mean residual penetration remains between $0.164\%$ and $0.225\%$, while the mean obstacle-constraint residual remains between $0.167\%$ and $0.203\%$. Both the NN augmentation and hyper-reduction therefore produce only modest changes in primal feasibility and overall obstacle-constraint consistency. The most visible feasibility limitation remains multiplier nonnegativity: the corresponding indicator ranges from $6.544\%$ for the ROM to $7.413\%$ for the NN-HR-ROM. Most of this increase is associated with the NN augmentation, whereas hyper-reduction has a smaller effect within each model family.

The mean mechanical-response error ranges from $2.703\%$ for the ROM to $3.269\%$ for the HR-ROM. Hyper-reduction increases this error from $2.703\%$ to $3.269\%$ for the projection-based models and from $2.757\%$ to $3.245\%$ for the NN-augmented models. In contrast, the NN augmentation has no systematic effect on this quantity, slightly increasing the error without hyper-reduction and slightly decreasing it in the hyper-reduced case. The systematic increase is therefore associated primarily with hyper-reduction. Nevertheless, the mean error remains below $3.3\%$ in every case. Thus, within the displacement directions represented by each model, the multiplier-force error produces only a limited local response relative to the response generated by the HDM obstacle force.

\subsection{Three-dimensional frictional contact problem}
\label{app:contact_constraint_analysis}

Let $\{\Gamma_j\}_{j=1}^{R}$ denote the contact faces. For a face-wise scalar or two-component tangential field $Z$, we use the area-weighted interface norm
\[ \lVert Z\rVert_\Gamma=\left(\sum_{j=1}^{R}|\Gamma_j|\,\lVert Z_j\rVert_2^2\right)^{1/2}. \]
Similarly, we denote by $U_\mu^{\mathrm{red}}$, $\Lambda_{\mu,n}^{\mathrm{red}}$, and $\Lambda_{\mu,\tau}^{\mathrm{red}}$ the reconstructed reduced displacement, normal multiplier, and tangential multiplier, respectively, and by $U_\mu$, $\Lambda_{\mu,n}$, and $\Lambda_{\mu,\tau}$ the corresponding HDM fields.

\paragraph{\textbf{Relative residual penetration.}} The first metric quantifies how effectively the reconstructed reduced solution suppresses the interpenetration that the applied loading would generate in the absence of contact. For each parameter $\mu$, it is defined as
\[ 100\,\frac{\left\lVert\left[B_nU_\mu^{\mathrm{red}}\right]_+\right\rVert_\Gamma}{\left\lVert\left[B_nY_\mu\right]_+\right\rVert_\Gamma}. \]
Since the HDM solution satisfies the nonpenetration condition, $\left[B_nU_\mu\right]_+$ vanishes up to the numerical solver tolerance and therefore cannot provide a meaningful normalization. We instead use the exact unconstrained elastic predictor $Y_\mu$ introduced in Section~\ref{sec:contact_correction_pod}. The denominator measures the penetration that the applied loading would generate in the absence of contact, whereas the numerator measures the penetration remaining in the reconstructed reduced solution. The reported percentage therefore represents the fraction of the unconstrained attempted penetration that remains after the reduced contact correction.

\paragraph{\textbf{Relative negative normal traction.}}
The second metric measures the inadmissible negative part of the reconstructed normal traction relative to the HDM normal-traction scale:
\[ 100\, \frac{\left\lVert\left[-\Lambda_{\mu,n}^{\mathrm{red}}\right]_+\right\rVert_\Gamma} {\left\lVert\Lambda_{\mu,n}\right\rVert_\Gamma}. \]

\paragraph{\textbf{Relative friction-bound excess.}} The third metric checks whether the reconstructed tangential traction exceeds the Coulomb friction bound determined by the reconstructed normal traction. For each parameter $\mu$, it is defined as
\[ 100\,\frac{\left\lVert\left[\left\lVert\Lambda_{\mu,\tau}^{\mathrm{red}}\right\rVert_2-F_\mu\left[\Lambda_{\mu,n}^{\mathrm{red}}\right]_+\right]_+\right\rVert_\Gamma}{\left(\left\lVert\Lambda_{\mu,n}\right\rVert_\Gamma^2+\left\lVert\Lambda_{\mu,\tau}\right\rVert_\Gamma^2\right)^{1/2}}. \]
The Euclidean tangential norm and the positive-part operators in the numerator are applied face by face. On each contact face, the numerator retains only the portion of the reconstructed tangential traction that exceeds the Coulomb bound $F_\mu[\Lambda_{\mu,n}^{\mathrm{red}}]_+$. The interface norm aggregates this excess over the complete contact interface, while the denominator provides a reference scale based on the total HDM contact traction.

\paragraph{\textbf{Relative mechanical-response error.}}
The primal equation of the three-dimensional contact problem is
\[ KU_\mu +B_n^T\Lambda_{\mu,n} +B_\tau^T\Lambda_{\mu,\tau} =L_\mu. \]
Thus, the normal and tangential multipliers enter the equilibrium through the combined contact action
\[ B_n^T\Lambda_{\mu,n} +B_\tau^T\Lambda_{\mu,\tau}. \] 
Defining $\Delta\Lambda_{\mu,n} = \Lambda_{\mu,n}^{\mathrm{red}}-\Lambda_{\mu,n}$ and $\Delta\Lambda_{\mu,\tau} = \Lambda_{\mu,\tau}^{\mathrm{red}}-\Lambda_{\mu,\tau},$ the corresponding contact-force error is
\[ B_n^T\Delta\Lambda_{\mu,n} +B_\tau^T\Delta\Lambda_{\mu,\tau}. \]

Using the primal ROM basis $V$ and the reduced stiffness matrix $K_n=V^TKV$, the relative tangent-resolved contact-action error is therefore defined as follows for the linear ROM:
\[ 100\, \frac{ \left[ \left( V^T \left( B_n^T\Delta\Lambda_{\mu,n} +B_\tau^T\Delta\Lambda_{\mu,\tau} \right) \right)^T K_n^{-1} \left( V^T \left( B_n^T\Delta\Lambda_{\mu,n} +B_\tau^T\Delta\Lambda_{\mu,\tau} \right) \right) \right]^{1/2} }{ \left[ \left( V^T \left( B_n^T\Lambda_{\mu,n} +B_\tau^T\Lambda_{\mu,\tau} \right) \right)^T K_n^{-1} \left( V^T \left( B_n^T\Lambda_{\mu,n} +B_\tau^T\Lambda_{\mu,\tau} \right) \right) \right]^{1/2} }. \]
This metric includes errors in both the normal and tangential contact forces. A small value means that their combined effect on the displacement directions represented by the ROM is small relative to the response generated by the HDM contact force.
As in the $2D$ test case, for the NN-augmented models the same expression is evaluated by replacing $V$ with the primal reconstruction tangent $T_{\mathcal U_\mu}(q_\mu^r)$ at the computed reduced solution, so that $K_n=T_{\mathcal U_\mu}(q_\mu^r)^TKT_{\mathcal U_\mu}(q_\mu^r)$.

By construction, this indicator isolates the model-resolved mechanical response induced by the contact-multiplier reconstruction error; it does not include the displacement-reconstruction error $U_\mu^{\mathrm{red}}-U_\mu$ or the error of the surrogate predictor $\widehat{Y}_\mu$.

\begin{table}[pos=H]
\centering
\begin{tabular}{lcccc}
\toprule
& \multicolumn{3}{c}{\textbf{Constraint feasibility}}
& \multicolumn{1}{c}{\textbf{Mechanical effect}} \\
\cmidrule(lr){2-4}\cmidrule(lr){5-5}
\textbf{Method}
& \shortstack{\textbf{Residual}\\\textbf{penetration}}
& \shortstack{\textbf{Negative normal}\\\textbf{traction}}
& \shortstack{\textbf{Friction-bound}\\\textbf{excess}}
& \shortstack{\textbf{Mechanical-response}\\\textbf{error}} \\
\midrule
ROM & $0.300$ & $0.277$ & $1.174$ & $1.057$ \\
ROM + surr           & $0.299$ & $0.277$ & $1.174$ & $1.062$ \\
HR-ROM               & $0.321$ & $0.277$ & $2.734$ & $1.819$ \\
NN-ROM               & $0.288$ & $1.049$ & $4.109$ & $0.885$ \\
NN-HR-ROM            & $0.288$ & $1.070$ & $4.347$ & $0.981$ \\
\bottomrule
\end{tabular}
\caption{Mean full-interface constraint-feasibility and mechanical-response diagnostics over the test set. All entries are percentages. The first row uses the exactly assembled load and exact elastic predictor $Y_\mu$ while the remaining accelerated variants use the POD--RBF predictor surrogate.}
\label{tab:contact_constraint_analysis}
\end{table}

The mean residual penetration is nearly independent of the approximation strategy and remains between $0.288\%$ and $0.321\%$. Neither the predictor surrogate, the NN augmentation, nor hyper-reduction therefore produces an appreciable change in average primal feasibility. In particular, replacing the exact elastic predictor by the POD--RBF surrogate changes the four tabulated ROM indicators by at most $0.005$ percentage points.

For the ROM and HR-ROM, the negative-normal-traction indicator is identical at $0.277\%$. Their main feasibility difference concerns the friction bound, for which the mean excess increases from $1.174\%$ for the ROM to $2.734\%$ for the HR-ROM. The NN-ROM and NN-HR-ROM exhibit larger negative-normal-traction indicators, $1.049\%$ and $1.070\%$, respectively, together with friction-bound excesses of $4.109\%$ and $4.347\%$. The small difference between the two NN-augmented models indicates that these larger traction-feasibility deviations are associated primarily with the NN augmentation, whereas hyper-reduction has only a limited additional effect in this family.

Despite these larger traction-feasibility indicators, the mechanical-response errors remain small for all five reduced variants, ranging from $0.885\%$ for the NN-ROM to $1.819\%$ for the HR-ROM. Replacing the exact predictor by its surrogate changes this indicator only from $1.057\%$ to $1.062\%$. Hyper-reduction then increases it from $1.062\%$ to $1.819\%$ for the projection-based surrogate models and from $0.885\%$ to $0.981\%$ for the NN-augmented models. Furthermore, larger local traction-feasibility deviations of the NN-augmented models are not accompanied by a larger model-resolved mechanical response. Within the primal directions represented by each model, the error in the combined normal and tangential contact force therefore produces only a limited displacement response relative to the HDM contact force.

\section{Selection of the HDM snapshot budget}
  \label{app:snapshot_budget}

  This appendix details the procedure used to select the number of HDM snapshots used to construct the POD bases. The same procedure is applied to both test cases. The training parameters are generated from an ordered Sobol sequence. For a given budget $N_{\mathrm{snap}}$, the training set consists of the first $N_{\mathrm{snap}}$ Sobol samples. The budgets are therefore nested.

  For each value of $N_{\mathrm{snap}}$, primal POD bases are constructed from the corresponding HDM snapshots. The approximation quality of the primal POD basis is then assessed directly on the validation set through the projection error of the reconstructed primal solution. The projection is evaluated for the POD dimensions
  \[ r \in \{50,100,150,200\}. \]
  For each validation parameter, the primal snapshot is projected onto the $r$-dimensional POD space, the corresponding primal solution is reconstructed, and its error with respect to the HDM solution is measured in the energy norm. 

    Denoting by $V_r^{(N_{\mathrm{snap}})}$ the first $r$ POD modes obtained from $N_{\mathrm{snap}}$ training snapshots, the reconstructed projection and its relative energy error are written, for the two-dimensional problem, as
    \[ U_{\mu,r}^{\mathrm{proj}}=V_r^{(N_{\mathrm{snap}})} \left(V_r^{(N_{\mathrm{snap}})}\right)^T U_\mu, \qquad e_{\mathrm{proj}}(\mu;N_{\mathrm{snap}},r)= \frac{\left\|U_{\mu,r}^{\mathrm{proj}}-U_\mu\right\|_K}{\left\|U_\mu\right\|_K}.\]
    For the three-dimensional contact problem, the same procedure is applied to the correction snapshots introduced in Section~4.2.3, yielding
    \[ U_{\mu,r}^{\mathrm{proj}}=Y_\mu+V_r^{(N_{\mathrm{snap}})} \left(V_r^{(N_{\mathrm{snap}})}\right)^T \left(U_\mu-Y_\mu\right).\]
    
  Let $E(N_{\mathrm{snap}},r)$ denote the mean relative primal projection energy error over the validation set for a basis constructed from $N_{\mathrm{snap}}$ snapshots and truncated to $r$ modes.

  For two consecutive budgets $N_{k-1}$ and $N_k$, we define the relative error decrease
  \[ \delta_k(r) = 100\, \frac{E(N_{k-1},r)-E(N_k,r)}      {E(N_{k-1},r)} . \]

  The snapshot budget is then selected sequentially. Starting from $250$ snapshots, the budget is increased until the improvement obtained by adding new snapshots becomes small for all considered POD dimensions. More precisely, the first budget $N_k$ satisfying
  \[ \max_r \delta_k(r)\leq\varepsilon_{\max}, \qquad E(N_k,r)\leq E(N_{k-1},r) \quad \forall r, \]
  is retained. If the latter monotonicity condition fails for any $r$, the preceding budget $N_{k-1}$ is retained instead. In the numerical experiments, we use $\varepsilon_{\max}=5\%$.

\begin{figure}[pos=H]
    \centering
    \includegraphics[width=0.6\textwidth]{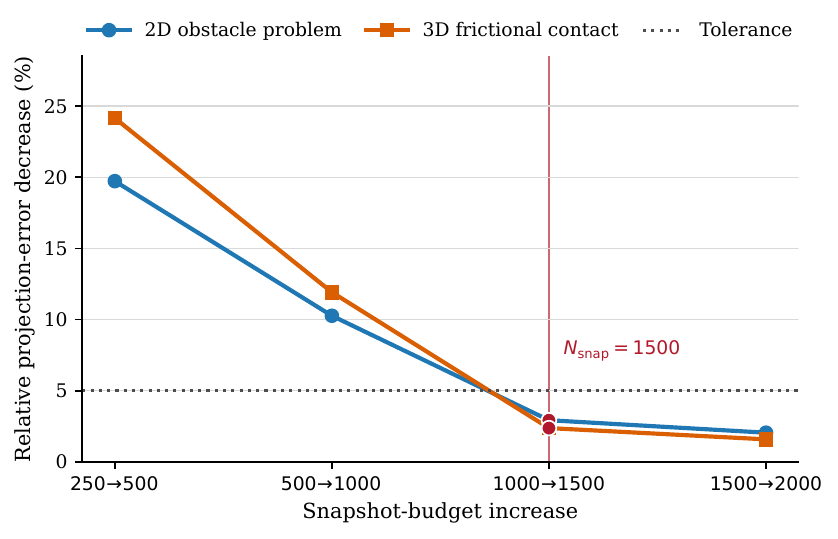}
    \caption{Maximum relative decrease of the primal validation projection energy error over the considered POD dimensions for the two test cases when increasing the snapshot budget. The two curves exhibit a similar saturation behavior and both happen to fall below the prescribed $5\%$ tolerance for the transition from $1000$ to $1500$ snapshots.}
    \label{fig:snapshot_budget_saturation}
  \end{figure}

  Figure~\ref{fig:snapshot_budget_saturation} illustrates the resulting criterion for the two considered problems. The tested snapshot budgets are $N_{\mathrm{snap}}\in\{250,500,1000,1500,2000\}$. It shows, for each consecutive pair of budgets, the maximum relative decrease of the validation projection energy error over the considered POD dimensions. The two test cases happen to give very similar curves. For both, the criterion is first satisfied for the transition from $1000$ to $1500$ snapshots, so $N_{\mathrm{snap}}=1500$ is selected. The $2000$-snapshot budget is shown for visualization but is not used in the numerical experiments.

\section{Neural-network training procedure}
\label{app:nn_architecture_selection}

This appendix summarizes the numerical settings used to train the neural correction maps and the procedures used to select their architectures. The loss functions are those introduced in Section~\ref{sec:nn_offline}.
In both test cases, each input and output component is standardized using its mean and standard deviation computed on the training set. After training, the network weights corresponding to the lowest validation loss are retained.

All considered networks are trained in double precision using the full-batch L-BFGS optimizer with a zoom line search, a tolerance of $10^{-6}$, and at most $5000$ iterations. The weights are initialized with the Glorot-uniform scheme and the biases are initialized to zero. For the final architectures, each configuration was trained once with a fixed random seed, so variability across random initializations was not quantified.

\subsection{Two-dimensional obstacle problem}

The enlarged POD dimensions are fixed to $n_{\mathrm{tot}}=m_{\mathrm{tot}}=150$, while $n_r=m_r=25$ leading coordinates are retained as online variables. The primal and dual networks therefore receive the $25$ retained coordinates together with the seven parameter features
\[
\eta_\mu=\bigl(\widehat{\gamma},c_x,c_y,\cos\theta,\sin\theta,\alpha,\kappa\bigr),
\]
and predict the remaining $125$ POD coordinates. The low-fidelity dataset contains $5000$ Galerkin-ROM samples, of which $4500$ are used for fitting and $500$ for validation.

The architecture is selected through an exhaustive Optuna grid search. To limit the size of the search, the primal and dual networks are constrained to share the same hidden-layer architecture, although their weights are trained independently. We consider the activation functions SiLU and $\tanh$ and the width patterns reported in Table~\ref{tab:nn_architecture_search_space}.

\begin{table}[pos=H]
\centering
\begin{tabular}{ll}
\toprule
\textbf{Number of hidden layers} & \textbf{Hidden widths} \\
\midrule
$2$ & $(256,256)$, $(512,512)$, $(1024,512)$, $(512,1024)$ \\
$3$ & $(256,512,256)$, $(512,512,512)$, $(1024,512,256)$, $(1024,512,512)$ \\
$4$ & $(256,512,512,256)$, $(512,512,256,256)$, $(1024,512,512,256)$ \\
\bottomrule
\end{tabular}
\caption{Architecture search space used for the two-dimensional NN-ROM correction maps.}
\label{tab:nn_architecture_search_space}
\end{table}
The candidate architectures are ranked according to
\[
J_{\mathrm{search}}
=
\frac{1}{2}
\left(
\widehat J_{\mathrm{val}}^u
+
\widehat J_{\mathrm{val}}^\lambda
\right),
\]
where $J_{\mathrm{val}}^u$ and $J_{\mathrm{val}}^\lambda$ are the primal and dual losses of Eqs.~\eqref{eq:primal_nn_loss}--\eqref{eq:dual_nn_loss}, evaluated on the validation set. Since these two losses have different physical scales, each component is normalized for comparability as
\[
\widehat J_{\mathrm{val}}^u
=
\frac{J_{\mathrm{val}}^u}{E_{u,\mathrm{val}}^c},
\qquad
\widehat J_{\mathrm{val}}^\lambda
=
\frac{J_{\mathrm{val}}^\lambda}{E_{\lambda,\mathrm{val}}^c},
\]
with
\[
E_{u,\mathrm{val}}^c
=
\frac{1}{N_{\mathrm{val}}}
\sum_{i=1}^{N_{\mathrm{val}}}
\left\|V^c q_{\mu^i}^c\right\|_K^2,
\qquad
E_{\lambda,\mathrm{val}}^c
=
\frac{1}{N_{\mathrm{val}}}
\sum_{i=1}^{N_{\mathrm{val}}}
\left\|B_{\mu^i}^T W^c\xi_{\mu^i}^c\right\|_{K^{-1}}^2.
\]
The best candidates according to this criterion are then evaluated at the NN-ROM level by solving the full reduced SSN system on the validation set. This second step is used because a small coordinate-prediction loss does not necessarily imply the best performance once the neural correction is embedded in the nonlinear solve.

Table~\ref{tab:nn_top_architectures} reports the five architectures retained after this solver-level comparison. Their validation errors are close, and all five use SiLU activations, suggesting that SiLU performs better than $\tanh$ for this particular test case. We therefore select the $(256,512,256)$ architecture as a compact compromise between accuracy and model size. This architecture is used for both the primal and dual correction maps.

\begin{table}[H]
\centering
\begin{tabular}{llcc}
\toprule
\textbf{Activation}
& \textbf{Hidden widths}
& \multicolumn{1}{c}{\shortstack{\textbf{Number of parameters}\\\textbf{per network}}}
& \multicolumn{1}{c}{\shortstack{\textbf{Mean validation}\\\textbf{energy error}}} \\
\midrule
SiLU & $(512,512,512)$       & $606\,333$ & $2.69\%$ \\
SiLU & $(1024,512,512)$      & $885\,373$ & $2.73\%$ \\
SiLU & $(256,512,512,256)$   & $566\,141$ & $2.75\%$ \\
\textbf{SiLU} & $\mathbf{(256,512,256)}$ & $\mathbf{303\,485}$ & $\mathbf{2.72\%}$ \\
SiLU & $(512,512)$           & $343\,677$ & $2.73\%$ \\
\bottomrule
\end{tabular}
\caption{Best NN-ROM candidate architectures after validation at the reduced-solver level. The parameter count is reported per network, and the selected architecture is highlighted in bold.}
\label{tab:nn_top_architectures}
\end{table}

\subsection{Three-dimensional frictional-contact problem}

For the three-dimensional problem, the enlarged dimensions are $n_{\mathrm{tot}}=150$, $m_n=50$, and $m_{\tau,\mathrm{tot}}=150$. The neural correction is applied only to the primal and tangential-dual coordinates. We retain $25$ leading coordinates for each map and learn the remaining $125$, while the normal-dual coefficients remain explicit online variables. With
\[
\eta_\mu=(Y_{CL},Z_{CL},Y_{CR},Z_{CR},F_\mu),
\]
both neural maps therefore have input dimension $30$ and output dimension $125$.

The low-fidelity dataset contains $7500$ successful Galerkin-ROM solutions generated from an independent scrambled Sobol design. The first $7000$ samples are used to fit the network weights, while the last $500$ are reserved for validation and checkpoint selection.

For the three-dimensional test case, the network architectures were chosen empirically from preliminary validation runs, allowing different capacities for the primal and tangential-dual maps. The selected primal network uses two hidden layers of widths $(512,512)$, while the tangential-dual network uses widths $(64,64)$. Both networks use SiLU activations.

\section{Sensitivity to the number of retained NN input coordinates in the 2D obstacle problem}
\label{app:nn_input_dimension_sensitivity}

We examine the influence of the number of retained reduced coordinates used as inputs to the neural networks in the case of the $2D$ obstacle problem.
To assess the influence of the retained reduced dimension, we fix the network architecture to the one selected in Appendix~\ref{app:nn_architecture_selection}, namely a three-hidden-layer SiLU network of size $(256, 512, 256)$.

For each value
\[
n_r=m_r \in \{5,10,15,25,35,45\},
\]
the enlarged POD dimensions are kept fixed to $n_{\mathrm{tot}}=m_{\mathrm{tot}}=150$. The primal and dual networks receive the retained reduced coordinates together with the same parameter-feature vector $\eta_\mu$ so that the actual input dimension is $n_r+7$ for the primal network and $m_r+7$ for the dual network. All networks are trained with the same low-fidelity training set. 

Because the networks are retrained independently for each retained dimension, the observed error differences mainly reflect the effect of the retained dimension, but some variability is also introduced by the NN training.

\begin{figure}[pos=H]
    \centering
    \includegraphics[width=0.6\linewidth]{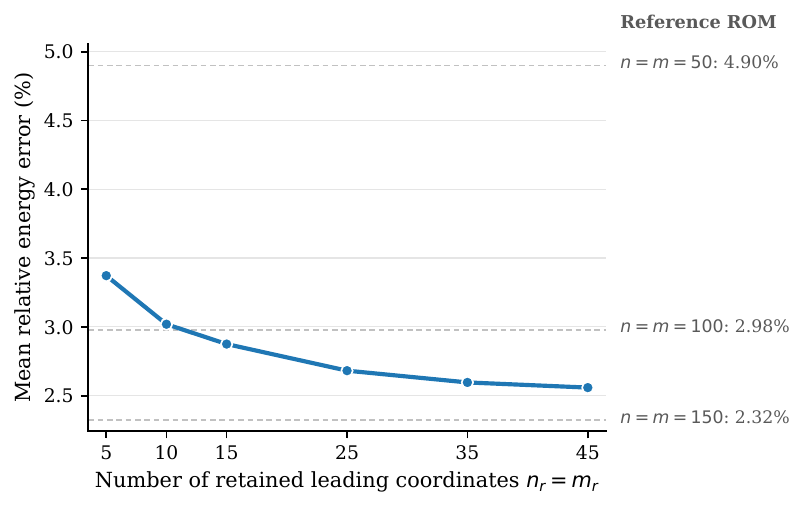}
    \caption{Sensitivity of the NN-ROM validation error with respect to the number of retained leading reduced coordinates used as neural-network inputs. The network architecture is fixed to three hidden layers of widths $256$, $512$, and $256$, with SiLU activations, and the same parameter features are used in all cases. The gray dashed lines report the validation errors of the standard Galerkin ROM for increasing dimensions $n=m$, and are shown as reference levels.}
    \label{fig:nn_input_dimension_sensitivity}
\end{figure}

Each trained pair of networks is then embedded in the NN-ROM solver and evaluated on the validation set. Figure~\ref{fig:nn_input_dimension_sensitivity} reports the mean relative energy error as a function of the number of retained leading coordinates.

The error decreases rapidly when increasing the number of retained coordinates from $5$ to about $25$, and then reaches a plateau. Already with $n_r=m_r=10$, the NN-ROM reaches an accuracy comparable to that of a standard Galerkin ROM with $n=m=100$ modes, highlighting the gain brought by the nonlinear correction. The choice $n_r=m_r=25$, used in the numerical experiments, therefore provides a good compromise between accuracy and online reduced dimension.

\section{Cubature convergence}
\label{app:cubature_convergence}
This appendix reports the convergence of the greedy cubature selection for both test cases. For each hyper-reduced contribution, we monitor the relative cubature residual
\[
\eta_\ell
=
\frac{\|\mathbf y-\mathbf G_{\mathcal S_\ell}\omega_\ell\|_2}
     {\|\mathbf y\|_2},
\]
as a function of the number of selected cubature points or interface faces $|\mathcal S_\ell|$.

For both test cases, the cubature rules are constructed from a cubature training set of $N_{\mathrm{cub}}=250$ converged reduced-model states. The corresponding parameter values are taken in order from the same scrambled Sobol training design used for the HDM snapshots.

For the two-dimensional obstacle problem, cubature is applied to the cubic nonlinear term and to the complementarity projection. The corresponding tolerances are $ \varepsilon_{\mathrm{cub}}^{\mathrm{cubic}}=10^{-2}$ and $\varepsilon_{\mathrm{cub}}^{\mathrm{comp}}=2\times10^{-2}$, and the maximum number of selected points is $c_{\max}=1250$.

For the three-dimensional frictional-contact problem, only the tangential Coulomb contact contribution is hyper-reduced. The cubature tolerance is set to $\varepsilon_{\mathrm{cub}}^{\mathrm{contact}}=2\times10^{-2}$, with at most $c_{\max}=750$ selected interface faces.

\begin{figure}[pos=H]
  \centering
  \includegraphics[width=0.85\textwidth]{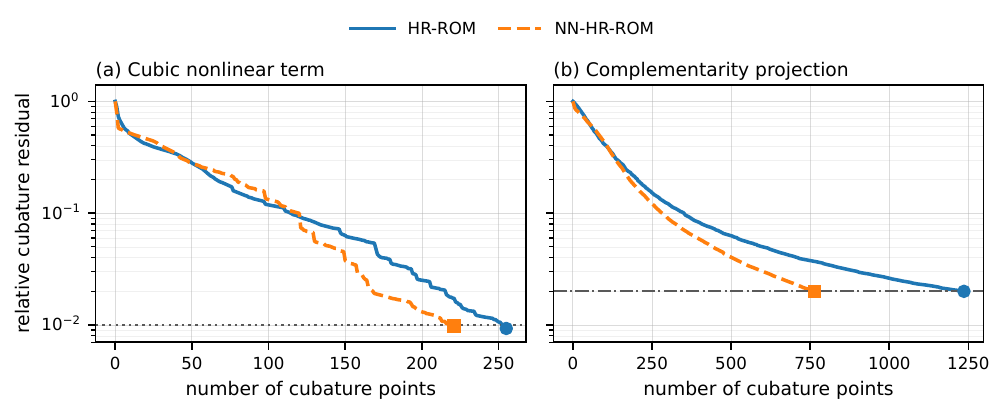}
\caption{Convergence of the greedy cubature residual for the two-dimensional obstacle problem. The curves correspond to the cubic nonlinear contribution and the complementarity projection for the HR-ROM and NN-HR-ROM cubature rules. The corresponding tolerances are $\varepsilon_{\mathrm{cub}}^{\mathrm{cubic}}=10^{-2}$ and $\varepsilon_{\mathrm{cub}}^{\mathrm{comp}}=2\times10^{-2}$. The final numbers of selected cubature points are $255/1236$ for the HR-ROM and $221/763$ for the NN-HR-ROM, for the cubic and complementarity contributions, respectively.}
\label{fig:cubature_convergence_2d}
\end{figure}

\begin{figure}[pos=H]
  \centering
  \includegraphics[width=0.48\textwidth]{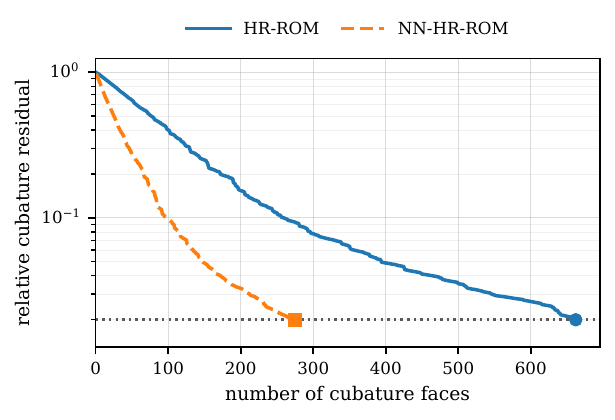}
\caption{Convergence of the greedy cubature residual for the tangential Coulomb contact contribution in the three-dimensional frictional-contact problem, for the HR-ROM and NN-HR-ROM cubature rules. The cubature tolerance is $\varepsilon_{\mathrm{cub}}^{\mathrm{contact}}=2\times10^{-2}$. The final numbers of selected interface faces are $662$ for the HR-ROM and $275$ for the NN-HR-ROM.}
\label{fig:cubature_convergence_3d}
\end{figure}

Figures~\ref{fig:cubature_convergence_2d} and~\ref{fig:cubature_convergence_3d} show the corresponding greedy-cubature convergence curves for the two test cases.

\section{Shifted reduced system for the three-dimensional contact problem}
\label{app:contact_shifted_rom}

This appendix gives the reduced system associated with the correction-based representation introduced in Section~\ref{sec:contact_correction_pod}. Let $Y_\mu$ be the unconstrained elastic predictor,
\[K Y_\mu=L_\mu,\]
and write
\[ U_\mu=Y_\mu+\delta U_\mu . \]
The load contribution is then eliminated from the equilibrium equation, while the contact laws remain imposed on the total displacement. We define the known interface quantities
\[ g^Y_{\mu,n}=B_nY_\mu, \qquad g^Y_{\mu,\tau}=B_\tau Y_\mu . \]
With the reduced approximations
\[ \delta U_\mu\approx Vq_\mu,\qquad \Lambda_{\mu,n}\approx W_n\xi_{\mu,n},\qquad \Lambda_{\mu,\tau}\approx W_\tau\xi_{\mu,\tau}, \]
the shifted Galerkin ROM reads: find $(q_\mu,\xi_{\mu,n},\xi_{\mu,\tau}) \in \mathbb R^n\times\mathbb R^{m_n}\times\mathbb R^{m_\tau} $ such that
\begin{equation}
\label{eq:contact_shifted_rom}
  \begin{aligned} 
  V^T\left( KVq_\mu +B_n^T W_n\xi_{\mu,n} +B_\tau^T W_\tau\xi_{\mu,\tau} \right)&=0,\\ 
  W_n^T\left[ W_n\xi_{\mu,n} - \left[ W_n\xi_{\mu,n} +\rho\left(g^Y_{\mu,n}+B_nVq_\mu\right) \right]_+ \right]&=0,\\ 
  W_\tau^T\left[ W_\tau\xi_{\mu,\tau} - \Pi_{c_\mu(W_n\xi_{\mu,n})} \left( W_\tau\xi_{\mu,\tau} +\rho\left(g^Y_{\mu,\tau}+B_\tau Vq_\mu\right) \right) \right]&=0. 
\end{aligned}
\end{equation}
Here $\Pi_{c_\mu(\cdot)}$ is the face-wise Coulomb projection defined in Section~\ref{sec:3D_contact_problem}. 

When the elastic predictor is approximated by the POD--RBF surrogate, the contact equations are approximated by replacing the interface traces with their surrogate values:
\[ g^Y_{\mu,n} \approx \widehat g^Y_{\mu,n}=B_n\widehat Y_\mu,\qquad g^Y_{\mu,\tau}\approx \widehat g^Y_{\mu,\tau}=B_\tau\widehat Y_\mu . \]

Note that the shifted equilibrium equation~\eqref{eq:contact_shifted_rom} is derived using the exact predictor $Y_\mu$, so that $KY_\mu=L_\mu$ cancels the load term exactly. The surrogate approximation is introduced only afterwards in the shifted contact terms, and the projected surrogate equilibrium defect $V^T(K\widehat Y_\mu-L_\mu)$ is therefore not included in the reduced equilibrium equation. Thus, the surrogate is used to approximate the shifted interface traces rather than as a direct substitute for $Y_\mu$ in the equilibrium equation.

To quantify the accuracy of this approximation and the magnitude of the omitted equilibrium defect, Table~\ref{tab:contact_predictor_surrogate_diagnostics} reports the corresponding mean relative errors over the validation set. For the projected equilibrium defect, $V$ denotes the $150$-dimensional primal correction basis used in the ROM. The interface traces used in the contact equations are reproduced with mean relative errors below $0.1\%$, while the omitted projected equilibrium defect is also below $0.1\%$ on average.

\begin{table}[pos=H]
\centering
\begin{tabular*}{0.7\textwidth}{@{\extracolsep{\fill}}lcc}
\toprule
\textbf{Indicator} & \textbf{Definition} & \textbf{Mean}\\
\midrule
\addlinespace[0.5ex]
Predictor energy error & $\displaystyle \frac{\|\widehat Y_\mu-Y_\mu\|_K}{\|Y_\mu\|_K}$ & $0.404 \%$\\
\addlinespace[1.5ex]
Normal-trace error & $\displaystyle \frac{\|\widehat g^Y_{\mu,n}-g^Y_{\mu,n}\|_\Gamma}{\|g^Y_{\mu,n}\|_\Gamma}$ & $0.078 \%$\\
\addlinespace[1.5ex]
Tangential-trace error & $\displaystyle \frac{\|\widehat g^Y_{\mu,\tau}-g^Y_{\mu,\tau}\|_\Gamma}{\|g^Y_{\mu,\tau}\|_\Gamma}$ & $0.090 \%$\\
\addlinespace[1.5ex]
Projected equilibrium defect & $\displaystyle \frac{\|V^T(K\widehat Y_\mu-L_\mu)\|_2}{\|V^TL_\mu\|_2}$ & $0.098 \%$\\
\addlinespace[0.5ex]
\bottomrule
\end{tabular*}
\caption{Mean accuracy diagnostics of the POD--RBF elastic-predictor surrogate over the validation set. }
\label{tab:contact_predictor_surrogate_diagnostics}
\end{table}

For the NN-ROM, the primal correction and tangential multiplier are reconstructed on the nonlinear manifolds introduced in Section~\ref{sec:nn_online}, while the normal multiplier retains its linear representation:
\[
\delta U_\mu\approx\mathcal U_\mu(q_\mu^r),\qquad
\Lambda_{\mu,n}\approx W_n\xi_{\mu,n},\qquad
\Lambda_{\mu,\tau}\approx\mathcal D_{\mu,\tau}(\xi_{\mu,\tau}^r).
\]
Denoting by $T_{\mathcal U_\mu}(q_\mu^r)$ and $T_{\mathcal D_{\mu,\tau}}(\xi_{\mu,\tau}^r)$ the corresponding reconstruction tangents, the shifted NN-ROM reads
\[
\begin{aligned}
T_{\mathcal U_\mu}(q_\mu^r)^T\left(K\mathcal U_\mu(q_\mu^r)+B_n^TW_n\xi_{\mu,n}+B_\tau^T\mathcal D_{\mu,\tau}(\xi_{\mu,\tau}^r)\right)&=0,\\
W_n^T\left[W_n\xi_{\mu,n}-\left[W_n\xi_{\mu,n}+\rho\left(\widehat g^Y_{\mu,n}+B_n\mathcal U_\mu(q_\mu^r)\right)\right]_+\right]&=0,\\
T_{\mathcal D_{\mu,\tau}}(\xi_{\mu,\tau}^r)^T\left[\mathcal D_{\mu,\tau}(\xi_{\mu,\tau}^r)-\Pi_{c_\mu(W_n\xi_{\mu,n})}\left(\mathcal D_{\mu,\tau}(\xi_{\mu,\tau}^r)+\rho\left(\widehat g^Y_{\mu,\tau}+B_\tau\mathcal U_\mu(q_\mu^r)\right)\right)\right]&=0.
\end{aligned}
\]
The tangent-projected system is linearized using the quasi-Newton approximation introduced in Section~\ref{sec:nn_online}, in which derivatives of the reconstruction tangents with respect to the reduced coordinates are omitted.


\bibliographystyle{cas-model2-names}

\bibliography{cas-refs}



\end{document}